\documentclass{sn-jnl} 

\usepackage{anyfontsize}
\usepackage{times,amsmath,amssymb}
\usepackage{hyperref}
\usepackage{amsthm}
\usepackage{bbm}
\usepackage{bm}
\usepackage[OT2,T1]{fontenc}

\newtheorem{theorem}{Theorem}
\newtheorem{ Proposition }{Proposition}
\newtheorem{Corollary}{Corollary}
\newtheorem{Lemma}{Lemma}
\newtheorem{Definition}{Definition}

\newcommand\norm[1]{\left\lVert#1\right\rVert}
\newcommand\bb[1]{\mathbf{#1}}

\newcommand{\R}{\mathbb{R}}

\newcommand{\C}{\mathcal{C}}
\newcommand{\Zz}{\mathcal{Z}}
\newcommand{\Yy}{\mathcal{Y}}
\newcommand{\Ww}{\mathcal{W}}
\newcommand{\M}{\mathcal{M}}

\newcommand{\D}{\mathrm{D}}
\newcommand{\Dd}{\mathcal{D}}
\newcommand{\Ls}{\bm{\mathrm{L}}}
\newcommand{\ls}{\mathrm{L}}

\newcommand{\N}{\mathrm{N}}

\begin{document}

\title{Adaptive Vector-Valued Splines for the Resolution of Inverse Problems}

\author[1]{\fnm{Vincent} \sur{Guillemet}}\email{vincent.guillemet@epfl.ch}

\author*[1]{\fnm{Michael} \sur{Unser}}\email{michael.unser@epfl}

\affil[1]{\orgdiv{Biomedical Imaging Group}, \orgaddress{\street{Station 17}, \orgname{EPFL}, \city{Lausanne}, \postcode{CH-1015}, \country{Switzerland}}}

\abstract{
 We introduce a general framework for the reconstruction of vector-valued functions from finite and possibly noisy data, acquired through a known measurement operator. The reconstruction is done by the minimization of a loss functional formed as the sum of a convex data fidelity functional and a total-variation-based regularizer involving a suitable matrix $\Ls$ of differential operators. Here, the total variation is a norm on the space of vector measures. These are split into two categories: inner, and outer norms. The minimization is performed over an infinite-dimensional Banach search space. When the measurement operator is $\text{weak}^{\star}$-continuous over the search space, our main result is that the solution set of the loss functional is the closed convex hull of adaptive $\Ls$-splines, with fewer knots than the number of measurements. We reveal the effect of the total-variation norms on the structure of the solutions and show that inner norms yield sparser solutions. We also provide an explicit description of the class of admissible measurement operators.}

\vspace{1em}
\keywords{vector splines, total-variation norm,  representer
 theorems, optimization on measure spaces}

\pacs[MSC Classification]{41A15, 47A52, 46G10, 47N10, 93B28, 46N10}

\maketitle

\section{Introduction}
In this paper, we consider the resolution of vector-valued continuous-domain inverse problems (IP). Here, continuous-domain means that we search for a solution inside a continuum $\mathcal{X}$, an infinite-dimensional vector space. In its interpolation version, the task is to find a vector-valued function $\bb{f}=(f_{d})_{d=1}^D:\R\to\R^D$ such that 
\begin{equation}
\label{eq:1.0.2}
    \forall1\leq m\leq M,\quad
    y_m\approx\bb{c}_m^{\intercal}\bb{f}(t_m)=\sum_{d=1}^Dc_{m,d}f_d(t_m),\quad\text{without overfitting.}
\end{equation}
There, only the data $\bb{y}=(y_m)_{m=1}^M$ and the measurement parameters 
$(\bb{c}_m,t_m)_{m=1}^M=((c_{m,d})_{d=1}^D,t_m)_{m=1}^M\in\R^{(D+1)\times M}$ are known. The vector formulation used here is motivated by the fact that, in many applications, the components of $\bb{f}$ are either correlated or bound to a physical equation. Hence, they must be recovered jointly. For example, they could be the position $f_{1}$ and speed $f_{2}$, in a one-dimensional space (line) of a spring subject to an unknown external force \cite{egerstedt2009control}[page 27]. In this case, the (positional) measurement parameters in \eqref{eq:1.0.2} at time $t_m$ are given by $\bb{c}_m=(1,0)$. More generally, such IPs appear in system theory  \cite{egerstedt2009control}, time series \cite{durbin2012time}, and optimal control \cite{nagahara2020sparsity}, and $\bb{f}$ is recovered as a minimizer 

\begin{equation}
\label{eq:1.0.2.2}
    \bb{f}^{\star}\in\underset{\bb{f}\in\mathcal{X}}{\text{argmin }}\sum_{m=1}^M\left(y_m-\bb{c}_m^{\intercal}\bb{f}(t_m)\right)^2.
\end{equation}
The major concern is that the recovery of $\bb{f}$ from \eqref{eq:1.0.2} is an ill-posed problem. Indeed, there are infinitely many unknown parameters encoding $\bb{f}$ but only finitely many (possibly noisy) data points. One resolution strategy is to impose some regularity on the solutions. 

\subsection{Regularisation Strategies for the Scalar Interpolation Problem}
In this section, we consider the simpler case $D=1$ and drop some indices in the interest of conciseness. Then, regularization is done in a variational framework: a norm (or a seminorm) is added to the data fidelity functional in \eqref{eq:1.0.2.2} that is to be minimized. One possibility is to use the idea of the minimization of the energy of $\ls\{f\}$, where $\ls$ is an appropriate ordinary differential operator. This leads to Hilbertian-norm regularization with the optimization problem in
\begin{equation}
\label{eq:1.1.0}
    f^{\star}\in\underset{f\in\mathcal{X}}{\text{argmin }}\left(\sum_{m=1}^M\left(y_m-f(t_m)\right)^2+\lambda\int_{\R}\ls\{f\}(t)^2f\mathrm{d}t\right),
\end{equation}
where the search space $\mathcal{X}$ is an adequate continuous-domain Hilbert space. It is known that the unique solution $f^{\star}$ of this problem is a (scalar) $\ls^{\star}\ls$-spline of the form \cite{de1966splines,wahba1990spline}
\begin{equation}
    f^{\star}(\cdot)=q(\cdot)+\sum_{m=1}^M\alpha_mg_{\ls^{\star}\ls}(\cdot-t_m),\quad\quad q\in\mathcal{N}_{\ls}=\{f\in\Dd'(\R):\quad\ls\{f\}=0\},
\end{equation}
where $\ls^{\star}$ is the adjoint operator of $\ls$, $g_{\ls^{\star}\ls}$ is the causal Green's function of $\ls^{\star}\ls$, and $q$ is the signal component that is in the null space $\mathcal{N}_{\ls}$ of $\ls$, which is typically formed of polynomials. Observe that the nodes $(t_m)_{m=1}^M$ of the spline are exactly the measurement locations. 

Stemming in the infinite-dimensional compressed sensing of \cite{adcock2016generalized,adcock2018infinite}, 
another possibility is to favor solutions for which $\ls\{f\}$ is parsimonious. This parsimony is promoted by the total-variation norm $\norm{\cdot}_{\M}$ from measure theory \cite{rudin}[Chapter 6]. It is an extension of the $\mathcal{L}_1$ norm and the continuous-domain counterpart of the $\ell_1$ norm. This leads to a Banach-space-based theory of regularization with the associated optimization problem in
\begin{equation}
 \label{eq:1.1.7} 
    f^{\star}\in\underset{f\in\mathcal{X}}{\text{argmin }}\left(\sum_{m=1}^M\left(y_m-f(t_m)\right)^2+\lambda\norm{\ls\{f\}}_{\M}\right),
\end{equation}
where the search space $\mathcal{X}$ is an adequate continuous-domain Banach space. Then, the extreme points of the solution set in \eqref{eq:1.1.7} are known to be scalar adaptive $\ls$-splines \cite{fisher1975spline,unser2017splines} of the form
\begin{equation}
\label{eq:1.1.5}
  f^{\star}(\cdot)=q(\cdot)+\sum_{k=1}^K\alpha_k g_{\ls}(\cdot-x_k),\quad\quad K\leq M-N,\quad q\in\mathcal{N}_{\ls},
\end{equation}
where $(\alpha_{k})_{k=1}^{K}$ are the weights and  $(x_{k})_{k=1}^{K}$ are the knots of the spline. In contrast to the Hilbert-space theory, the spline knots $x_k$ may differ from the sampling locations $t_m$ and may be fewer as well. While going from the Hilbert to the Banach realm, we gained in solution adaptivity.

\subsection{Regularization Strategy for the Vector-Valued Interpolation Problem: A Hilbert-Space Theory}
In this paper, we consider the vector-valued counterpart of Section 1.1. Rooted in the problem of curve fitting \cite{weinert1978statistical}, the vector-valued interpolation problem is assumed to derive from the physical model, expressed in its state-space formulation \cite{egerstedt2009control, willems1997introduction} as
\begin{equation}
\label{eq:1.1.1}
    (\bb{f},\bb{x})\text{ s.t. }\begin{cases}
    \Ls\{\bb{f}\}(t)=\bb{x}(t),\forall t\in[0,T]\\
    y(t)=\bb{c}^{\intercal}(t)\bb{f}(t),
    \end{cases}
\end{equation}
where $\bb{x}:[0,T]\to\R^D$ is the control function to be determined. The operator $\Ls$ is a known square matrix of ordinary differential operator (ODO) that represents the dynamics of $\bb{f}$. The observation interval $[0,T]$, the measurement function $\bb{c}(\cdot):\R\to\R^D$,  and the measurements $(y(t))_{t\in[0,T]}$ are known. For the example of the spring previously mentioned, $\Ls$ encapsulates Hooke's and Newton's second laws. Often, only a discrete noisy subset of the data is available as $(y_m)_{m=1}^M=(\bb{c}_m^{\intercal}\bb{f}(t_m)+\epsilon_m)_{m=1}^M$, with $\bb{c}_m=\bb{c}(t_m)$. In this case, one minimizes the energy of the unknown control function $\bb{x}(\cdot)$. Thus, \eqref{eq:1.1.1} is replaced by the variational solution \cite{sidhu1984vector}
\begin{equation}
\label{eq:1.1.2} \bb{f}^{\star}=\underset{\bb{f}\in\mathcal{X}}{\text{argmin }} \left(\sum_{m=1}^M\left(y_m-\bb{c}_m^{\intercal}\bb{f}(t_m)\right)^2+\lambda\int_{0}^{T}\langle\Ls\{\bb{f}\}(t),\Ls\{\bb{f}\}(t)\rangle\mathrm{d}t\right),
\end{equation}
where the optimal control $\bb{x}^{\star}$ is obtained as $\Ls\{\bb{f}^{\star}\}=\bb{x}^{\star}$. In \eqref{eq:1.1.2}, $\lambda\in\R^+$ balances the strength of the regularization with respect to the data fidelity and $\mathcal{X}$ is an adequate continuous-domain Hilbert space. Under Assumption 3 (see \ref{ass:2}), the solution $\bb{f}^{\star}$ to \eqref{eq:1.1.2} is unique \cite{sidhu1984vector}[Theorem 2.2] and is called vector-valued \emph{smoothing spline}.
The spline is of the form \cite{sidhu1984vector}\footnote{Equation \eqref{eq:1.2.4} does not appear in \cite{sidhu1984vector} but it is the combination of Equations 4.23, 4.15 and 4.4.}
\begin{equation}
\label{eq:1.2.4}\bb{f}^{\star}(\cdot)=\bb{q}(\cdot)+\sum_{m=1}^M\alpha_m\bb{G}_{\Ls^{\star}\Ls}(\cdot-t_m)\bb{c}_m,\quad\quad\bb{q}(\cdot)\in\mathcal{N}_{\Ls}.
\end{equation}
The associated control function $\bb{x}^{\star}(\cdot)$ \eqref{eq:1.1.2} is $\bb{x}^{\star}(\cdot)=\sum_{m=1}^M\alpha_m\bb{G}_{\Ls^{\star}}(\cdot-t_m)\bb{c}_m,$ where  $\bb{G}_{\Ls^{\star}\Ls}$ is the Green's matrix of $\Ls^{\star}\Ls$ (see Section \ref{sec:3.2}), and $\mathcal{N}_{\Ls}$ is the finite-dimensional null space of $\Ls$. The spline knots are exactly the measurement locations 
$(t_m)_{m=1}^M$. To each $\Ls$ is associated a specific vector-spline \eqref{eq:1.2.4}. These vector-valued splines are generally more elaborate than the simple vectorial concatenation of one-dimensional splines, as they may exhibit combinations of polynomials, exponentials, and trigonometric splines \cite{zhang1997splines}.

Equation \eqref{eq:1.1.2} is the vector-valued counterpart of \eqref{eq:1.1.0} in the Hilbertian setting. 
To the best of our knowledge, there is no vector-valued counterpart of \eqref{eq:1.1.7} available in the literature for the Banach-space setting. Furthermore, we recall in Section 1.3 that the Banach-space formulation has a crucial advantage over the Hilbert-space-based one.

\subsection{Vector-Valued Generalized Interpolation: Splines as Universal Solutions}
\hspace{+0.5cm}
So far, we only discussed the interpolation case. However real-world IP measurements typically involve more complex operators such as the Radon transform or the Fourier transform \cite{bertero2021introduction}. For the optimization problem to be well-posed, with a guaranteed existence of solutions, the measurement operator $\bb{V}=(\bm{\nu}_m)_{m=1}^M:\mathcal{X}\to\R^M$ must satisfy an appropriate form of continuity. For the case $D=1$ in the Banach-space setting, it was shown in \cite{unser2021unifying} that adaptive $\ls$-splines of the form given by \eqref{eq:1.1.5} are solutions of regularized IP, as long as the measurement operator $\bm{\nu}$ is $\text{weak}^{\star}$-continuous. For this reason, scalar adaptive $\ls$-splines are said to be universal solutions of regularized IP \cite{unser2017splines}.  For  $D\in\mathbb{N}=\{1,2,3,...\}$ in the Hilbert-space setting, we can show that the solution takes the form 
\begin{equation}
\label{eq:1.2.2.8}\bb{q}(\cdot)+\sum_{m=1}^M\alpha_m\left(\bb{G}_{\Ls^{\star}\Ls}\ast\bm{\nu}_m\right)(\cdot),\quad\bb{q}\in\mathcal{N}_{\Ls},
\end{equation}
which, in contrast to \eqref{eq:1.2.4}, depends on $\bm{\nu}_m$. This solution can be hard to calculate when
$\bm{\nu}_m$ is nonlocal.
In addition to the two advantages discussed in the end of Section 1.1, regularization with a Banach-space norm has the convenience of yielding solutions of the form \eqref{eq:1.1.5} which are easy to calculate.

\subsection{Contributions}
In this paper, we define and study the $D$-dimensional generalization of \eqref{eq:1.1.7}:
\begin{equation}
\label{eq:1.4.11}
    \underset{\bb{f}\in\mathcal{X}}{\text{argmin }}\left(E(\bb{y},\bb{V}\{\bb{f}\})+\begin{cases}
        \lambda\norm{\Ls\{\bb{f}\}}_{\Zz',\M}\\
        \lambda\norm{\Ls\{\bb{f}\}}_{\M,\Zz'}
    \end{cases}\right),
\end{equation}
where $E(\cdot,\cdot)$ is the data fidelity functional, $\bb{V}:\mathcal{X}\to\R^M$ represents the measurement operator, and $\Ls\{\bb{f}\}\in\M(\R,\R^D)$ the Banach space of $\R^D$-valued measures (see Section \ref{sec:2}). The quantities $\norm{\cdot}_{\Zz',\M}$ and $\norm{\cdot}_{\M,\Zz'}$ are norms on $\M(\R,\R^D)$ called inner and outer norm, respectively, whose specifities are described in Section 2. The technical requirements on \eqref{eq:1.4.11} are listed in Section \ref{sec:4.1}.

Inspired by \eqref{eq:1.1.5}, \eqref{eq:1.2.4} we define an $\Ls$-spline $\bb{f}$ as an $\R^D$-valued distribution of the form 
\begin{equation}
\label{eq:1.3.10}
\bb{f}=\bb{G}_{\Ls}\ast\bb{m}+\bb{q},\quad\bb{q}\in\mathcal{N}_{\Ls},
\end{equation}
where $\bb{m}$ is an $\R^D$-valued atomic measure, $\Ls$ is a matrix of differential operators (MDO) with null space $\mathcal{N}_{\Ls}$ and $\mathrm{dim}(\mathcal{N}_{\Ls})=N$, determinant $\mathrm{det}(\Ls)$, and Green's matrix $\bb{G}_{\Ls}$ (see Section \ref{sec:3.2}).

We now summarize our main contributions, that fall into two categories. First, the contributions to the construction of the optimization problem.

\begin{itemize}

    \item \textbf{Framework for the Inversion of MDO.}
    Theorem \ref{prop:greens2} in Section \ref{sec:3.2} reveals that the Green's matrix $\bb{G}_{\Ls}$ of $\Ls$ exists as long as $\mathrm{det}(\Ls)\neq0$. This shows that $\Ls$-splines are well-defined. 

    \item \textbf{Identification of the Appropriate Topology.} Theorem \ref{prop:MLprop} in Section \ref{sec:3.3} states that the search space $\mathcal{X}=\M_{\Ls}(\R,\R^D)$ is a Banach space, where all $\Ls$-splines live, and reveals that $\norm{\Ls\{\cdot\}}_{\Zz',\M}$ and $\norm{\Ls\{\cdot\}}_{\M,\Zz'}$ are equivalent seminorms.

    \item \textbf{Characterization of Admissible Measurement Operators.} Theorem \ref{th:predual} in Section \ref{sec:3.4} reveals that  $\M_{\Ls}(\R,\R^D)$ is the dual space of $\C_{\Ls}(\R,\R^D)$ and Corollary \ref{coro:diracpredual} provides an easy-to-verify, sufficient condition for a compactly supported distribution to be in $\C_{\Ls}(\R,\R^D)$. This allows us to classify and describe which measurement operators $\bb{V}$ are $\text{weak}^{\star}$-continuous on $\M_{\Ls}(\R,\R^D)$.

\end{itemize}

Second, the contributions to the analysis of the optimization problem.
\begin{itemize}
    
    \item \textbf{Well-Posedness and Description of the Solutions of the Optimization Problem.}
    We show that $\Ls$-splines are universal solutions to vector-valued IP. Theorem \ref{th:optweakstar} reveals that, when $\bb{V}$ is $\text{weak}^{\star}$-continuous, the solution set of \eqref{eq:1.4.11} is the closed convex hull of its extreme points which all are adaptive $\Ls$-splines. In addition, we describe the effects of the norm $\norm{\cdot}_{\M(\R,\R^D)}$ on the $\Ls$-spline weights, knots, and sparsity. The most important point is that inner norms generate sparser $\Ls$-spline solutions, with fewer degrees of freedom.

\end{itemize}

\subsection{Related Work}
\begin{itemize}
    \item \textbf{System Theory.} The examination of (time invariant) linear systems  \cite{willems1997introduction, rugh1996linear} leads one to study the existence and structure of the solutions of linear systems of differential equations as
    \begin{equation}
    \label{eq:1.5.13}
        \Ls\{\bb{f}\}=\bb{Q}\{\bb{x}\},
    \end{equation}
    where $\bb{f}$, $\bb{x}$ are, a priori, vector-valued distributions and $\Ls$, $\bb{Q}$ are MDO of appropriate dimensions. Observe that \eqref{eq:1.5.13} is more general than \eqref{eq:1.1.1}. This structure is used, for example, for the building of monotone smoothing splines in \cite{nagahara2013monotone}. We reserve this general case for a later work. 
    
    \item \textbf{Control Theory-Control Splines.} Several problems from control theory
    are known to be closely related to the classic numerical theory of splines \cite{egerstedt2009control} . Some of these problems directly involve a regularized IP of the form \eqref{eq:1.1.2}. Nevertheless, to the best of our knowledge, the regularization is always done with a Hilbertian norm. The authors of \cite{nagahara20141} focus on $\mathcal{L}_1$-control smoothing splines but these problems are only formulated on discrete-domain Banach spaces. Some $\mathcal{L}_1$-control interpolating splines were also shown to be solutions of some different optimal control problem in \cite{nagahara2020sparsity}[Theorem 8.5].
    
    \item \textbf{Banach-Space Regularization.} Previous investigations of regularized continuous-domain IP on Banach spaces fall into two categories. The first one focuses on general continuous-domain IP \cite{flinth2019exact,unser2021unifying,bredies2020sparsity,boyer2019representer}, regularized by an abstract seminorm. Although these results typically provide an abstract description of the set of minimizers, they often require significant additional work to be applicable in any specific scenario.  The second line of work focuses on IP on specific Banach spaces, for example on the scalar signals defined over the Euclidean space $\R^d$\cite{unser2017splines,9655475}, the torus $\mathbb{T}^d$\cite{fageot2020tv}, or time-dependent signals \cite{bredies2023generalized}. In the scenario $\Ls=\bb{I}$, the author of \cite{fernandez2016super} considered the group total variation norm gTV as a continuous-domain generalization of the group sparsity norm \cite{tropp2006algorithms,yuan2006model}.  Here, gTV corresponds to the inner norm $\norm{\cdot}_{2,\M}$. The authors of \cite{boyer2019representer} subsequently described the extreme points of the solution  set of gTV-regularized IP, where they exhibited the sparsifying effect of this inner norm.
    
    \item \textbf{Reproducing-Kernel Hilbert Spaces (RKHS) and Reproducing Kernel Banach Spaces (RKBS).}
     Interpolating and smoothing splines for $\mathcal{L}_2$-regularized problems were introduced in \cite{sidhu1979vector,sidhu1984vector}. For the interpolation problem, a more modern take on the subject would leverage the RKHS-theory \cite{scholkopf2001generalized,alvarez2012kernels}, which is by nature limited to Hilbertian norms. Some Banach counterpart of these methods (RKBS) also exists (see \cite{zhang1997splines}, more recently \cite{bartolucci2023understanding} and,  for vector-valued neural networks \cite{shenouda2024variation}). Nevertheless, these results are often less general than the one developed for regularization on abstract Banach spaces and are limited to interpolation and smoothing problems.

\end{itemize}

\subsection{Notations}
\begin{itemize}
    \item [] \emph{The Letter $d$} The capital letters $D,D'$ are dimensions. The straight letter $\D$ is the derivative operator. The set $\mathcal{D}(\R,\R^{D\times D'})$ is the space of $\R^{D\times D'}$-valued, compactly supported, and $\mathcal{C}^{\infty}$ functions on $\R$. The set  $\mathcal{D}'(\R,\R^{D\times D'})$ is the space of $\R^{D\times D'}$-valued distributions. The small letter $d$ will always be an iterable.
    \vspace{0.2cm}
    \item [] \emph{Calligraphic Letters} Calligraphic letters are used for vector spaces or, occasionally, sets. For example, the set of minimizers $\mathcal{V}$, the open set $O$.
    \vspace{0.2cm}
    \item [] \emph{Operators} The action of an operator on an element is denoted by the bracket notation $\{\cdot\}$. The composition of operators is denoted by $\circ$. The convolution between two functions or distributions is denoted by $\ast$. For shift-invariant operators, we may switch from the notation $\{\cdot\}$ or $\circ$ to $\ast$.  Finally, the pointwise product between the components of two vectors is denoted by $\odot.$
    \vspace{0.2cm}
    \item [] \emph{Font of Letters} Capital and bold letters such as $\bb{G}$ refer to a matrix. The entries of the matrix are denoted by capital, straight letters such as $\mathrm{G}_{d,d'}$, for the entry $(d,d')$. ODO are also denoted by capital, straight letters. Vectors are denoted by lowercase, bold letters such as $\bm{\phi}$. The components of a vector are denoted by lowercase italic such as $\phi_d$, for the entry $(d)$. Finally, uppercase italic letters such as $D$ or $N$ are dimensions. 
    \vspace{0.2cm}
    \item [] \emph{Operations on Operators} The adjoint operation is denoted by $^{\star}$. The transpose operation is denoted by $^{\intercal}$. The time-reversal operation $\cdot\mapsto(-\cdot)$ is denoted by $^{\vee}$.
    \vspace{0.2cm}
    \item [] \emph{Duality Bracket} The duality bracket is denoted by $\langle\cdot,\cdot\rangle$. We do not specify on which space the operation happens as it is often clear from the context. Likewise, dependent on the context, $\langle\cdot,\cdot\rangle$ may denote an inner product.
\end{itemize}

\section{Vector-Valued Measures}
\label{sec:2}

 For the sake of flexibility in the specification of the norm on the space of vector-valued measures, we introduce the finite-dimensional vector space $\Zz=(\R^D,\norm{\cdot}_\Zz)$ and its dual $\Zz'=(\R^D,\norm{\cdot}_{\Zz'})$. 
 
 There exist two philosophies for the norming of the space of vector measures.
\begin{itemize}
    \item [(i)] One first integrates the measures $\bm{\mu}=(\mu_d)_{d=1}^D$ (one takes its component-wise total variation) to find the $\R^D$-vector $\int_{\R}\vert\bm{\mu}\vert=\left(\int_{\R}\vert\mu_d\vert\right)_{d=1}^D$ and then norms this vector to get the (informal) norm $\norm{\int_{\R}\vert\bm{\mu}\vert}_{\Zz'}$. These are called \emph{outer} norms because the finite-dimensional norm $\norm{\cdot}_{\Zz'}$ is taken \emph{out} of the integration process.
    \item [(ii)] One first takes "pointwise" the norm $\norm{\cdot}_{\Zz'}$ of a measure $\bm{\mu}$ to reduce it to the scalar-valued measure $\norm{\bm{\mu}}_{\Zz'}(\cdot)$, One then evaluates its total variation to get the (informal) norm $\int_{\R}\norm{\bm{\mu}}_{\Zz'}$. These are called \emph{inner} norms because the finite-dimensional norm $\norm{\cdot}_{\Zz'}$ is taken \emph{inside} the integrate.
\end{itemize}

 \subsection{Vector-Valued Measures with Outer Norms}
 \label{sec:2.1}
We use $\M(\R)$ to denote the space of finite Radon measures equipped with the total-variation norm
\begin{equation}
\label{eq:2.1.1.11}
    \norm{\mu}_{\M}:=\underset{\pi}{\text{sup}}\sum_{E\in\pi}\vert{\mu}(E)\vert<\infty,
\end{equation}
where the supremum is taken over all finite partitions $\pi$ of $\R$. By the Riesz-Markov-Kakutani theorem \cite{rudin}[Theorem 6.19], it is the functional dual of the
space $\C_0(\R)$ of continuous vanishing functions endowed with the sup norm $\norm{\cdot}_{\infty}$. Thus, the total-variation norm on $\M(\R)$, for
which it forms a Banach space \cite{rudin}[Chapter 6], is given by duality as
\begin{equation}
\label{eq:2.1.1.12}
\norm{\mu}_{\M}=\underset{g\in\C_0(\R),\norm{g}_{\infty}\leq1}{\text{sup}}\vert\langle \mu,g\rangle\vert.
\end{equation}
The norm $\norm{\cdot}_{\M}$ is a generalization of the $\norm{\cdot}_{\mathcal{L}_1}$ norm. Indeed if $\mu\in\mathcal{L}_1(\R)$, then $\norm{\mu}_{\mathcal{L}_1}=\norm{\mu}_{\M}$. The vector space $\left(\M(\R,\R^D),\norm{\cdot}_{\M,\Zz'}\right)$ is the Banach space of (countably additive) $\R^D$-valued measures $\bm{\mu}$ such that 
\begin{equation}
\label{eq:2.1.1.13}
    \norm{\bm{\mu}}_{\M,\Zz'}:=\norm{\norm{\bm{\mu}}_\M}_{\Zz'}<\infty,\quad\quad\mathrm{where}\quad\norm{\bm{\mu}}_{\M}:=\left(\norm{[\bm{\mu}]_d}_{\M}\right)_{d=1}^D.
\end{equation}
Here, the total-variation $\norm{\cdot}_{\M}$ is taken component-wise on $\bm{\mu}$. It is also a vectorial extension of the classic $\M(\R)$ because there exists $\lambda>0$ such that $\norm{\cdot}_{\M,\Zz'}=\lambda\norm{\cdot}_{\M}$ when $D=1$. Moreover, if $\bm{\mu}$ is an absolutely continuous measure such that $\forall1\leq d\leq D:[\bm{\mu}]_d\in\mathcal{L}_1(\R)$, then 
    \begin{equation}
    \norm{\bm{\mu}}_{\M,\Zz'}=\norm{\norm{\bm{\mu}(\cdot)}_{\mathcal{L}_1(\R)}}_{\Zz'}=\norm{\left(\int_{\R}\vert[\bm{\mu}]_d(x)\vert\mathrm{d}x\right)_{d=1}^D}_{\Zz'},
    \end{equation}
where we are identifying $\bm{\mu}$ to its Radon-Nikodym derivative with a slight abuse of notation. In addition, we observe that the definition of $\M(\R,\R^D)$ by \eqref{eq:2.1.1.13} is independant of the choice of outer norm $\norm{\cdot}_{\M,\Zz'}$. That implies that $\M(\R,\R^D)=\M(\R)^D$ without the need to specify how $\M(\R,\R^D)$ was initially defined. In Proposition \ref{prop:1}, whose proof is in Appendix \ref{app:2.1}, we list four properties of $\left(\M(\R,\R^D),\norm{\cdot}_{\M,\Zz'}\right)$ based on the following concept of monotone norm.
\begin{Definition}
\label{def:1}
    A norm $\norm{\cdot}$ on $\R^D$ is said to be monotone if the following holds:
    \begin{equation}
        \forall\{\bb{w},\bb{v}\}\subset\R^D,\bb{w}=(w_d)_{d=1}^D,\bb{v}=(v_d)_{d=1}^D:\quad\forall1\leq d\leq D\quad \vert w_d\vert\leq \vert v_d\vert\Rightarrow\norm{\bb{w}}\leq\norm{\bb{v}}.
    \end{equation}  
\end{Definition}
The dual norm of a monotone norm is again monotone \cite{bauer1961absolute}. All classic, eventually re-weighted norms $\norm{\cdot}_p$ for $1\leq p\leq\infty$ are monotone.

\begin{ Proposition }
\label{prop:1}
Let $\norm{\cdot}_{\Zz}$ be a monotone norm. The Banach space $\left(\M(\R,\R^D),\norm{\cdot}_{\M,\Zz'}\right)$ has the following properties.
\begin{itemize}
    \item [1.] The Banach space $\left(\M(\R,\R^D),\norm{\cdot}_{\M,\Zz'}\right)$ is the functional dual of $\left(\C_0(\R,\R^D),\norm{\cdot}_{\infty,\Zz}\right),$ where $\C_0(\R,\R^D)$ is the space of continuous functions $\bb{g}:\R\to\R^D$ that vanish at $\infty$, with its norm being
    \begin{equation}
        \label{2.1.1.14}
        \forall\bb{g}\in\C_0(\R,\R^D):\norm{\bb{g}}_{\infty,\Zz}:=\norm{\norm{\bb{g}}_{\infty}}_{\Zz},\quad\quad\mathrm{where}\quad\norm{\bb{g}}_{\infty}:=\left(\norm{[\bb{g}]_d}_{\infty}\right)_{d=1}^D.
    \end{equation}
    \item [3.] All norms on $\M(\R,\R^D)$ of the type \eqref{eq:2.1.1.13} are equivalent. More precisely, let $\norm{\cdot}_A$ and $\norm{\cdot}_B$ be two norms on $\R^D$, and let $\norm{\cdot}_{A'}$ and $\norm{\cdot}_{B'}$ be their dual norms. Denote by $\norm{\cdot}_{\M,A'}$ and $\norm{\cdot}_{\M,B'}$ the corresponding norms on $\M(\R,\R^D)$ according to \eqref{eq:2.1.1.13}. Then, there exist constants $C_1,C_2>0$ such that 
    \begin{equation}
    \label{eq:2.1.15}
    \forall\bm{\mu}\in\M(\R,\R^D):\quad C_1\norm{\bm{\mu}}_{\M,A'}\leq\norm{\bm{\mu}}_{\M,B'}\leq C_2\norm{\bm{\mu}}_{\M,A'}.
    \end{equation}

        \item [4.] The extreme points $\bb{e}$ of the centered unit ball $\mathcal{B}_{\M(\R,\R^D)}$ are exactly of the form
    \begin{equation}
    \label{eq:2.1.14}
        \bb{e}=\bb{a}\odot\bm{\delta}_{\bb{x}},
    \end{equation}
    where $\bb{x}\in\R^D$, $\bb{a}\in\R^D$ is an extreme point of the centered unit ball in $\Zz'$, and $\odot$ denotes the component-wise multiplication. The measure $\bm{\delta}_{\bb{x}}$ is defined as $[\bm{\delta}_{\bb{x}}]_d=\delta(\cdot-[\bb{x}]_d)$. Moreover, if $\bb{x}_k\neq \bb{x}_{k'}$ (for all rows) for all $k\neq k'$, then  
    \begin{equation}
    \label{eq:2.1.16}
    \norm{\sum_{k=1}^{K}\bb{a}_k\odot\bm{\delta}_{\bb{x}_k}}_{\M,\Zz'}=\norm{\left(\sum_{k=1}^K\vert [\bb{a}_k]_d\vert\right)_{d=1}^D}_{\Zz'}.
    \end{equation}
\end{itemize}
\end{ Proposition }

\subsection{Vector-Valued Measures with Inner Norms}
\label{sec:2.2}
$\quad$  The vector space $\left(\M(\R,\R^D), \norm{\cdot}_{\Zz',\M}\right)$ is the Banach space of (countably additive) $\R^D$-valued measure $\bm{\mu}$ such that 
\begin{equation}
\label{eq:2.1.11}
    \norm{\bm{\mu}}_{\Zz',\M}:=\underset{\pi}{\text{sup}}\sum_{E\in\pi}\norm{\bm{\mu}(E)}_{\Zz'}<\infty,
\end{equation}
where the supremum is taken over all finite partitions of $\R$. It is a Banach space \cite{dobrakov1971representation}. In particular, when $D=1$, $\left(\M(\R,\R^D),\norm{\cdot}_{\Zz',\M}\right)$ is the space $\left(\M(\R),\norm{\cdot}_{\M}\right)$ of all finite Radon measures, equipped with the total-variation norm. Moreover, if $\bm{\mu}$ is an absolutely continuous measure such that $\forall1\leq d\leq D:[\bm{\mu}]_d\in\mathcal{L}_1(\R)$, then 
\begin{equation}
        \norm{\bm{\mu}}_{\Zz',\M}=\norm{\norm{\bm\mu(\cdot)}_{\Zz'}}_{\mathcal{L}_1(\R)}=\int_{\R}\norm{\bm{\mu}(x)}_{\Zz'}\mathrm{d}x.
\end{equation}
As for outer norms, we observe that the definition of $\M(\R,\R^D)$ by \eqref{eq:2.1.11} is independant of the choice of inner norm $\norm{\cdot}_{\Zz',\M}$. That implies that $\M(\R,\R^D)=\M(\R)^D$ without the need to specify how $\M(\R,\R^D)$ was initially defined.
The vector space $\left(\C_0(\R,\R^D),\norm{\cdot}_{\Zz,\infty}\right)$ of $\R^D$-valued continuous functions $\bb{g}=([\bb{g}]_d)_{d=1}^D$ that vanish at $\infty$, equipped with the norm
\begin{equation}
\label{eq:2.1.12}
\norm{\bb{g}}_{\Zz,\infty}:=\underset{x\in\R}{\text{sup}}\norm{\bb{g}(x)}_{\Zz},
\end{equation}
is a Banach space whose dual is well-identified, as revealed by Proposition \ref{prop:2},  a generalization of the famous Riesz-Markov-Kakutani theorem which is taken from \cite{dobrakov1971representation}[Theorem 2 and Remark page 26]

\begin{ Proposition }
\label{prop:2}
The space $\left(\M(\R,\R^D),\norm{\cdot}_{\Zz',\M}\right)$ is the functional dual of $\left(\C_0(\R,\R^D),\norm{\cdot}_{\Zz,\infty}\right)$.
\end{ Proposition }
Consequently, the norm $\norm{\cdot}_{\Zz',\M}$ can be re-written as a dual norm:
\begin{equation}
    \label{eq:2.1.13}
\norm{\bm{\mu}}_{\Zz',\M}=\underset{\bb{g}\in\C_0(\R,\R^D)\setminus\{0\}:\norm{\bb{g}}_{\Zz,\infty}\leq1}{\text{sup}}\vert\langle\bm{\mu},\bb{g}\rangle\vert,\quad\quad\mathrm{where}\quad\langle\bm{\mu},\bb{g}\rangle:=\sum_{d=1}^D\langle[\bm{\mu}]_d,[\bb{g}]_d\rangle.
\end{equation}
Thus, the norm $\norm{\cdot}_{\Zz',\M}$ can either be seen from the measure-theoretic point of view \eqref{eq:2.1.11} or from the duality point of view \eqref{eq:2.1.13}, depending on what is the most practical. In Proposition \ref{prop:3}, whose proof is in Appendix \ref{app:2.2}, we state four properties of the space $\M(\R,\R^D)$ and the measures $\bm{\mu}$ that live in. 
\begin{ Proposition }
\label{prop:3}
The Banach space $\left(\M(\R,\R^D),\norm{\cdot}_{\Zz',\M}\right)$ has the following properties.
\begin{itemize}
    \item [1.] The norms $\norm{\cdot}_{1,\M}$ and $\norm{\cdot}_{\M,1}$, defined in \eqref{eq:2.1.1.13} and \eqref{eq:2.1.11} for $\norm{\cdot}_{\Zz'}=\norm{\cdot}_1$, are the same: $\forall\bm{\mu}\in\M(\R,\R^D),\norm{\bm{\mu}}_{1,\M}=\norm{\bm{\mu}}_{\M, 1}$.
    \item [2.] All norms on $\M\left(\R,\R^D\right)$ of the type \eqref{eq:2.1.11} are equivalent. More precisely, let $\norm{\cdot}_A$ and $\norm{\cdot}_B$ be two norms on $\R^D$ and let $\norm{\cdot}_{A'}$ and $\norm{\cdot}_{B'}$ be their dual norms. Denote by $\norm{\cdot}_{A',\M}$ and $\norm{\cdot}_{B',\M}$ the corresponding norms on $\M(\R,\R^D)$ according to \eqref{eq:2.1.11}. Then, there exist constants $C_1,C_2>0$ such that 
    \begin{equation}
    \forall\bm{\mu}\in\M(\R,\R^D):\quad C_1\norm{\bm{\mu}}_{A',\M}\leq\norm{\bm{\mu}}_{B',\M}\leq C_2\norm{\bm{\mu}}_{A',\M}.
    \end{equation}
    
    \item [3.]  The extreme points $\bb{e}$ of the centered unit ball $\mathcal{B}_{\M(\R,\R^D)}$ are exactly of the form
    \begin{equation}
        \bb{e}=\bb{a}\delta_{x},
    \end{equation}
    where $x\in\R$ and $\bb{a}$ is an extreme point of the centered unit ball in $\Zz'$. In addition, for $\{\bb{a}_k\}_{k=1}^{K}$ and $\{x_k\}_{k=1}^{K}$ two lists of coefficients $\bb{a}_k\in\R^D$ and $x_k\in\R$, such that $x_k\neq x_{k'}$ for all $k\neq k'$, one has that
    \begin{equation}
    \norm{\sum_{k=1}^{K}\bb{a}_k\delta_{x_k}}_{\Zz',\M}=\sum_{k=1}^K\norm{\bb{a}_k}_{\Zz'}.
    \end{equation}
\end{itemize}
\end{ Proposition }

The combination of Item 2 of Proposition \ref{prop:1} and Items 1 and 2 of Proposition \ref{prop:3} implies that all inner and outer norms are topologically equivalent on $\M(\R,\R^D)$. The calculation of the norm of atomic measures, together with the characterization of extreme points (Item 3), will be particularly useful to characterize the solution set \eqref{eq:3.1.30} (see Theorem \ref{th:optweakstar}).

\section{The Space of \texorpdfstring{$\R^D$}{RD}-Valued Splines}
\label{sec:3}
\hspace{+0.5cm}In this section, we define a suitable notion of $\R^D$-valued splines, associated to a matrix $\Ls$ of differential operators. We further define a bounded-variations type of space, shaped by $\Ls,$ which will become our search space in Section \ref{sec:4}.

\subsection{Matrix Differential Operators and Convolutions}
\label{sec:3.1}

\hspace{+0.5cm}An \emph{ordinary differential operator} (ODO) $\mathrm{L}$ is an operator of the form 
\begin{equation}
    \ls=\sum_{n=0}^{N}a_n\D^n,\quad(a_n)_{n=0}^N\in\R^{N+1},\quad\text{with}\quad\D^0=\text{Identity},
\end{equation}
where $\D^n$ is the $n$th derivative operator and $N$ is the \emph{degree} of $\ls$. Although ODO are not injective on $\mathcal{D}'(\R)$, they all admit a Green's function whose existence is guaranteed by the Malgrange–Ehrenpreis theorem \cite{reed2003methods}. The causal Green's function is unique \cite{unser2005cardinal} and will be denoted by $\ls^{-1}\delta$ or $g_{\ls}$.

A \emph{matrix (ordinary) differential operator} (MDO) $\bb{L}:\mathcal{D}'(\R,\R^{D\times D'})\to\mathcal{D}'(\R,\R^{D\times D'})$ is an operator of the form  
\begin{equation}
\Ls=\begin{bmatrix}
\ls_{1,1} & \cdots & \ls_{1,D} \\
\vdots & \ls_{r,c} & \vdots \\
\ls_{D,1} & \cdots & \ls_{D,D} \\
\end{bmatrix},\quad\quad\text{with}\quad[\bb{L}]_{r,c}=\ls_{r,c}:\mathcal{D}'(\R)\to\mathcal{D}'(\R).
\end{equation}
Its action on a matrix-valued distribution $\bb{G}=(\mathrm{G}_{d,d'})_{d,d'=1,1}^{D,D'}\in\mathcal{D}'(\R,\R^{D\times D'})$ is defined as 
\begin{equation}
    \forall r\in[1\cdots D],\forall c\in[1\cdots D']:\quad [\Ls\{\bb{G}\}]_{r,c}=\sum_{d=1}^D\ls_{r,d}\{\mathrm{G}_{d,c}\},
\end{equation}
where each $\mathrm{G}_{d,c}\in\mathcal{D}'(\R)$ is a distribution and where each $\ls_{r,c}$ is an ODO whose degree is $N_{r,c}$. First note that $\Ls$ is linear and continuous because each differential operator $[\Ls]_{r,d}$ is linear and continuous on $\mathcal{D}'(\R)$. Moreover, the space of such matrix operators $\Ls$ is a ring for the addition of operators and the matrix multiplication, where the formal multiplication is replaced by the composition of operators. The usual product in the matrix multiplication is replaced by the composition of operators. The adjoint $\Ls^{\star}$ of $\Ls$ is such that $\forall r,c\in[1\cdots D],$ $[\Ls^{\star}]_{r,c}=\ls_{c,r}^{\star}$, where we recall that $(\D^n)^{\star}=(-1)^n\D^n.$

The \emph{convolution} between two matrix-valued test functions $\bm{\psi}\in\mathcal{D}(\R,\R^{D\times D'})$ and $\tilde{\bm{\psi}}\in\mathcal{D}(\R,\R^{D'\times D''})$ is the $\R^{D\times D''}$-valued test function given by 
\begin{equation}
    (\bm{\psi}\ast\tilde{\bm{\psi}})(y):=\int_{\R}\bm{\psi}(y-x)\cdot\tilde{\bm{\psi}}(x)\mathrm{d}x,\quad\quad\mathrm{where}\quad\bm{\psi}(y-x)\cdot\tilde{\bm{\psi}}(x)\in\R^{D\times D''}.
\end{equation}
This definition extends naturally to the convolution between a matrix-valued test function and a matrix-valued distribution. It also extends to the convolution of two matrix-valued distributions, if it is well-defined, which happens when the two distributions are compactly supported or one-sided. A brief summary of these matrix-valued convolutions is given in Appendix \ref{app:B}.

\subsection{MDO, Green's Matrix, and Null Space}
\label{sec:3.2}
In the sequel, we shall simplify the notation $\mathcal{D}'(\R,\R^{D\times D})$ to $\mathcal{D}'$ when the context is clear. The \emph{Green's matrix} $\bb{G}_{\Ls}(\cdot)$ of $\Ls$ is a distribution $\bb{G}_{\Ls}\in\mathcal{D}'$ such that 
\begin{equation}
\label{eq:3.2.19}
\Ls\{\bb{G}_{\Ls}\}=\bb{I}\delta,
\end{equation}
where $\bb{I}\in\R^{D\times D}$ is the identity matrix. It follows from \eqref{eq:3.2.19} that 
\begin{equation}
\forall \bb{t}\in\R^D:\quad \Ls\{\bb{G}_{\Ls}\ast\mathbf{diag}(\bm{\delta}_{\bb{t}})\}=\mathbf{diag}(\bm{\delta}_{\bb{t}}).
 \end{equation}
We show in Proposition \ref{prop:greens2}, using the Smith matrix form, that a Green's matrix exists whenever
\begin{equation}
\mathrm{det}(\Ls)=\sum_{\sigma\in S(D)}\mathrm{sgn}(\sigma)\ls_{1,\sigma(1)}\circ\cdots\circ\ls_{D,\sigma(D)}\neq0,
\end{equation}
where $S(D)$ is the group of all permutations of the set $\{1,...,D\}$. The operator $\mathrm{det}(\Ls):\mathcal{D}'(\R)\to\mathcal{D}'(\R)$ is an ODO and we say that $\Ls$ is invertible if $\mathrm{det}(\Ls)\neq0$. We also recall the notion of unimodular matrices.
\begin{Definition}
An MDO $\bb{V}:\mathcal{D}'\to\mathcal{D}'$ is said to be unimodular if there exists another MDO $\bb{V}^{-1}:\mathcal{D}'\to\mathcal{D}'$ such that $\bb{V}\circ\bb{V}^{-1}=\bb{V}^{-1}\circ\bb{V}=\bb{I}$.
\end{Definition}
Observe that if $\bb{V}$ is unimodular then so is $\bb{V}^{-1}$. In Theorem \ref{prop:greens2}, we leverage the property that, for any invertible $\Ls$, there exist two unimodular MDO $\bb{U},\bb{V}$ and a diagonal MDO $\bb{D}=\textbf{Diag}(\D_1,...,\D_D),$ where $\D_d$ are all ODO \cite{willems1997introduction}[Theorem 2.5.15]. The operators $\bb{U},\bb{D}$, and $\bb{V}$ are all of size $(D\times D)$, such that $\bb{L}=\bb{U}\circ\bb{D}\circ\bb{V}$ and  $\mathrm{det}(\bb{D})=\pm\mathrm{det}(\bb{L})$.
\begin{theorem}
\label{prop:greens2}
An invertible MDO $\Ls$ has a unique causal Green's matrix $\bb{G}_{\Ls}$, which can be written as $\bb{G}_{\Ls}=\bb{V}^{-1}\circ\bb{D}^{-1}\circ\bb{U}^{-1}\{\delta\}$, where $\bb{D}^{-1}$ is defined for all $\forall\bm{\phi}=(\phi_d)_{d=1}^D\in\mathcal{D}'(\R,\R^D)$ as
\begin{equation}
    \forall1\leq d\leq D,\quad[\bb{D}^{-1}\{\bm{\phi}\}]_d=g_{\D_d}\ast\phi_d,
\end{equation}
where $g_{\D_d}$ is the unique causal Green's function of $\D_d$ and $\bb{U}^{-1}$,$\bb{V}^{-1}$ are unimodular matrices. In addition, the adjoint $\bb{G}_{\Ls}^{\star}$ of the Green's matrix verifies that 
\begin{equation}
\label{eq:3.2.37}
\bb{G}_{\Ls}^{\star}(\cdot)=\bb{G}_{\Ls}^{\intercal,\vee}(\cdot)=\bb{G}_{\Ls}^{\intercal}(-\cdot),
\end{equation}
and it is the unique anti-causal Green's matrix of $\Ls^{\star}$, with
\begin{equation}
\label{eq:3.2.38}
    \Ls^{\star}\{\bb{G}_{\Ls}^{\star}\}=\bb{I}\delta.
\end{equation}
\end{theorem}

The assumption that $\Ls$ is invertible is necessary. To illustrate this requirement, we consider the MDO
\begin{equation}
    \Ls=\begin{bmatrix}
        \D & \D\\
        0 & 0
    \end{bmatrix},\quad\text{s.t.}\quad\mathrm{det}(\Ls)=0,
\end{equation}
for which there exists no Green's matrix. As an additional issue, the null space $\mathcal{N}_{\Ls}$ is now infinite-dimensional. One remedy would be to use a "pseudo Green's matrix" that features some zeros on the diagonal, together with a restricted null space. In Section \ref{sec:4.3}, we consider only the case of non-invertible MDOs with a specific structure as the study of general non-invertible MDOs is outside of the scope of this paper. The practicality of the  assumption that $\Ls$ is invertible is discussed in Section \ref{ex:firstorder} for first-order systems.

Observe that some of the entries of the Green's matrix may be singular distributions, even though the order of the determinant is high. For example, consider the following MDO and its Green's matrix
\begin{align}
    \Ls=\begin{bmatrix}
        \D^2 & \D^4\\
        0 & \D^2
    \end{bmatrix},\quad&\text{with}\quad
    \bb{U}\circ\bb{D}\circ\bb{V}=
    \begin{bmatrix}
        1+\D^2 & \D^2\\
    1 & 1
    \end{bmatrix}\circ
    \begin{bmatrix}
        \D^2 & 0\\
        0 & \D^2
    \end{bmatrix}\circ
    \begin{bmatrix}
        1 & 0\\
        -1 & 1
    \end{bmatrix},\\
    \bb{G}_{\Ls}(\cdot)=    \begin{bmatrix}
        (\cdot)_+ & -\delta(\cdot)\\
        0 & (\cdot)_+
    \end{bmatrix},\quad&\text{with}\quad
    \bb{V}^{-1}\circ\bb{D}^{-1}\circ\bb{U}^{-1}\{\delta(\cdot)\}=
    \begin{bmatrix}
        1 & 0\\
    1 & 1
    \end{bmatrix}\circ
    \bb{D}^{-1}\circ
    \begin{bmatrix}
        1 & -\D^2\\
       -1 & 1+\D^2
    \end{bmatrix}\{\delta(\cdot)\}.
\end{align}
One directly finds that $\mathrm{det}(\Ls)=\D^4$. This shows that the degree of the determinant is not enough to infer the regularity of the entry $[\bb{G}_{\Ls}]_{r,c}$. Proposition \ref{prop:reg}, whose proof is in Appendix \ref{app:3.2}, provides an a priori (loose) lower-bound on the regularity of each entry $[\bb{G}]_{r,c}$. (We recall that $N_{c,r}$ is the order of the differential operator $[\Ls]_{c,r}$.)
\begin{ Proposition }
\label{prop:reg}
Let $\Ls$ be an invertible MDO such that $\forall d\in[1\cdots D],$
\begin{align}
    &\forall d'\in[1\cdots D]\setminus \{d\},\quad N_{d,d'}<N_{d,d},\quad\quad\text{or}\\
    &\forall d'\in[1\cdots D]\setminus \{d\},\quad N_{d',d}<N_{d,d}.
\end{align}
Then,
\begin{align}
   &\mathrm{if}\quad N_{c,r}\geq2\quad\mathrm{then} \quad[\bb{G}_{\Ls}]_{r,c}\in\mathcal{C}^{N_{c,r}-2}(\R);\\
    &\mathrm{if}\quad N_{c,r}=1\quad\mathrm{then} \quad[\bb{G}_{\Ls}]_{r,c}\in\mathcal{L}_{\infty,loc}(\R),
\end{align}
where $\mathcal{C}^{N_{c,r}-2}(\R)$ is the space of $(N_{c,r}-2)$-times continuously differentiable functions and $\mathcal{L}^{\infty}_{loc}$ is the space of locally, essentially bounded functions.
\end{ Proposition }
Let $\mathcal{N}_{\Ls}$ be the null space of $\Ls$ such that
\begin{equation}
    \mathcal{N}_{\Ls}=\{\bb{q}\in\mathcal{D}'(\R,\R^D):\quad\Ls\{\bb{q}\}=\bm{0}\in\R^D\}.
\end{equation}
It follows from \eqref{eq:3.2.19} that $\bb{G}$ is a Green's matrix of $\Ls$ if and only if it is of the form 
\begin{equation}
\bb{G}=\bb{G}_{\Ls}+\bb{Q},\quad\mathrm{where}\quad\forall1\leq d\leq D,[\bb{Q}]_{\cdot,d}\in\mathcal{N}_{\Ls}.
\end{equation}
The null space $\mathcal{N}_{\Ls}$ is characterized in the Proposition \ref{prop:kernel}.
\begin{ Proposition }[\cite{willems1997introduction}[Theorem 3.2.16]]
\label{prop:kernel}
    The determinant $\mathrm{det}(\Ls)$ of an invertible MDO $\Ls=\bb{U}\circ\bb{D}\circ\bb{V}$, where $\bb{U},$$\bb{V}$ are unimodular operators and $\bb{D}$ a diagonal operator, is a polynomial in $\D$. Further, it can always be written as
    \begin{equation}
        \mathrm{det}(\Ls)=\pm\mathrm{det}(\bb{D})=c(\D-\alpha_1)^{n_1}\circ\cdots\circ(\D-\alpha_D)^{n_D},
    \end{equation}
    where $\alpha_r\in\mathbb{C}$ is a root of $\mathrm{det}(\Ls)$ and $n_r$ is its multiplicity. Let $R$ be the number of roots; then, $\mathrm{dim}\left(\mathcal{N}_{\Ls}\right)=\sum_{r=1}^Rn_r$ with $\mathcal{N}_{\bb{L}}=\bb{V}^{-1}\mathcal{N}_{\bb{D}}$, where $\mathcal{N}_{\bb{D}}$ is the null space of $\bb{D}$. 
\end{ Proposition }
Calculations of Green's matrices and null spaces will be given in Section \ref{section:5}.

\subsection{Vector-Valued Splines and Their Native Spaces}
\label{sec:3.3}
We introduce the notion of $\Ls$-splines. These $\Ls$-splines will turn-out to be solutions to our optimization problem.
\begin{Definition}
    Let $\bb{f}\in\mathcal{D}'(\R,\R^{D})$ and $\Ls:\mathcal{D}'\to\mathcal{D}'$ be an invertible MDO. The distribution $\bb{f}$ is an $\Ls$-spline if 
    \begin{equation}
        \Ls\{\bb{f}\}=\sum_{k=1}^K\bb{a}_k\odot\bm{\delta}_{\bb{x}_k},
    \end{equation} where  $\{\bb{a}_k\}_{k=1}^K$ and $\{\bb{x}_k\}_{k=1}^K$ are two lists of coefficients $\bb{a}_k\in\R^D$ and $\bb{x}_k\in\R^D$. In addition, $\bm{\delta}_{\bb{x}_k}\in\mathcal{D}'$ is such that $[\bm{\delta}_{\bb{x}_k}]_d=\delta_{[\bb{x}_k]_d}$. The locations $\bb{x}_k$ are called the knots of the spline with $K$ being the number of knots. If we define the innovation $\bb{m}=\sum_{k=1}^K\bb{a}_k\odot\bm{\delta}_{\bb{x}_k}\in\mathcal{D}'$, we can generate such a spline by
    \begin{equation}
        \bb{f}(\cdot)=\bb{q}(\cdot)+\left(\bb{G}_{\Ls}\ast\bb{m}\right)(\cdot),
    \end{equation}
    where $\bb{q}\in\mathcal{N}_{\Ls}$, and $\bb{G}_{\Ls}$ is the causal Green's function of $\Ls.$
\end{Definition}
In Definition \ref{def:M}, we introduce a $BV$-type Banach space which will be used as a search space.
\begin{Definition}
\label{def:M}
Let $\Ls$ be an invertible MDO. Then, the vector space $\M_{\Ls}(\R,\R^D)$ as a set is defined as
    \begin{equation}
        \M_{\Ls}(\R,\R^D)=\left\{\bb{f}\in\mathcal{D}'(\R)^D:\quad\Ls\{\bb{f}\}\in\M(\R,\R^D)\right\}.
    \end{equation}
\end{Definition}
This space can be equipped with the outer seminorm $\norm{\Ls\{\cdot\}}_{\M,\Zz'}$ or the inner seminorm $\norm{\Ls\{\cdot\}}_{\Zz',\M}$. It is a native space for $\Ls$-splines in the sense that all $\Ls$-splines live inside. This property, as well as other important ones, are gathered in Theorem \ref{prop:MLprop}, whose proof is in Appendix \ref{app:3.3}.
Prior to its statement, we remark that, to norm the space $\M_{\Ls}(\R,\R^D)$, one must be able to project an element $\bb{f}$ onto the null space $\mathcal{N}_{\Ls}\subset\M_{\Ls}(\R,\R^D)$ of $\Ls$, as the existence of the latter is the sole reason why $\norm{\Ls\{\cdot\}}_{\M,\Zz'}$ (or  $\norm{\Ls\{\cdot\}}_{\Zz',\M}$) is not a norm. To this extent, we provide Definition \ref{def:admissible}.

\begin{Definition}
\label{def:admissible}
    Let $\bb{p}=(\bb{p}_n)_{n=1}^N$ be a basis of $\mathcal{N}_{\Ls}$ and $\bm{\phi}=(\bm{\phi}_{n})_{n=1}^N$ be a complementary set of analysis functions. Then, the system $(\bb{p}, \bm{\phi})$ is said to be $\Ls$-admissible if 
    \begin{itemize}
        \item [1.] the basis functions are biorthogonal, with $\forall m,n\in[1\cdots N],\quad\langle \bb{p}_n,\bm{\phi}_m\rangle=\delta_{m,n}$;
        \item [2.] for all $ 1\leq n\leq N,\bm{\phi}_n\in\mathcal{D}(\R,\R^D)$;
        \item [3.] all functions $\bm{\phi}_n$ are compactly supported and causal, such that there exists $\phi^+\in\R^+$ with
        \begin{equation}
            \forall 1\leq d\leq D,\forall 1\leq n\leq N:\quad\mathrm{supp}([\bm{\phi}_n]_d)\subset[0,\phi^+].
        \end{equation}
    \end{itemize}
\end{Definition}
Such systems exist, as proven in Proposition \ref{prop:testfunctions}, whose demonstration is in Appendix \ref{app:3.3}.
\begin{ Proposition }
\label{prop:testfunctions}
There exists a family $(\bm{\phi}_n)_{n=1}^N$ of functions that verify Items 1, 2, and 3 of Definition \ref{def:admissible}.
\end{ Proposition }
We recall that $\mathrm{deg}(\mathrm{det}(\Ls))$ is a polynomial in $\D$ and define $N=\mathrm{dim}(\mathcal{N}_{\Ls})=\mathrm{deg}(\mathrm{det}(\Ls))$. In order to offer flexibility in the specification of the norm on the null space, we introduce the finite-dimensional vector space $\Yy=(\R^N,\norm{\cdot}_\Yy)$ and its dual $\Yy'=(\R^N,\norm{\cdot}_{\Yy'})$. 
\begin{theorem}
\label{prop:MLprop}
Let $\Ls$ be a an invertible MDO, $\norm{\cdot}_{\Zz}$ a monotone norm, and $(\bb{p}, \bm{\phi})$ an $\Ls$-admissible system. Then, the following properties hold.
\begin{itemize}
    \item [1.] The space $\M_{\Ls}(\R,\R^D)$ is a Banach space for the norm 
    \begin{equation}
    \label{eq:norm}
        \norm{\bb{f}}_{\M_{\Ls}}:=\norm{\Ls\{\bb{f}\}}_{\Zz',\M}+\norm{\left( \langle\bb{f},\bm{\phi}_n\rangle\right)_{n=1}^N}_{\mathcal{Y}'}.
    \end{equation} In adequation with our terminology, we call $\norm{\cdot}_{\M_{\Ls}}$ an inner-norm. In the same way, by substituting $\norm{\cdot}_{\M,\Zz'}$ for  $\norm{\cdot}_{\Zz',\M}$, we define an outer-norm on $\M_{\Ls}(\R,\R^D)$.
    \item [2.] The space $\M_{\Ls}(\R,\R^D)$ admits a predual $\C_{\Ls}(\R,\R^D)$, as detailed in Theorem \ref{th:predual}. 
    \item [3.] All seminorms on $\M_{\Ls}(\R,\R^D)$ of the form $\norm{\Ls\{\cdot\}}_{\M,\Zz'}$ or $\norm{\Ls\{\cdot\}}_{\Zz',\M}$  are equivalent. Moreover, the seminorms $\norm{\Ls\{\cdot\}}_{1,\M}$ and $\norm{\Ls\{\cdot\}}_{\M,1}$ are the same.
    \item [4.] The extreme points $\bb{e}$ of the centered unit regularization ball $\mathcal{B}=\left\{[\bb{f}]\in\M_{\Ls}(\R,\R^D)/\mathcal{N}_{\bb{L}}:\quad\norm{\Ls\{\bb{f}\}}_{\M,\Zz'}\leq1\right\}$ are all of the form
    \begin{equation}
    \label{eq:2.2.28}
    [\bb{e}]=\bb{G}_{\Ls}\ast(\bb{w}\odot\bm{\delta}_{\bb{x}})+\mathcal{N}_{\Ls},\quad\quad\bb{x}\in\R^D,
    \end{equation}
    where $\bb{w}$ with $\norm{\bb{w}}_{\Zz'}=1$ is an extreme point of the centered unit ball in $\Zz'$. Correspondingly, the seminorm of an $\Ls$-spline $\bb{f}$ of the generic form $\bb{f}=\sum_{k=1}^K\bb{G}_{\Ls}\ast(\bb{a}_k\odot\bm{\delta}_{\bb{x}_k})+\bb{q}$ with $\{\bb{x}_k\}_{k=1}^K\subset\bb{K}$, such that $\bb{x}_k\neq \bb{x}_{k'}$ (for all rows) for all $k\neq k'$, is given by
    \begin{equation}
    \label{eq:3.3.51}
\norm{\Ls\{\bb{f}\}}_{\M,\Zz'}=\norm{\left(\sum_{k=1}^K\vert [\bb{a}_k]_d\vert\right)_{d=1}^D}_{\Zz'}.
    \end{equation}
    \item [5.] The extreme points $\bb{e}$ of the centered unit regularization ball $\mathcal{B}=\left\{[\bb{f}]\in\M_{\Ls}(\R,\R^D)/\mathcal{N}_{\bb{L}}:\quad\norm{\Ls\{\bb{f}\}}_{\Zz',\M}\leq1\right\}$ are all of the form
    \begin{equation}
        \label{eq:2.2.29}
[\bb{e}]=\bb{G}_{\Ls}\ast(\bb{w}\delta_x)+\mathcal{N}_{\Ls},\quad\quad x\in\R,
    \end{equation}
     where $\bb{w}$ with $\norm{\bb{w}}_{\Zz'}=1$ is an extreme point of the centered unit ball in $\Zz'$. Correspondingly, the seminorm of an $\Ls$-spline $\bb{f}$ of the generic form $\bb{f}=\sum_{k=1}^K\bb{G}_{\Ls}\ast(\bb{a}_k\delta_{x_k})+\bb{q}$ with $\{x_k\}_{k=1}^K\subset\R$ such that $x_k\neq x_{k'}$, is given by
    \begin{equation}
    \label{eq:3.3.53}
\quad\quad\norm{\Ls\{\bb{f}\}}_{\Zz',\M}=\sum_{k=1}^K\norm{\bb{a}_k}_{\Zz'}.
\end{equation}
\end{itemize}
\end{theorem}

Theorem \ref{prop:MLprop} reveals that the search space $\M_{\Ls}(\R,\R^D)$ (Item 1) is not affected by the choice of norm (Item 3), but its geometry is (Items 4 and 5). However, the choice of the norm has an impact on both the weight $\bb{a}$ and the support ($\bb{x}$ or $x$, Items 4 and 5) of the extreme points of the unit ball. By creating an inter-dimensional dependence, inner norms favor sparser extreme points ($D$ free parameters in $\bb{a}$ and $1$ free parameters in $x$) than do outer norms ($D$ free parameters in $\bb{a}$ and $D$ free parameters in $\bb{x}$). In \eqref{eq:3.3.51} and \eqref{eq:3.3.53}, we have calculated the seminorm of linear combinations of extreme points for the outer and inner norms. We shall see that this has an impact on the solutions of our optimization problem. Finally, Item 2 reveals the important fact that $\M_{\Ls}(\R,\R^D)$ admits a $\text{weak}^{\star}$ topology, which is the subject of Section \ref{sec:3.4}. In fact, this topology is the one used for optimization and convergence analysis.

We conclude this section with Proposition \ref{prop:stabilitymat}, where we study the stability of $\M_{\Ls}(\R,\R^D)$ when $\Ls$ is multiplied by an invertible matrix. The proof is given in Appendix \ref{app:3.3}.
\begin{ Proposition }
    \label{prop:stabilitymat}
Let $\bb{A}\in\R^{D\times D}$ be an invertible matrix. 
\begin{itemize}
    \item [1.] If $\tilde{\Ls}=\bb{A}\bb{L}$, then $\bb{G}_{\tilde{\Ls}}=\bb{G}_{\Ls}\bb{A}^{-1}$. In addition, it holds that $\M_{\Ls}(\R,\R^D)=\M_{\tilde{\Ls}}(\R,\R^D)$ and, in particular, $\mathcal{N}_{\Ls}=\mathcal{N}_{\tilde{\Ls}}.$
    \item [2.] If $\tilde{\Ls}=\bb{L}\bb{A}$, then $\bb{G}_{\tilde{\Ls}}=\bb{A}^{-1}\bb{G}_{\Ls}$. In addition, there exist $\bb{S}$ and $\Ls$ such that $\M_{\Ls}(\R,\R^D)\neq\M_{\tilde{\Ls}}(\R,\R^D)$ and, in particular, $\mathcal{N}_{\Ls}\neq\mathcal{N}_{\tilde{\Ls}}.$
\end{itemize}
\end{ Proposition }
We observe that a multiplication on the right is unstable. Intuitively, this is because this multiplication may lead to a permutation of the rows of the Green's matrix, hence modifying the regularity. The stability of the multiplication on the left will be central in Section \ref{sec:4.3}.

\subsection{\texorpdfstring{$\text{Weak}^{\star}$}{WeakStar} Continuity on \texorpdfstring{$\mathcal{M}_{\bb{L}}(\mathbb{R},\mathbb{R}^D)$}{ML(R,RD)}}
\label{sec:weakstar}
\label{sec:3.4}
We shall now equip the space $\M_{\Ls}(\R,\R^D)$ with a suitable topology and notion of continuity. Theorem \ref{th:predual}, whose proof is in Appendix \ref{app:3.4}, makes the $\text{weak}^{\star}$ topology explicit by revealing that $\C_{\Ls}(\R,\R^D)$ is the predual of $\M_{\Ls}(\R,\R^D)$ (Items 1, 2, and 3) and by giving a simple test to verify that a distribution $\bm{\rho}$ is in $\C_{\Ls}(\R,\R^D)$ (Corollary \ref{coro:diracpredual}).
\begin{theorem}
\label{th:predual}
Let $\Ls:\M_{\Ls}(\R,\R^D)\to\M(\R,\R^D)$ be a an invertible MDO, $\norm{\cdot}_{\Zz}$ be a monotone norm, and $(\bb{p}, \bm{\phi})$ be an $\Ls$-admissible system. Then, the function space
\begin{equation}
    \C_{\Ls}(\R,\R^D)=\left\{\bb{g}=\Ls^{\star}\{\bb{v}\}+\sum_{n=1}^Na_n\bm{\phi}_n:\bb{v}\in\C_{0}(\R,\R^D),\bb{a}=(a_n)_{n=1}^N\in\R^N\right\}
\end{equation}
has the following properties:
\begin{itemize}
    \item [1.]it has the direct-sum representation $\C_{\Ls}(\R,\R^D)=\Ls^{\star}\left(\C_0(\R,\R^D)\right)\oplus\mathrm{span}\{\bm{\phi}_n\}_{n=1}^N$;
    \item [2.] it is a Banach space for the norm $ \norm{\bb{g}}_{\C_{\Ls}}=\mathrm{max}(\norm{\bb{v}}_{\Zz,\infty},\norm{\bb{a}}_{\Yy})$ and it is the predual of $\M_{\Ls}(\R,\R^D)$, in the sense that $\C_{\Ls}(\R,\R^D)'=\M_{\Ls}(\R,\R^D)$ when $\M_{\Ls}(\R,\R^D)$ is equipped with the norm
    \begin{equation}
        \norm{\Ls\{\cdot\}}_{\Zz',\M}+\norm{\left(\langle \cdot,\bm{\phi}_n\rangle\right)_{n=1}^N}_{\Yy'};
    \end{equation}
    \item [3.] it is also a Banach space for the equivalent norm $\norm{\bb{g}}_{\C_{\Ls}}=\mathrm{max}(\norm{\bb{v}}_{\infty,\Zz},\norm{\bb{a}}_{\Yy})$ and is the predual of $\M_{\Ls}(\R,\R^D)$, in the sense that $\C_{\Ls}(\R,\R^D)'=\M_{\Ls}(\R,\R^D)$ when $\M_{\Ls}(\R,\R^D)$ is equipped with the norm 
    \begin{equation}
        \norm{\Ls\{\cdot\}}_{\M,\Zz'}+\norm{\left(\langle \cdot,\bm{\phi}_n\rangle\right)_{n=1}^N}_{\Yy'}.
    \end{equation}
\end{itemize}
\end{theorem}

Theorem \ref{th:predual} if of theoretical importance but lacks practicality, as it does not directly tell us which "common" functions or distributions are inside $\C_{\Ls}(\R,\R^D)$. A more practical criterion is provided by Corollary \ref{coro:diracpredual}, whose proof is in Appendix \ref{app:3.4}. In essence, Corollary \ref{coro:diracpredual} tells us that the admissible irregularity of $[\bm{\rho}]_d$, for $\bm{\rho}\in\C_{\Ls}(\R,\R^D)$, is proportional to the regularity of the entries of $[\bb{G}_{\Ls}]_{d,\cdot}$.

\begin{Corollary}
    \label{coro:diracpredual}
    If $\bm{\rho}\in\mathcal{D}'(\R,\R^D)$ is a \emph{compactly supported} distribution such that the vector-valued distribution
$\bb{v}=\bb{G}_{\Ls}^{\star}\ast\bm{\rho}$   
    is representable by a function $\bb{v}\in\C(\R,\R^D)$, then $\bm{\rho}\in\C_{\Ls}(\R,\R^D)$.
\end{Corollary}

\section{Resolution of Vector-Valued Continuous-Domain Inverse Problems}
\label{sec:4}
\subsection{Definition of the Optimization Problem and Set of Assumptions}
\label{sec:4.1}
We are interested in general optimization problems of the form
 
\begin{equation}
\label{eq:3.1.30}
\mathcal{V}=\underset{\bb{f}\in\M_{\Ls}(\R,\R^D)}{\text{argmin}}\quad\mathcal{J}(\bb{f}),\quad\quad\mathcal{J}^{\star}=\underset{\bb{f}\in\M_{\Ls}(\R,\R^D)}{\text{inf}}\mathcal{J}(\bb{f}).
\end{equation}
with the $\emph{loss functional}$
\begin{equation}
\label{eq:3.1.31}
\mathcal{J}(\bb{f})=E(\bb{y},\bb{V}\{\bb{f}\})+\lambda\norm{\Ls\{\bb{f}\}}_{\Zz',\M},\quad\text{or}\quad\mathcal{J}(\bb{f})=E(\bb{y},\bb{V}\{\bb{f}\})+\lambda\norm{\Ls\{\bb{f}\}}_{\M,\Zz'}
\end{equation}
with the following specificities. 
\begin{itemize}
    \item The vector $\bb{y}\in\R^M$, of finite dimension $M$, represents the available measurements.
    \item The operator $\Ls$ encodes the underlying dynamics of the signal to be recovered.
    \item The operator $\bb{V}:\M_{\Ls}(\R,\R^D)\to\R^{M}$ is the forward measurement operator, whose null space is $\mathcal{N}_{\bb{V}}$. It models the data-acquisition process. In finite dimensions, $\bb{V}$ is represented by a vector of vector-valued operators. Here, it is such that $\bb{V}=(\bm{\nu}_m)_{m=1}^M$, with the marginal measurement operator $\bm{\nu}_m=(\nu_{m,d})_{d=1}^D:\M_{\Ls}(\R,\R^D)\to\R$. Thus, the action of $\bb{V}$ on the function $\bb{f}=(f_d)_{d=1}^D$ is described as 
    \begin{equation}
        \bb{V}\left\{\begin{bmatrix}
            f_1\\
            \vdots\\
            f_D
        \end{bmatrix}\right\}=\begin{bmatrix}
            \langle \bm{\nu}_1,\bb{f}\rangle\\
            \vdots\\
            \langle \bm{\nu}_M,\bb{f}\rangle
        \end{bmatrix}=
        \begin{bmatrix}
            \sum_{d=1}^D\langle \nu_{1,d},f_{d}\rangle\\
            \vdots\\
            \sum_{d=1}^D\langle \nu_{M,d},f_{d}\rangle
        \end{bmatrix}=
        \begin{bmatrix}
            \langle \nu_{1,1},\cdot\rangle&...&\langle \nu_{1,D},\cdot\rangle\\
            \vdots& &\vdots\\
            \langle \nu_{M,1},\cdot\rangle&...&\langle \nu_{M,D},\cdot\rangle
        \end{bmatrix}
        \left\{\begin{bmatrix}
            f_1\\
            .\\
            .\\
            f_D
        \end{bmatrix}\right\}.
    \end{equation}
    \item The functional $E$ is the data fidelity functional that measures the discrepancy between the data $\bb{y}$ and the simulated measurement $\bb{V}\{\bb{f}\}$. It is usually tuned to the underlying measurement noise. A popular $E$ has been illustrated in the introduction by the classic squared loss $\norm{\cdot-\cdot}_2^2$.
    \item The metrics $\norm{\Ls\{\cdot\}}_{\Zz',\M}$ and $\norm{\Ls\{\cdot\}}_{\M,\Zz'}$ are seminorms for the search space $\M_{\Ls}(\R,\R^D)$ .
\end{itemize}
In the remaining of this paper, we shall always assume the following set of assumptions to hold.
\begin{itemize}
    \item []\textbf{Assumption 1} The operator $\Ls$ is an invertible MDO.
    \item[] \textbf{Assumption 2} \label{ass:1}The measurement operator $\bb{V}:\M_{\Ls}(\R,\R^D)\to\R^M$ is linear and $\text{weak}^{\star}$-continuous and surjective. 
    \item []\textbf{Assumption 3} \label{ass:2} The measurement operator $\bb{V}$ is injective on the null space $\mathcal{N}_{\Ls}$ of the operator $\Ls$. Equivalently,  $\mathcal{N}_{\bb{V}}\cap\mathcal{N}_{\Ls}=\{0\}$.
    \item []\textbf{Assumption 4} \label{ass:6} If an outer norm $\norm{\cdot}_{\M,\Zz'}$ is used, the underlying norm $\norm{\cdot}_{\Zz'}$ is monotone.
    \item [] \textbf{Assumption 5} \label{ass:4} The loss functional $E:\R^M\times\R^M\to\R^+\cup\{\infty\}$ is proper, coercive, lower semi-continuous, and convex in its second argument.
\end{itemize}

\subsection{Representer Theorem for Outer and Inner Norms}
\label{sec:4.2}
\subsubsection{Main Result}
Functions of the form \eqref{eq:3.3.1.33} will take an important role in Theorem \ref{th:optweakstar}. This form is
\begin{equation}
     \label{eq:3.3.1.33}
        \bb{e}(\cdot)=\bb{q}(\cdot)+\left(\bb{G}_{\bb{\Ls}}\ast\bb{m}\right)(\cdot),\quad\text{with}\quad
        \begin{cases}
        \bb{m}=\sum_{k=1}^{M-N}\bb{a}_k\delta_{x_k},\quad\text{for inner norms},\\
        \bb{m}=\sum_{k=1}^{M-N}\bb{a}_k\odot\bm{\delta}_{\bb{x}_k},\quad\text{for outer norms},
        \end{cases}
\end{equation}
where $\bb{x}_k\in\R^D$ and $x_k\in\R$. In addition, $\bb{q}\in\mathcal{N}_{\Ls}$, $\{\bb{a}_k\}_{k=1}^{M-\N}$ is a list of coefficients $\bb{a}_k\in\R^D$, and each $\frac{\bb{a}_k}{\norm{\bb{a}_k}_{\Zz'}}$ is an extreme point of the centered unit ball in $\Zz'.$ Observe that this is an $\Ls$-spline. We are now in position to state the main theorem, whose proof is in Appendix \ref{app:4.3}.

\begin{theorem}
\label{th:optweakstar}
If Assumptions 1 to 5 are verified, then 
 \begin{equation}
 \label{eq:4.3.2.51}
     \mathcal{V}=\underset{\bb{f}\in\M_{\Ls}(\R,\R^D)}{\text{argmin}}\quad\mathcal{J}(\bb{f})
 \end{equation}
is nonempty, $\mathrm{weak}^{\star}$-compact, and is the $\mathrm{weak}^{\star}$-closed convex hull of its extreme points, which are all of the form \eqref{eq:3.3.1.33}.
\end{theorem}
Theorem \ref{th:optweakstar} is an upgrade to $\R^D$-valued functions of \cite{unser2017splines}[Theorem 1]. It is a strong existence result that shows the optimality of adaptive vector-valued $\Ls$-splines.  Although the optimization problem is deﬁned over the continuum $\M_{\Ls}(\R,\R^D)$, the remarkable takeaway message is that the solution set is the closed convex hull of solutions that are intrinsically sparse, with the level of sparsity equal to the number of measurements minus $N$. Observe further that the sparse solutions are composed of two elements: an adaptive one, specified by the amplitudes $\bb{a}_k$ and the knots $\bb{x}_k$, and a linear regression term $\bb{q}$ that describes the component in the null space of the operator. Since $\bb{q}$ does not contribute to $\norm{\Ls\{\cdot\}}_{\M,\Zz'}$ or $\norm{\Ls\{\cdot\}}_{\Zz', \M}$, the optimization tends to maximize the contribution of the null space component.

\subsubsection{Discussion on the Role of the Norm}
\label{sec:4.2.2}

Theorem \ref{th:optweakstar} reveals the consequence of the choice of the norm. Inner norms, because of their inter-dimensional dependence, yield sparse solutions with knots $\bb{x}_k$ that are only parameterized by $x_{k}\in\R$, with $\bb{x}_k=(x_k)_{d=1}^D$. Hence, spline solutions for inner-norm regularization are parameterized by 
\begin{align}
    &(M-N)D+(M-N)+N=(M-N)(D+1)+N\\
    &\text{(\# of parameters for weights + knots + null space components)}\nonumber
\end{align} parameters. Meanwhile, spline solutions for outer-norm regularization are parameterized by 
\begin{equation}
    (M-N)D+(M-N)D+N=(M-N)2D+N
\end{equation} parameters. This yields a difference of $(M-N)(D-1)$ parameters.  Observe that there is no difference in sparsity when $D=1$, as inner and outer norms are the same in this case. Furthermore, the choice of norm $\norm{\cdot}_{\Zz'}$ has an impact on the amplitudes $\bb{a}_k$, whose normalized versions $\frac{\bb{a}_k}{\norm{\bb{a}_k}}_{\Zz'}$ must be extreme points of the centered unit ball in  $\Zz'$. It follows that, if this ball in $\Zz'$ only has sparse extreme points with $I$ nonzero components, then the extreme-point minimizers for inner norms are parameterized by $(M-N)(I+1)+N$ parameters and, for outer norms, by $(M-N)2I+N$ parameters. The straightforward example is $\norm{\cdot}_1$ with $I=1$.

In summary, the norm $\norm{\cdot}_{\Zz'}$ can sparsify the weights $\bb{a}_k$, while the choice of inner or outer norms sparsifies the knots $\bb{x}_k$. Observe that norms between inner and outer ones can be built in order to sparsify only certain dimensions. To do so, one just needs to take one inner norm and one outer norm acting on complementary sets of indices. Then, one builds a composite norm for which spline knots will be shared (inner-norm behavior) in specific dimensions and free (outer-norm behavior) in others.

\subsection{Extension to Non Invertible MDOs}
\label{sec:4.3}

We extend the representer theorem of Section \ref{sec:4.2} to the class of non invertible MDOs of the form 
\begin{equation}
    \bb{Q}\Ls:\M_{\Ls}(\R,\R^D)\to\M(\R,\R^{D'}),
\end{equation}
where $\Ls$ is a $(D\times D)$ invertible MDO and $\bb{Q}\in\R^{D'\times D}$ is a matrix with $D'\leq D$. If $\bb{Q}$ is invertible and $D'=D$, then so is $\bb{Q}\Ls$ and the results from previous sections are applicable. 
When $\bb{Q}$ is not invertible, we require Assumption 6.
\begin{itemize}
    \item [] \textbf{Assumption 6} The matrix $\bb{Q}\in\R^{D'\times D}$ is such that  $\mathrm{rank}(\bb{Q})=D'$ (i.e., the matrix $\bb{Q}$ is full-rank).
\end{itemize}
The lack of invertibility leads to the issue that, now, $\bb{Q}\bb{\Ls}$ has an infinite-dimensional null space over $\M_{\Ls}(\R,\R^D)$. As remedy, we shall restrict the analysis to that of an appropriate subspace of $\M_{\Ls}(\R,\R^D)$. 
To identify this subspace, we use the Smith decomposition $\bb{Q}=\bb{S}\bb{A}\bb{T},$ where $\bb{T}\in\R^{D\times D},$  $\bb{S}\in\R^{D\times D'}$ are invertible and where, by Assumption 6, we may state without loss of generality that 
$\bb{A}=\begin{bmatrix}
    \bb{I}&\bb{0}
\end{bmatrix},$ with $\bb{I}$ being the $(D'\times D')$ identity matrix. Let $\tilde{\Ls}=\bb{T}\Ls$, which is an invertible MDO with respect to Proposition \ref{prop:stabilitymat}, such that $\M_{\tilde{\Ls}}(\R,\R^{D})=\M_{\Ls}(\R,\R^{D})$.
\begin{Definition}
\label{def:Mrestricted}
The vector space $\M_{\tilde{\Ls}\vert\bb{A}}(\R,\R^{D})$ is, as a set, defined as
\begin{equation}
    \M_{\tilde{\Ls}\vert\bb{A}}(\R,\R^{D})=\left\{\bb{f}\in\M_{\Ls}(\R,\R^{D}):\quad(\bb{I}-\bb{A}^{\intercal}\bb{A})\{\tilde{\Ls}
    \{\bb{f}\}\}=\bb{0}\right\}.
\end{equation}
\end{Definition}
\noindent Clearly, $\M_{\tilde{\Ls}\vert\bb{A}}(\R,\R^{D})$ is a subspace of $\M_{\Ls}(\R,\R^{D})$. It further solves the null-space issue in the sense that
\begin{equation}
  \left\{\bb{f}\in\M_{\tilde{\Ls}\vert\bb{A}}(\R,\R^D):\quad(\bb{Q}\Ls)\{\bb{f}\}=0\right\}=\mathcal{N}_{\Ls}. 
\end{equation}
In Proposition \ref{prop:banachrestricted}, whose proof is in Appendix \ref{app:4.4}, we state that $ \M_{\tilde{\Ls}\vert\bb{A}}(\R,\R^{D})$ is a closed subspace.
\begin{ Proposition }
\label{prop:banachrestricted}
    The vector space $ \M_{\tilde{\Ls}\vert\bb{A}}(\R,\R^{D})$ is closed in the strong topology of $ \M_{\Ls}(\R,\R^{D})$. Therefore, it is a Banach subspace of $\M_{\Ls}(\R,\R^D)$ for the norm \eqref{eq:norm}.
\end{ Proposition }
In Corollary \ref{coro:restrictedRT1}, we study the solution set of regularized optimization problems posed over $\M_{\tilde{\Ls}\vert\bb{A}}(\R,\R^{D})$. The proof is given in Appendix \ref{app:4.4}. For convenience, we introduce the vector space $\Ww'=(\R^{D'},\norm{\cdot}_{\Ww'}),$ where $\norm{\bb{w}}_{\Ww'}:=\norm{\bb{A}^{\intercal}\bb{w}}_{\Zz'}$. Inner norms $\norm{\cdot}_{\Ww',\M}$ and outer  norms $\norm{\cdot}_{\M,\Ww'}$ on $\M(\R,\R^{D'})$ are defined likewise.
\begin{Corollary}
\label{coro:restrictedRT1}
 If Assumptions 1 to 6 are verified, then the solution set
\begin{equation}
    \label{eq:4.3.75}
    \mathcal{V}=\underset{\bb{f}\in\M_{\tilde{\Ls}\vert\bb{A}}(\R,\R^{D})}{\mathrm{argmin}}\quad E(\bb{y},\bb{V}\{\bb{f}\})+\lambda
     \begin{cases}
        \norm{(\bb{Q}\Ls)\{\bb{f}\}}_{\Ww',\M}\\
        \norm{(\bb{Q}\Ls)\{\bb{f}\}}_{\M,\Ww'}
     \end{cases}
\end{equation}
is non empty, $\mathrm{weak}^{\star}$-compact, and is the $\mathrm{weak}^{\star}$-closed convex hull of its extreme points, which are all of the form
 \begin{equation}
    \label{eq:4.3.76}
    \left((\bb{G}_{\Ls}\bb{Q}^{\dagger})\ast\bb{m}\right)(\cdot)+\bb{p}(\cdot)\quad\text{with}\quad
    \begin{cases}
        \bb{m}=\sum_{k=1}^{K}\bb{a}_k\delta_{x_k},\quad\quad\hspace{0.1cm} x_k\in\R\\
        \bb{m}=\sum_{k=1}^{K}\bb{a}_k\odot\delta_{\bb{x}_k},\quad \bb{x}_k\in\R^{D'},
    \end{cases}
\end{equation}
where $\bb{Q}^{\dagger}=\bb{T}^{-1}\bb{A}^{\intercal}\bb{S}^{-1}$ is a right-inverse of $\bb{Q}$, $K\leq (M-N)$, $\bb{p}\in\mathcal{N}_{\Ls}$, and $\frac{\bb{a}_k}{\norm{\bb{a}_k}_{\Ww'}}\in\R^{D'}$ is an extreme point of the centered unit ball in  $\Ww'.$
\end{Corollary}
Corollary \ref{coro:restrictedRT1} is an extension of Theorem \ref{th:optweakstar} when the regularizing MDO is not invertible. In this case, it is $\bb{G}_{\Ls}\bb{Q}^{\dagger}$, a right-inverse of $\bb{Q}\Ls$ that describes the dictionary in which solutions are expressed. The number of degrees of freedom of the extreme-point solutions is adaptive in the sense that, now, $\bb{a}_k\in\R^{D'}.$ Likewise, the comments of Section \ref{sec:4.2.2} are still valid, with the replacement of $D$ by $D'$. When $\bb{Q}=\bb{A}$, only the first $D'$ columns of $\bb{G}_{\Ls}$ are in the description of extreme points. We conclude this section with a system-theory-based result.
In the prequel we have used the equation $\bb{\Ls}\{\bb{f}\}=\bb{m}$ to structure our theory and to study the optimization problem 
\begin{equation}
\label{eq:4.3.opt}
\underset{(\bb{f},\bb{m})\in\M_{\Ls}(\R,\R^{D})\times\M(\R,\R^D)}{\mathrm{argmin}} \mathcal{J}(\bb{f},\bb{m})\quad\quad\text{subject to}\quad\Ls\{\bb{f}\}=\bb{m},
\end{equation}
where 
\begin{equation}
    \mathcal{J}(\bb{f},\bb{m})=E(\bb{y},\bb{V}\{\bb{f}\})+\lambda
     \begin{cases}
        \norm{\bb{m}}_{\Ww',\M}\\
        \norm{\bb{m}}_{\M,\Ww'}.
     \end{cases}
\end{equation}
The operator $\Ls$ represents the dynamics and  $\bb{m}\in\M(\R,\R^D)$ is the unknown control, assumed to be sparse. Nevertheless, in system theory \cite{kailath1980linear,willems1997introduction}, one often has the more intricate structure in
\begin{equation}
\label{eq:4.3.69}
    \bb{\Ls}\{\bb{f}\}=\bb{Q}^{\dagger}\{\bb{m}\},
\end{equation}
where $\bb{Q}^{\dagger}\in\R^{D\times D'}$ is a matrix and $\bb{m}\in\M(\R,\R^{D'})$ is the unknown control. In this case, there are fewer controls than the number of dimensions in the range of $\bb{f}$. In other words, the measure $\Ls\{\bb{f}\}$ lives in a vector subspace of $\M(\R,\R^D)$. In Corollary \ref{coro:restrictedRT2}, we describe the solution set of \eqref{eq:4.3.opt} subject to the new dynamics \eqref{eq:4.3.69}. The proof is given in Appendix \ref{app:4.4}.
\begin{Corollary}
\label{coro:restrictedRT2}
 If Assumptions 1 to 6 are verified, then the solution set in \eqref{eq:4.3.75} is equal to
 \begin{align}
 \label{eq:4.3.2.51}
     \underset{(\bb{f},\bb{m})\in\M_{\tilde{\Ls}\vert\bb{A}}(\R,\R^{D})\times\M(\R,\R^{D'})}{\mathrm{argmin}} \mathcal{J}(\bb{f},\bb{m})\quad\quad\mathrm{subject}\hspace{0.1cm}\mathrm{to}\quad\Ls\{\bb{f}\}=\bb{Q}^{\dagger}\{\bb{m}\}.
 \end{align}
\end{Corollary}

\section{Examples}
\label{section:5}
The null space and the Green's function of an ordinary differential operator $\ls$ are well-known \cite{unser2005cardinal}. In this section, the regularity $\varrho$ of a distribution $f\in\mathcal{D}'(\R)$ is such that
\begin{itemize}
    \item if $f\in\C^{N}(\R)$, then $\varrho=N;$
    \item if 
    $f\in\mathcal{L}_{\infty,loc}(\R)$ then $\varrho=(-1);$
    \item if $f$ is composed of the $N$th derivative of a Dirac mass, then $\varrho=-(N+2).$
\end{itemize}
We provide two examples for the theory developed in this paper: i) the case when $\Ls$ is a diagonal MDO; and ii) the case when $\Ls$ is a $(2\times2)$ MDO. Observe that, in the case of a general MDO $\Ls$ of dimensions $(D\times D)$, the procedure to calculate the Green's matrix, the null space, and the spline solutions remains the same. But the calculations complexity increases in $D$. 

Observe also that, under the identification $\R^2\cong\mathbb{C},$ the special case of the reconstruction of $\R^2$-valued functions translates into the reconstruction of $\mathbb{C}$-valued functions. Likewise, $\mathbb{C}$-valued splines can also be defined. 
\subsection{Diagonal \texorpdfstring{$\mathbf{L}$}{L}}
\label{sec:5.1}
Let $\Ls$ be the diagonal invertible MDO characterized by $[\Ls]_{d,d'}=0,$ $\forall d\neq d'$ and call the $d$th diagonal element $\ls_d$. It follows directly from the definition of a Green's matrix $\left(\bb{L}\{\bb{G}_{\Ls}\}=\bb{I}\delta\right)$ that 
\begin{equation}
   \bb{G}_{\Ls}=\textbf{Diag}(g_{1},...,g_{D}), 
\end{equation}
where $g_{d}$ is the Green's function of $\ls_d$. The null space $\mathcal{N}_{\Ls}$ is the Cartesian product of the marginal null spaces, so that
\begin{equation}
    \mathcal{N}_{\Ls}=\mathcal{N}_{\ls_1}\times...\times\mathcal{N}_{\ls_D}=\prod_{d=1}^D\mathcal{N}_{\ls_d},
\end{equation}
with $\mathcal{N}_{\ls_d}$ the null space of $\ls_d$. For $\bb{a}_{k}=(a_{k,d})_{d=1}^D$, the spline solutions to the inner-norm-regularized problem \eqref{eq:3.1.30} are of the form
\begin{equation}
    \bb{q}(\cdot)+\sum_{k=1}^{M-N}\left(a_{k,d}g_{\ls_d}(\cdot-x_{k})\right)_{d=1}^D,\quad\bb{q}\in\prod_{d=1}^D\mathcal{N}_{\ls_d}.
\end{equation}
For $\bb{a}_{k}=(a_{k,d})_{d=1}^D$ and $\bb{x}_{k}=(x_{k,d})_{d=1}^D$, the spline solutions to the outer-norm-regularized problem \eqref{eq:3.1.30} are of the form 
\begin{equation}
    \bb{q}(\cdot)+\sum_{k=1}^{M-N}\left(a_{k,d}g_{\ls_d}(\cdot-x_{k,d})\right)_{d=1}^D,\quad\bb{q}\in\prod_{d=1}^D\mathcal{N}_{\ls_d}.
\end{equation}

\subsection{\texorpdfstring{($2\times2$)}{(2x2)} Matrix \texorpdfstring{$\mathbf{L}$}{L}}
Let $\Ls$ be the $(2\times2)$ invertible MDO
    $\Ls=\begin{bmatrix}
        \ls_{1,1}&\ls_{1,2}\\
        \ls_{2,1}&\ls_{2,2}
    \end{bmatrix}$.
\subsubsection{Green's Matrix}
The invertibility condition $\mathrm{det}(\Ls)\neq0$
translates into $(\ls_{1,1}\ls_{2,2}-\ls_{1,2}\ls_{2,1})\neq0$. The Green's matrix is then found by inversion of a $(2\times2)$ matrix, as
\begin{equation}
    \bb{G}_{\Ls}(\cdot)=\begin{bmatrix}
        \ls_{2,2}&-\ls_{1,2}\\
        -\ls_{2,1}&\ls_{1,1}
    \end{bmatrix}\left\{
    \begin{bmatrix}
        (\ls_{1,1}\ls_{2,2}-\ls_{1,2}\ls_{2,1})^{-1}\delta(\cdot)\\
        (\ls_{1,1}\ls_{2,2}-\ls_{1,2}\ls_{2,1})^{-1}\delta(\cdot)
    \end{bmatrix}
    \right\}
    =\begin{bmatrix}
        \ls_{2,2}&-\ls_{1,2}\\
        -\ls_{2,1}&\ls_{1,1}
    \end{bmatrix}\left\{\begin{bmatrix}
        g_{\mathrm{det}(\Ls)}(\cdot)\\
        g_{\mathrm{det}(\Ls)}(\cdot)
    \end{bmatrix}\right\},
\end{equation}
where $g_{\mathrm{det}(\Ls)}=(\ls_{1,1}\ls_{2,2}-\ls_{1,2}\ls_{2,1})^{-1}\delta$ is the causal Green's function of the operator $\mathrm{det}(\Ls)$. The regularity  of the entries of $\bb{G}_{\Ls}$ is given by 
\begin{equation}
\label{eq:5.2.58}
    \begin{bmatrix}
        \varrho([\bb{G}_{\Ls}]_{1,1})&\varrho([\bb{G}_{\Ls}]_{1,2})\\
        \varrho([\bb{G}_{\Ls}]_{2,1})&\varrho([\bb{G}_{\Ls}]_{2,2})
    \end{bmatrix}=N-2-\begin{bmatrix}
        N_{2,2}&N_{1,2}\\
        N_{2,1}&N_{1,1}
    \end{bmatrix},
\end{equation}
where $N$ is the order of $\mathrm{det}(\Ls)$ and $N_{d,d'}$ is the order of $\ls_{d,d'}.$
\subsubsection{Null Space}
The null space $\mathcal{N}_{\Ls}$ of $\Ls$ is
\begin{equation}
    \tilde{\mathcal{N}_{\Ls}}=\begin{bmatrix}
        \ls_{2,2}&-\ls_{1,2}\\
        -\ls_{2,1}&\ls_{1,1}
    \end{bmatrix}\left\{\mathcal{N}_{\mathrm{det}(\Ls)}^2\right\}=\left\{\begin{bmatrix}
        \ls_{2,2}\{f\}-\ls_{1,2}\{g\}\\
        -\ls_{2,1}\{f\}+\ls_{1,1}\{g\}
    \end{bmatrix}:\quad f,g\in\mathcal{N}_{\mathrm{det}(\Ls)}
    \right\}\label{eq:5.2.59},
\end{equation}
where $\mathcal{N}_{\mathrm{det}(\Ls)}$ is the null space of $\mathrm{det}(\Ls)$. To prove that $\mathcal{N}_{\Ls}=\tilde{\mathcal{N}_{\Ls}}$, we first observe that $\begin{bmatrix}
        \ls_{2,2}&-\ls_{1,2}\\
        -\ls_{2,1}&\ls_{1,1}
    \end{bmatrix}$ is the transpose of the comatrice $\mathbf{com}(\Ls)$ of $\Ls$, whose definition is given in Proposition \ref{prop:Greens}. This implies that, $\forall 
    f,g\in\mathcal{N}_{\mathrm{det}(\Ls)},$
\begin{equation}
    \Ls\circ\mathbf{com}(\Ls)^{\intercal}\left\{\begin{bmatrix}
        f\\
        g
    \end{bmatrix}\right\}=\textbf{diag}(\mathrm{det}(\Ls))\left\{\begin{bmatrix}
        f\\
        g
    \end{bmatrix}\right\}=\begin{bmatrix}
        \mathrm{det}(\Ls)\{f\}\\
        \mathrm{det}(\Ls)\{g\}
    \end{bmatrix}=\bm{0}\quad\Rightarrow\quad\tilde{\mathcal{N}_{\Ls}}\subset\mathcal{N}_{\Ls}.
\end{equation}
To conclude, we show that $\mathrm{dim}(\tilde{\mathcal{N}_{\Ls}})=\mathrm{dim}(\mathcal{N}_{\Ls})$. To do so, observe that $\mathrm{dim}\left(\mathcal{N}_{\mathrm{det}(\Ls)}^2\right)=2\mathrm{dim}\left(\mathcal{N}_{\Ls}\right)$ and that the determinant of $\mathbf{com}(\Ls)^{\intercal}$ is $\mathrm{det}(\Ls)$. It follows from Proposition \ref{prop:kernel} that the dimension of the null space of $\mathbf{com}(\Ls)^{\intercal}$ is $\mathrm{dim}(\mathcal{N}_{\mathrm{det}(\Ls)})=\mathrm{dim}(\mathcal{N}_{\Ls})$. From the definition $\tilde{\mathcal{N}_{\Ls}}=\mathbf{com}(\Ls)^{\intercal}\left\{\mathcal{N}_{\mathrm{det}(\Ls)}^2\right\}$, it follows that $\mathrm{dim}(\tilde{\mathcal{N}_{\Ls}})\geq\left(2\mathrm{dim}\left(\mathcal{N}_{\Ls}\right)-\mathrm{dim}\left(\mathcal{N}_{\Ls}\right)\right)=\mathrm{dim}\left(\mathcal{N}_{\Ls}\right)$. This proves that $\mathcal{N}_{\Ls}=\tilde{\mathcal{N}_{\Ls}}$.

\subsubsection{Spline Solutions}
According to our previous calculations and Theorem \ref{th:optweakstar}, the spline solutions to the inner-norm-regularized problem \eqref{eq:3.1.30} are of the form 
\begin{equation}
    \bb{q}(\cdot)+\sum_{k=1}^{M-N}\begin{bmatrix}
        \ls_{2,2}&-\ls_{1,2}\\
        -\ls_{2,1}&\ls_{1,1}
    \end{bmatrix}\left\{\bb{a}_kg_{\mathrm{det}(\Ls)}(\cdot-x_k)\right\},\quad\bb{q}\in\mathcal{N}_{\Ls},
\end{equation}
with $\frac{\bb{a}_k}{\norm{\bb{a}_k}_{\Zz'}}$ an extreme point of the centered unit ball in $\Zz'$. Likewise, for $\bb{a}_{k}=(a_{k,d})_{d=1}^2$ and $\bb{x}_{k}=(x_{k,d})_{d=1}^2$, the spline solutions to the outer-norm-regularized problem \eqref{eq:3.1.30} are of the form
\begin{equation}
    \bb{q}(\cdot)+\sum_{k=1}^{M-N}\begin{bmatrix}
        \ls_{2,2}&-\ls_{1,2}\\
        -\ls_{2,1}&\ls_{1,1}
    \end{bmatrix}\left\{\begin{bmatrix}
        a_{k,1}g_{\mathrm{det}(\Ls)}(\cdot-x_{k,1})\\
        a_{k,2}g_{\mathrm{det}(\Ls)}(\cdot-x_{k,2})
    \end{bmatrix}\right\},\quad\bb{q}\in\mathcal{N}_{\Ls},
\end{equation}
with $\frac{\bb{a}_k}{\norm{\bb{a}_k}_{\Zz'}}$ an extreme point of the centered unit ball in $\Zz'.$

\subsection{First-Order Systems}
\label{ex:firstorder}
\subsubsection{General Description}

We choose an MDO $\Ls$ that describes a first-order dynamical system as in 
\begin{equation}
\label{eq:firstorderdynamic}
    \Ls\{\bb{f}\}=\bb{P}\{\bb{m}\},\quad\text{with}\quad\Ls=(\bb{I}\D-\bb{A}),
\end{equation}
where $\bb{A}\in\R^{D\times D}$, $\bb{P}\in\R^{D\times D'}$ is a full rank matrix, and $\bb{m}\in\M(\R,\R^{D'})$. In Proposition \ref{prop:firstorder}, we state that our assumption that $\Ls$ is invertible is always verified when it describes a first-order system.
\begin{ Proposition }
\label{prop:firstorder}
For any matrix $\bb{A}\in\R^{D\times D}$, one has that 
\begin{equation}
    \mathrm{det}(\Ls)=\mathrm{det}(\bb{I}\D-\bb{A})\neq0.
\end{equation}
\end{ Proposition }
\begin{proof}[\textbf{Proof of Proposition \ref{prop:firstorder}}]
Let $\bb{A}=\bb{U}\bb{J}\bb{U}^{-1}$ with $\bb{J}$ the Jordan form of $\bb{A},$ with Jordan blocks $\{\bb{J}_k\}_{k=1}^K$. Each block $\bb{J}_k$ is of size $(N_k\times N_k)$, associated to the eigenvalue $\lambda_k$. Then, we calculate that 
\begin{equation}
\mathrm{det}(\Ls)=\mathrm{det}(\bb{U}(\bb{I}\D-\bb{J})\bb{U}^{-1})=\mathrm{det}(\bb{I}\D-\bb{J})=\Pi_{k=1}^K(\D-\lambda_k)^{N_k}\neq0.
\end{equation}
\end{proof}
\noindent Thus, for any matrix $\bb{A}\in\R^{D\times D}$, Corollaries \ref{coro:restrictedRT1} and \ref{coro:restrictedRT2} show that the extreme-point solutions of regularized IPs, constrained by the dynamics \eqref{eq:firstorderdynamic}, are of the form
\begin{equation}
    \label{eq:solEA}
    \bb{p}+(\bb{G}_{\Ls}\bb{P})\ast\bb{m}=\bb{p}+(u(\cdot)\mathrm{e}^{\bb{A}\cdot}\bb{P})\ast\bb{m},
\end{equation}
where $\mathrm{e}^{\bb{A}\cdot}$ is the exponential of the matrix $\bb{A}$ and $\bb{p}\in\mathcal{N}_{\Ls}$, with this null space being spanned by the columns of $\mathrm{e}^{\bb{A}\cdot}$. In addition to this elegant description of solutions, we note that the invertibility condition is weak. Typically, to study the system \eqref{eq:firstorderdynamic} when $D'=1$, one assumes that it is controllable in the sense that \cite{kailath1980linear}[page 86]
\begin{equation}
    \mathrm{dim}\left(\mathrm{span}\{\bb{A}^{d}\bb{P}\bb{1}\}_{d=0}^{D-1}\right)=D.
\end{equation}
The controllability assumption is also used to construct $\mathcal{L}_2$ and $\mathcal{L}_1$ smoothing splines \cite{nagahara20141,sun2000control}, which are (optimal) solutions of regularized IPs. Far from being trivial, it imposes implicit constraints on the couple $(\bb{A},\bb{P})$, as illustrated in Section \ref{sec:5.4}.

\subsubsection{Non-Controllability}
\label{sec:5.4}
We illustrate the fact that the invertibility of $\Ls$ is a condition that is weaker than that of controllability. Let a system be composed of two point masses $(p_1,p_2)$ with positions $(x_1,x_2)$, velocity $(\D x_1,\D x_2)$, acceleration $(\D^2x_1,\D^2x_2)$, and unit mass $(1,1)$. These two points are connected, over a frictionless surface, by a damper with damping coefficient $d$. Only $p_1$ is subject to a control $m$. The application of Newton's second law yields
\begin{equation}
\label{eq:ex1ODE}
    \begin{cases}
        \D^{2}x_1=-d(\D x_1-\D x_2)+m\\
        \D^{2}x_2=-d(\D x_2-\D x_1)
    \end{cases}\quad\Rightarrow\quad(\bb{I}\D-\bb{A})\{\bb{f}\}=\bb{P}\{m\},
\end{equation}
where the state function $\bb{f}$ and the matrices $\bb{A},\bb{P}$ are specified as  
\begin{equation}
    \bb{f}=\begin{bmatrix}
        f_1\\
        f_2\\
        f_3\\
        f_4
    \end{bmatrix}
    =\begin{bmatrix}
        x_1\\
        x_2\\
        \D x_1\\
        \D x_2
    \end{bmatrix},\quad\quad
    \bb{A}=
    \begin{bmatrix}
        0&0&1&0\\
        0&0&0&1\\
        0&0&-d&d\\
        0&0&d&-d
    \end{bmatrix},\quad\quad
    \bb{P}=
        \begin{bmatrix}
        0\\
        0\\
        1\\
        0
    \end{bmatrix}.
\end{equation}
We observe that $(\bb{A},\bb{P})$ is not controllable 
\begin{equation}
    \mathrm{dim}(\text{span}\{\bb{P},\bb{A}\bb{P},\bb{A}^2\bb{P},\bb{A}^3\bb{P}\})=\mathrm{dim}\left(\text{span}\left\{
    \begin{bmatrix}
        0\\
        0\\
        1\\
        0
    \end{bmatrix},
    \begin{bmatrix}
        1\\
        0\\
        -d\\
        d
    \end{bmatrix},
    \begin{bmatrix}
        -d\\
        d\\
        2d^2\\
        -2d^2
    \end{bmatrix}
    \begin{bmatrix}
        2d^2\\
        -2d^2\\
        -4d^3\\
        4d^3
    \end{bmatrix}
    \right\}\right)=3,
\end{equation}
because $(-2d\bb{A}^2\bb{P})=\bb{A}^3\bb{P}$, while $\Ls$ is invertible
\begin{equation}
    \Ls=
    \begin{bmatrix}
        \D & 0 & -1 & 0\\
        0 & \D & 0 & -1\\
        0 & 0 & \D + d & -d\\
        0 & 0 & -d & \D + d
    \end{bmatrix}\quad\Rightarrow\quad\mathrm{det}(\Ls)=\D^3(\D+2d)\neq0.
\end{equation}
We define $h_{\pm}(t)=\frac{1}{2d}(1\pm\mathrm{e}^{-2dt})$ and conclude this example by calculating that
\begin{align}
        \bb{G}_{\Ls}(t)&=
    u(t)\begin{bmatrix}
        1 & & 0 & & \frac{1}{2}(h_-(t)+t) & & \frac{1}{2}(-h_-(t)+t)\\
        0 & & 1 & & \frac{1}{2}(-h_-(t)+t) & &\frac{1}{2}(h_-(t)+t)\\
        0 & & 0 & &  dh_+(t) & & dh_-(t)\\
        0 & & 0 & & dh_-(t) & & dh_+(t)
    \end{bmatrix}.
\end{align}
When the control is null $(m=0)$, the solution of the system of ODE $\eqref{eq:ex1ODE}$ is given by an element of the null space $\mathcal{N}_{\Ls}$ whose coefficients are determined by the initial conditions of the system. These can be the position and velocity of $(p_1,p_2)$ at $t=0$. When the control $m\neq0$ is unknown, $m\in\M(\R)$ and outer and inner norms yield the same form of solutions. The latter is \eqref{eq:solEA} with $\bb{m}=m=\sum_{k=1}^Ka_k\delta_{x_k}$ and $K\leq(M-4)$.

\subsubsection{L1 versus L2}
In Section \ref{sec:5.4}, for a chosen example, we described the extreme-point solutions of the optimization problem \eqref{eq:4.3.2.51}, which is regularized with $\norm{\cdot}_{\M(\R)}$. We now provide, for the same example, the description of the unique solution when the regularization is done via $\norm{\cdot}_{\mathcal{L}_2(\R)}^2$. The solution  is the control-theoretic smoothing spline \cite{sun2000control}
\begin{equation}
\label{eq:solL2}
\bb{p}+(\bb{G}_{\Ls}\bb{P})\ast m,\quad m=\sum_{m=1}^{M}a_m\bb{P}^{\intercal}\bb{G}_{\Ls^{\star}}\ast\bm{\nu}_m^{\intercal},\quad\quad\bb{p}\in\mathcal{N}_{\Ls},
\end{equation}
where $\bb{p}$ and the $a_k$ are typically determined numerically. In \eqref{eq:solL2}, $m$ requires the knowledge of $\bb{G}_{\Ls^{\star}}$ to be evaluated. The MDO $\Ls^{\star}$ and its determinant are 
\begin{equation}
    \Ls^{\star}=
    -\begin{bmatrix}
        \D & 0 & 0 & 0\\
        0 & \D & 0 & 0\\
        1 & 0 & \D - d & d\\
        0 & 1 & d & \D - d
\end{bmatrix}\quad\Rightarrow\quad\mathrm{det}(\Ls^{\star})=\D^3(\D-2d)\neq0
\end{equation}
and, for $g_{\pm}(t)=\frac{-1}{2d}(1\pm\mathrm{e}^{2dt})$, its Green's matrix is
\begin{align}
        \bb{G}_{\Ls^{\star}}(t)&=
    u(t)\begin{bmatrix}
        -1 & & 0 & & 0 & & 0\\
        0 & & -1 & & 0 & &0\\
        \frac{1}{2}(g_-(t)+t) & & \frac{-1}{2}(g_-(t)-t) & &  dg_+(t) & & dg_-(t)\\
        \frac{-1}{2}(g_-(t)-t) & & \frac{1}{2}(g_-(t)+t) & & dg_-(t) & & dg_+(t)
    \end{bmatrix}.
\end{align}
The measurement functionals typically only involve the evaluation of the position $(x_1,x_2)$ and can therefore be described as 
$\bm{\nu}_m=(b_m\delta_{t_{1,m}},c_m\delta_{t_{2,m}},0,0).$ 
Thus, the unique solution is 
\begin{align}
\label{eq:l2sol}
    \bb{p}&+\sum_{m=1}^{M}a_mb_m    \begin{bmatrix}
        \frac{1}{2}(h_-(\cdot)u(\cdot)+(\cdot)_+)\\
       \frac{1}{2}(-h_-(\cdot)u(\cdot)+(\cdot)_+)\\
       dh_+(\cdot)u(\cdot)\\
        dh_-(\cdot)u(\cdot)
    \end{bmatrix}\ast\left(g_-(\cdot-t_{1,m})u(\cdot-t_{1,m})+(\cdot-t_{1,m})_+)\right)\nonumber\\
    &-\sum_{m=1}^{M}a_mc_m    \begin{bmatrix}
        \frac{1}{2}(h_-(\cdot)u(\cdot)+(\cdot)_+)\\
       \frac{1}{2}(-h_-(\cdot)u(\cdot)+(\cdot)_+)\\
       dh_+(\cdot)u(\cdot)\\
        dh_-(\cdot)u(\cdot)
    \end{bmatrix}\ast\left(g_-(\cdot-t_{2,m})u(\cdot-t_{2,m})-(\cdot-t_{2,m})_+)\right).
\end{align}
As concluding notes, we remark that the solution from the $\mathcal{L}_2$ regularization is smoother, in each dimension, due to the convolution in \eqref{eq:l2sol}. Further, the difficulty of the calculation of the solution, because of the convolution between $\bb{G}_{\Ls}$ and $\bb{G}_{\Ls^{\star}}$, can be prohibitive for large matrices.

\section{Conclusion}
In this paper we extended the regularized-IP-based framework for the reconstruction of scalar-valued functions to the reconstruction of vector-valued functions. Although our construction shares many similarities with the $D=1$ theory, it has additional fundamental complexities.

First of all, the construction of a pseudo-inverse operator $\Ls^{-1}$ through its Green's matrix $\bb{G}_{\Ls}$ appeared to be significantly more involved as its existence is typically not guaranteed when $\mathrm{det}(\Ls)=0$ (Theorem \ref{prop:greens2}), and as the null space structure is more intricate (Proposition \ref{prop:kernel}). Furthermore, because $\bb{G}_{\Ls}$ is richer than the concatenation of one-dimensional Green's functions, its calculus required additional tools such as matrix inversion (Section \ref{section:5}).

Second of all, with consideration to the varying degrees of regularity in each and across dimensions, we showed that functions in the search space $\M_{\Ls}(\R,\R^D)$ have a more complex functional representation (Section \ref{sec:3.3}).   

Finally, there exists for the choice of norm the same complexity gap from $\M_{\mathrm{L}}(\R)$ to $\M_{\Ls}(\R,\R^D)$ than from $\R$ to $\R^{D}$: although all norms (inner, outer or mixed) on $\M_{\Ls}(\R,\R^D)$ generate the same topology, they yield different geometries. Consequently, there is no canonical choice of norm and each of them have properties that may or may not be advantageous to the practitioner. With this prospect, we provided a specific analysis of the yielded geometries (Theorem \ref{prop:MLprop}) and their influences on the structure of the extreme point solutions (Theorem \ref{th:optweakstar}), allowing one to make an educated choice of norm.

\subsection{Computational Challenges}
The analogy with the one-dimensional theory raises the following fundamental question: Is it possible to devise an efficient numerical scheme for the resolution of the optimization problem in \eqref{eq:4.3.75}? The outcome of Theorem \ref{th:optweakstar} and Corollary \ref{coro:restrictedRT1} is that the solution of the optimization problem is given in a dictionary that is not countable. A natural discretization is the finite-dimensional set of equally spaced Green's matrices
\begin{equation}
\label{eq:6.1.96}
    \left\{\bb{G}_{\Ls}\left(\cdot-hk\right)\right\}_{k=K_1}^{K_2},
\end{equation}
to which classic optimization algorithms (for example, FISTA \cite{beck2009fast2}, primal-dual splitting \cite{condat2013primal}) could be applied to find the optimal coefficients $\bb{a}_k$. The $1$D theory suggests that this approach is feasible \cite{debarre2019b,guillemet2025convergence}. Nevertheless, this theory represents the $1$D counterpart of \eqref{eq:6.1.96} in its equivalent $B$-spline basis whose elements, because of their compact support, are much better conditioned than the Green's functions. 
For a numerical resolution, the key theoretical challenge is now to establish the existence and the expression of matrix $B$-splines.

\section{Acknowledgement}
Vincent Guillemet was supported by the Swiss National Science Foundation (SNSF) under grant 200020\_219356.

\bibliographystyle{abbrv}
\bibliography{ref}

\begin{thebibliography}{10}

\bibitem{adcock2018infinite}
B.~Adcock.
\newblock Infinite-dimensional compressed sensing and function interpolation.
\newblock {\em Foundations of Computational Mathematics}, 18(3):661--701, 2018.

\bibitem{adcock2016generalized}
B.~Adcock and A.~C. Hansen.
\newblock Generalized sampling and infinite-dimensional compressed sensing.
\newblock {\em Foundations of Computational Mathematics}, 16:1263--1323, 2016.

\bibitem{alvarez2012kernels}
M.~A. Alvarez, L.~Rosasco, N.~D. Lawrence, et~al.
\newblock Kernels for vector-valued functions: A review.
\newblock {\em Foundations and Trends{\textregistered} in Machine Learning}, 4(3):195--266, 2012.

\bibitem{bartolucci2023understanding}
F.~Bartolucci, E.~De~Vito, L.~Rosasco, and S.~Vigogna.
\newblock Understanding neural networks with reproducing kernel banach spaces.
\newblock {\em Applied and Computational Harmonic Analysis}, 62:194--236, 2023.

\bibitem{bauer1961absolute}
F.~L. Bauer, J.~Stoer, and C.~Witzgall.
\newblock Absolute and monotonic norms.
\newblock {\em Numerische Mathematik}, 3(1):257--264, 1961.

\bibitem{beck2009fast2}
A.~Beck and M.~Teboulle.
\newblock A fast iterative shrinkage-thresholding algorithm for linear inverse problems.
\newblock {\em SIAM journal on imaging sciences}, 2(1):183--202, 2009.

\bibitem{beffa2024weakly}
F.~Beffa.
\newblock {\em Weakly Nonlinear Systems: With Applications in Communications Systems}.
\newblock Springer Nature, 2024.

\bibitem{bertero2021introduction}
M.~Bertero, P.~Boccacci, and C.~De~Mol.
\newblock {\em Introduction to Inverse Problems in Imaging}.
\newblock CRC press, 2021.

\bibitem{boyer2019representer}
C.~Boyer, A.~Chambolle, Y.~D. Castro, V.~Duval, F.~De~Gournay, and P.~Weiss.
\newblock On representer theorems and convex regularization.
\newblock {\em SIAM Journal on Optimization}, 29(2):1260--1281, 2019.

\bibitem{bredies2020sparsity}
K.~Bredies and M.~Carioni.
\newblock Sparsity of solutions for variational inverse problems with finite-dimensional data.
\newblock {\em Calculus of Variations and Partial Differential Equations}, 59(1):14, 2020.

\bibitem{bredies2023generalized}
K.~Bredies, M.~Carioni, S.~Fanzon, and F.~Romero.
\newblock A generalized conditional gradient method for dynamic inverse problems with optimal transport regularization.
\newblock {\em Foundations of Computational Mathematics}, 23(3):833--898, 2023.

\bibitem{9655475}
J.~Campos, S.~Aziznejad, and M.~Unser.
\newblock Learning of continuous and piecewise-linear functions with hessian total-variation tegularization.
\newblock {\em IEEE Open Journal of Signal Processing}, 3:36--48, 2022.

\bibitem{condat2013primal}
L.~Condat.
\newblock A primal--dual splitting method for convex optimization involving lipschitzian, proximable and linear composite terms.
\newblock {\em Journal of optimization theory and applications}, 158(2):460--479, 2013.

\bibitem{de1966splines}
C.~De~Boor and R.~E. Lynch.
\newblock On splines and their minimum properties.
\newblock {\em Journal of Mathematics and Mechanics}, 15(6):953--969, 1966.

\bibitem{debarre2019b}
T.~Debarre, J.~Fageot, H.~Gupta, and M.~Unser.
\newblock B-spline-based exact discretization of continuous-domain inverse problems with generalized tv regularization.
\newblock {\em IEEE Transactions on Information Theory}, 65(7):4457--4470, 2019.

\bibitem{dobrakov1971representation}
I.~Dobrakov.
\newblock On representation of linear operators on $c_0(t,\bm{X})$.
\newblock {\em Czechoslovak Mathematical Journal}, 21(1):13--30, 1971.

\bibitem{durbin2012time}
J.~Durbin and S.~J. Koopman.
\newblock {\em Time Series Analysis by State Space Methods}, volume~38.
\newblock OUP Oxford, 2012.

\bibitem{egerstedt2009control}
M.~Egerstedt and C.~Martin.
\newblock {\em Control Theoretic Splines: Optimal Control, Statistics, and Path Planning}, volume~31.
\newblock Princeton University Press, 2009.

\bibitem{fageot2020tv}
J.~Fageot and M.~Simeoni.
\newblock Tv-based reconstruction of periodic functions.
\newblock {\em Inverse Problems}, 36(11):115015, 2020.

\bibitem{fernandez2016super}
C.~Fernandez-Granda.
\newblock Super-resolution of point sources via convex programming.
\newblock {\em Information and Inference: A Journal of the IMA}, 5(3):251--303, 2016.

\bibitem{fisher1975spline}
S.~D. Fisher and J.~W. Jerome.
\newblock Spline solutions to l1 extremal problems in one and several variables.
\newblock {\em Journal of Approximation Theory}, 13(1):73--83, 1975.

\bibitem{flinth2019exact}
A.~Flinth and P.~Weiss.
\newblock Exact solutions of infinite dimensional total-variation regularized problems.
\newblock {\em Information and Inference: A Journal of the IMA}, 8(3):407--443, 2019.

\bibitem{guillemet2025convergence}
V.~Guillemet, J.~Fageot, and M.~Unser.
\newblock Convergence analysis of the discretization of continuous-domain inverse problems.
\newblock {\em Inverse Problems}, 41(4):045008, 2025.

\bibitem{gupta2018continuous}
H.~Gupta, J.~Fageot, and M.~Unser.
\newblock Continuous-domain solutions of linear inverse problems with tikhonov versus generalized tv regularization.
\newblock {\em IEEE Transactions on Signal Processing}, 66(17):4670--4684, 2018.

\bibitem{kailath1980linear}
T.~Kailath.
\newblock {\em Linear Systems}, volume 156.
\newblock Prentice-Hall Englewood Cliffs, NJ, 1980.

\bibitem{nagahara2020sparsity}
M.~Nagahara.
\newblock {\em Sparsity Methods for Systems and Control}.
\newblock now Publishers, 2020.

\bibitem{nagahara2013monotone}
M.~Nagahara and C.~F. Martin.
\newblock Monotone smoothing splines using general linear systems.
\newblock {\em Asian Journal of Control}, 15(2):461--468, 2013.

\bibitem{nagahara20141}
M.~Nagahara and C.~F. Martin.
\newblock $l^1$ control theoretic smoothing splines.
\newblock {\em IEEE Signal Processing Letters}, 21(11):1394--1397, 2014.

\bibitem{reed1972methods}
M.~Reed.
\newblock Methods of modern mathematical physics i.
\newblock {\em Functional analysis}, 1972.

\bibitem{reed2003methods}
M.~Reed, B.~Simon, and S.~Reed.
\newblock {\em Methods of Modern Mathematical Physics: Fourier Analysis, Self-Adjointness}.
\newblock Academic Press Harcourt Brace Jovanovich, 1975.

\bibitem{rudin}
W.~Rudin.
\newblock {\em Real and Complex Analysis}.
\newblock New York: McGraw-Hill, 3rd ed. edition, 1987.

\bibitem{rugh1996linear}
W.~J. Rugh.
\newblock {\em Linear System Theory}.
\newblock Prentice-Hall, Inc., 1996.

\bibitem{scholkopf2001generalized}
B.~Sch{\"o}lkopf, R.~Herbrich, and A.~J. Smola.
\newblock A generalized representer theorem.
\newblock In {\em International Conference On Computational Learning Theory}, pages 416--426. Springer, 2001.

\bibitem{shenouda2024variation}
J.~Shenouda, R.~Parhi, K.~Lee, and R.~D. Nowak.
\newblock Variation spaces for multi-output neural networks: Insights on multi-task learning and network compression.
\newblock {\em Journal of Machine Learning Research}, 25(231):1--40, 2024.

\bibitem{sidhu1979vector}
G.~S. Sidhu and H.~L. Weinert.
\newblock Vector-valued lg-splines i. interpolating splines.
\newblock {\em Journal of Mathematical Analysis and Applications}, 70(2):505--529, 1979.

\bibitem{sidhu1984vector}
G.~S. Sidhu and H.~L. Weinert.
\newblock Vector-valued lg-splines ii. smoothing splines.
\newblock {\em Journal of Mathematical Analysis and Applications}, 101(2):380--396, 1984.

\bibitem{sun2000control}
S.~Sun, M.~B. Egerstedt, and C.~F. Martin.
\newblock Control theoretic smoothing splines.
\newblock {\em IEEE Transactions on automatic control}, 45(12):2271--2279, 2000.

\bibitem{treves2016topological}
F.~Treves.
\newblock {\em Topological Vector Spaces, Distributions and Kernels: Pure and Applied Mathematics}.
\newblock Elsevier, 2016.

\bibitem{tropp2006algorithms}
J.~A. Tropp.
\newblock Algorithms for simultaneous sparse approximation. part ii: Convex relaxation.
\newblock {\em Signal Processing}, 86(3):589--602, 2006.

\bibitem{unser2021unifying}
M.~Unser.
\newblock A unifying representer theorem for inverse problems and machine learning.
\newblock {\em Foundations of Computational Mathematics}, 21(4):941--960, 2021.

\bibitem{unser2022convex}
M.~Unser and S.~Aziznejad.
\newblock Convex optimization in sums of banach spaces.
\newblock {\em Applied and Computational Harmonic Analysis}, 56:1--25, 2022.

\bibitem{unser2005cardinal}
M.~Unser and T.~Blu.
\newblock Cardinal exponential splines: Part i-theory and filtering algorithms.
\newblock {\em IEEE Transactions on Signal Processing}, 53(4):1425--1438, 2005.

\bibitem{unser2017splines}
M.~Unser, J.~Fageot, and J.~P. Ward.
\newblock Splines are universal solutions of linear inverse problems with generalized tv regularization.
\newblock {\em SIAM Review}, 59(4):769--793, 2017.

\bibitem{wahba1990spline}
G.~Wahba.
\newblock {\em Spline Models for Observational Data}.
\newblock SIAM, 1990.

\bibitem{weinert1978statistical}
H.~L. Weinert.
\newblock Statistical methods in optimal curve fitting.
\newblock {\em Communications in Statistics-Simulation and Computation}, 7(4):417--435, 1978.

\bibitem{werner1984extreme}
D.~Werner.
\newblock Extreme points in spaces of operators and vector--valued measures.
\newblock {\em Proceedings of the 12th Winter School on Abstract Analysis}, pages 135--143, 1984.

\bibitem{willems1997introduction}
J.~C. Willems and J.~W. Polderman.
\newblock {\em Introduction to Mathematical Systems Theory: A Behavioral Approach}, volume~26.
\newblock Springer Science \& Business Media, 1997.

\bibitem{yuan2006model}
M.~Yuan and Y.~Lin.
\newblock Model selection and estimation in regression with grouped variables.
\newblock {\em Journal of the Royal Statistical Society Series B: Statistical Methodology}, 68(1):49--67, 2006.

\bibitem{zhang1997splines}
Z.~Zhang, J.~Tomlinson, and C.~Martin.
\newblock Splines and linear control theory.
\newblock {\em Acta Applicandae Mathematica}, 49:1--34, 1997.

\end{thebibliography}

\appendix
\section{Appendix}
\subsection{Proofs of Section 2.1}
\label{app:2.1}

\begin{proof}[\textbf{Proof of Proposition} \ref{prop:1}]
Item 1 follows from \cite{unser2022convex}[Lemma 1 Item 1] and the Riesz-Markov-Kakutani theorem. Item 2: Since all norms on $\R^D$ are equivalent, there exist constants $C_1,C_2>0$ such that 
\begin{equation}
        \forall z\in\R^D:\quad C_1\norm{z}_{A'}\leq\norm{z}_{B'}\leq C_2\norm{z}_{A'}.
\end{equation}
It follows that 
\begin{align}
        \forall \bm{\mu}\in\M(\R,\R^D):\quad &C_1\norm{\norm{\bm{\mu}}_\M}_{A'}\leq\norm{\norm{\bm{\mu}}_\M}_{B'}\leq C_2\norm{\norm{\bm{\mu}}_\M}_{A'}\\
        &\Rightarrow C_1\norm{\bm{\mu}}_{A',\M}\leq\norm{\bm{\mu}}_{B',\M}\leq C_2\norm{\bm{\mu}}_{A',\M}.
\end{align}
The extreme-point characterization in Item 3 follows from \cite{unser2022convex}[Lemma 1 Item 3], which uses the monotone norm assumption. Then, the norm calculation is a  consequence of the fact that, if $\mu=\sum_{k=1}^{K}a_k\delta_{x_k}\in\M(\R)$ with $x_k\neq x_{k'}$ for $k\neq k'$, then $\norm{\mu}_\M=\sum_{k=1}^K\vert a_k\vert$.
\end{proof}

\subsection{Proofs of Section 2.2}
\label{app:2.2}

\begin{proof}[\textbf{Proof of Proposition} \ref{prop:3}]
    \hfill\\
    \textbf{Item 1.} We know from Proposition \ref{prop:2} that
    $\norm{\cdot}_{\M,1}$ is the dual norm of
    \begin{equation}
    \forall \bb{g}\in\C_0(\R,\R^D):\underset{x\in\R}{\text{sup}}\norm{\bb{g}(x)}_{\infty}=\underset{x\in\R}{\text{sup}}\underset{1\leq d\leq D}{\text{sup}}\vert g_d(x)\vert=\underset{x\in\R,1\leq d\leq D}{\text{sup}}\vert g_d(x)\vert.    
    \end{equation}
    We also know from Item 1 of Proposition \ref{prop:1} that $\norm{\cdot}_{1,\M}$ is the dual norm of
    \begin{equation}
    \forall \bb{g}\in\C_0(\R,\R^D):\underset{1\leq d\leq D}{\text{sup}}\underset{x\in\R}{\text{sup}}\vert g_d(x)\vert=\underset{x\in\R,1\leq d\leq D}{\text{sup}}\vert g_d(x)\vert.    
    \end{equation}
    Since $\norm{\cdot}_{\M,1}$ and $\norm{\cdot}_{1,\M}$ are dual norms of the same norm, they are equal.
    \hfill\\\hfill\\
    \textbf{Item 2.} Since all norms on $\R^D$ are equivalent, there exist constants $C_1,C_2>0$ such that 
    \begin{equation}
        \forall z\in\R^D:\quad C_1\norm{z}_{A'}\leq\norm{z}_{B'}\leq C_2\norm{z}_{A'}.
    \end{equation}
    In addition, both the summation and the supremum operations, applied on positive quantities, preserve the ordering. It follows from the measure-theoretic definition of the norm $\norm{\cdot}_{\Zz',\M}$ that 
    \begin{equation}
      C_1\underset{\pi}{\text{sup}}\sum_{E\in\pi}\norm{\bm{\mu}(E)}_{A'}\leq\underset{\pi}{\text{sup}}\sum_{E\in\pi}\norm{\bm{\mu}(E)}_{B'}\leq C_2\underset{\pi}{\text{sup}}\sum_{E\in\pi}\norm{\bm{\mu}(E)}_{A'},  
    \end{equation}
    which yields to the desired result:
    \begin{equation}
      \forall\bm{\mu}\in\M(\R,\R^D):\quad C_1\norm{\bm{\mu}}_{\M,A'}\leq\norm{\bm{\mu}}_{\M,B'}\leq C_2\norm{\bm{\mu}}_{\M,A'}.  
    \end{equation}
        \hfill\\
    \textbf{Item 3.} The extreme-point characterization is a reformulation of Theorem 2 of \cite{werner1984extreme}. As for the norm calculation, on one hand we have that 
    \begin{equation}
\norm{\bm{\mu}}_{\Zz',\M}\leq\sum_{k=1}^K\norm{\bb{a}_k\delta_{x_k}}_{\Zz',\M}=\sum_{k=1}^K\norm{\bb{a}_k}_{\Zz'}.
    \end{equation}
    On the other hand, using the assumption that $x_k\neq x_{k'}$ for all $k'\neq k$, we find a partition $\{P_k\}_{k=1}^K$ of $\R$ such that $x_k\in P_k$ and $x_{k'}\not\in P_k$ for al $k'\neq k$. It follows that
    \begin{align}
       \norm{\sum_{k=1}^K\bb{a}_k\delta_{x_k}}_{\Zz',\M}=\underset{\pi}{\text{sup}}\sum_{E\in\pi}\norm{\sum_{k=1}^K\bb{a}_k\delta_{x_k}(E)}_{\Zz'}&\geq\sum_{k'=1}^K\norm{\sum_{k=1}^K\bb{a}_k\delta_{x_k}(P_{k'})}_{\Zz'}\\
       &=\sum_{k'=1}^K\norm{\bb{a}_{k'}\delta_{x_{k'}}(P_{k'})}_{\Zz'}=\sum_{k'=1}^K\norm{\bb{a}_{k'}}_{\Zz'}.
    \end{align}
\end{proof}

\subsection{Proofs of Section 3.2}
\label{app:3.2}
Observe that the action of an MDO $\Ls$ is the same as the convolution with $\Ls\{\bb{I}\delta\}$. For example, in one dimension, the action of the differential operator $\D$ is the same as the convolution with $\delta^{'}$. In the sequel, we use this fact and convolve with $\Ls$ instead of $\Ls\{\bb{I}\delta\}$, with a slight abuse of notation. We also make use, without any explicit statement, of the properties showed in Appendix \ref{app:B} of the convolution between matrix-valued distributions. We notate the convolution (composition) of $D$ MDO as $\circledast_{d=1}^D\Ls_d$

In Proposition \ref{prop:Greens}, we use the Laplace formula for inverting a matrix to construct the Green's matrix of an MDO. This construction, different from the one of the Smith normal form given in the core of this paper, will be useful later on. We denote by $\Ls_{r,c}$ the MDO found by removing the $r$th row and the $c$th column of $\Ls$.
\begin{ Proposition }
\label{prop:Greens}

    Let $\mathbf{com}(\Ls)$ be such that $[\mathbf{com}(\Ls)]_{r,c}=(-1)^{r+c}\mathrm{det}(\Ls_{r,c})$. Then, the following equality holds:
    \begin{equation}
    \label{eq:A.3.83}
    \Ls\ast\mathbf{com}(\Ls)^{\intercal}=\bb{I}\mathrm{det}(\Ls).
    \end{equation}
    If, in addition, $\mathrm{det}(\Ls)\neq0$, then it holds that
    \begin{equation}
    \label{eq:A.3.84}
    \Ls\ast\mathbf{com}(\Ls)^{\intercal}\ast\bb{I}g_{\mathrm{det}(\Ls)}=\bb{I}\delta,
    \end{equation}
    with $g_{\mathrm{det}(\Ls)}$ being the unique causal Green's function of $\mathrm{det}(\Ls)$. It follows that $\bb{G}_{\Ls}=\mathbf{com}(\Ls)^{\intercal}\ast \bb{I}g_{\mathrm{det}(\Ls)}$ is the causal Green's matrix of $\Ls$.
\end{ Proposition }
\begin{proof}[\textbf{Proof  of Proposition \ref{prop:Greens}}]
    For matrices (here $D\times D$) with coefficients in a commutative ring (here, the ring of ordinary differential operators), it is known that the equality $\Ls\ast\mathbf{com}(\Ls)^{\intercal}=\bb{I}\mathrm{det}(\Ls)$ always holds. This proves \eqref{eq:A.3.83}. Concerning \eqref{eq:A.3.84}, we calculate that 
    \begin{equation}
        \Ls\ast(\mathbf{com}(\Ls)^{\intercal}\ast\bb{I}g_{\mathrm{det}(\Ls)})=(\Ls\ast\mathbf{com}(\Ls)^{\intercal})\ast\bb{I}g_{\mathrm{det}(\Ls)}=\bb{I}\mathrm{det}(\Ls)\ast\bb{I}g_{\mathrm{det}(\Ls)}=\bb{I}\delta,
    \end{equation}
where we used the associativity of the convolution of matrix-valued distributions (Proposition \ref{prop:10}).
\end{proof}
\begin{proof}[\textbf{Proof of Theorem \ref{prop:greens2}}]
By inserting the factorization $\bb{U}\ast\bb{D}\ast\bb{V}$ of $\Ls$, we get
\begin{align}
\Ls\ast\bb{V}^{-1}\ast\bb{D}^{-1}\ast\bb{U}^{-1}\ast\bb{I}\delta&=\bb{U}\ast\bb{D}\ast\bb{V}\ast\bb{V}^{-1}\ast\bb{D}^{-1}\ast\bb{U}^{-1}\ast\bb{I}\delta\\
&=\bb{U}\ast\bb{D}\ast\bb{D}^{-1}\ast\bb{U}^{-1}\ast\bb{I}\delta=\bb{U}\ast\bb{U}^{-1}\ast\bb{I}\delta=\bb{I}\delta.
\end{align}
Hence, $\bb{G}_{\Ls}=\bb{V}^{-1}\ast\bb{D}^{-1}\ast\bb{U}^{-1}\ast\bb{I}\delta$ is a Green's matrix of $\Ls$. Since the convolution kernels $g_{\D_d}$ defining $\bb{D}^{-1}$ are all causal by definition, the distribution $\bb{G}_{\Ls}$ is also causal.

To show by contradiction that $\bb{G}_{\Ls}$ is the unique causal Green's matrix, we assume the existence of another causal Green's function $\bb{G}_{\Ls}'$. Consequently, $\Ls\ast(\bb{G}_{\Ls}-\bb{G}_{\Ls}')=\bm{0}$ and $[\bb{G}_{\Ls}-\bb{G}_{\Ls}']_{\cdot,d}\in\mathcal{N}_{\Ls}.$ From Proposition \ref{prop:kernel}, we find that $[\bb{G}_{\Ls}-\bb{G}_{\Ls}']_{\cdot,d}\in\bb{V}^{-1}\mathcal{N}_{\bb{D}}$. Since all functions in$\bb{V}^{-1}\mathcal{N}_{\bb{D}}$ are analytic and, by causality, $[\bb{G}_{\Ls}-\bb{G}_{\Ls}']_{\cdot,d}$ is 0 on $\R^-$, the only possibility is that $(\bb{G}_{\Ls}-\bb{G}_{\Ls}')=\bm{0}$. 

To establish \eqref{eq:3.2.37}, we calculate that, for $\bm{\psi}=(\psi_{d})_{d=1}^D\in\Dd(\R,\R^D)$ and $\bm{\phi}=(\phi_{d})_{d=1}^D\in\Dd(\R,\R^D)$,
\begin{align}
    \langle\bb{G}_{\Ls}\ast\bm{\phi},\bm{\psi}\rangle&=\sum_{d=1}^D\int_{\R}\psi_{d}(t)\int_{\R}\sum_{d'=1}^D[\bb{G}_{\Ls}]_{d,d'}(t-\tau)\phi_{d'}(\tau)\mathrm{d}\tau\mathrm{d}t\\
    &=\sum_{d'=1}^D\int_{\R}\phi_{d'}(\tau)\int_{\R}\sum_{d=1}^D[\bb{G}_{\Ls}]_{d,d'}(t-\tau)\psi_{d}(t)\mathrm{d}t\mathrm{d}\tau\\
    &=\sum_{d'=1}^D\int_{\R}\phi_{d'}(\tau)\int_{\R}\sum_{d=1}^D[\bb{G}_{\Ls}^{\intercal,\vee}]_{d',d}(\tau-t)\psi_{d}(t)\mathrm{d}t\mathrm{d}\tau=\langle\bm{\phi},\bb{G}_{\Ls}^{\intercal,\vee}\ast\bm{\psi}\rangle.
\end{align}
To get \eqref{eq:3.2.38}, we first calculate that 
\begin{align}
    \bb{G}_{\Ls}^{\star}=(\bb{V}^{-1}\ast\bb{D}^{-1}\ast\bb{U}^{-1}\ast\bb{I}\delta)^{\star}&=\bb{I}\delta^{\star}\ast\bb{U}^{-1,\star}\ast\bb{D}^{-1,\star}\ast\bb{V}^{-1,\star}\\
    &=\bb{U}^{-1,\star}\ast\bb{D}^{-1,\star}\ast\bb{V}^{-1,\star}\ast\bb{I}\delta,
\end{align}
where we used the fact that, since $\delta$ is the neutral element for the convolution, it is self-adjoint. Moreover $\bb{I}\delta$ commutes with all matrices because it has a constant diagonal. Then, we can write that 
\begin{align}
\Ls^{\star}\ast\bb{G}_{\Ls}^{\star}&=\bb{V}^{\star}\ast\bb{D}^{\star}\ast\bb{U}^{\star}\ast\bb{U}^{-1,\star}\ast\bb{D}^{-1,\star}\ast\bb{V}^{-1,\star}\ast\bb{I}\delta\\
&=\bb{V}^{\star}\ast\bb{D}^{\star}\ast\bb{U}^{\star}\ast\bb{U}^{\star,-1}\ast\bb{D}^{-1,\star}\ast\bb{V}^{\star,-1}\ast\bb{I}\delta\\
&=\bb{V}^{\star}\ast\bb{D}^{\star}\ast\bb{D}^{-1,\star}\ast\bb{V}^{\star,-1}\ast\bb{I}\delta=\bb{V}^{\star}\ast\bb{I}\delta\ast\bb{V}^{\star,-1}\ast\bb{I}\delta=\bb{I}\delta
\end{align}
because $\bb{U}^{-1,\star}=\bb{U}^{\star,-1}$ and $\bb{D}^{\star}\ast\bb{D}^{-1,\star}=\bb{I}\delta$. The latter comes from the equality $\D_d^{\star}\ast g_{\D_d}^{\star}=\D_d^{\star}\ast g_{\D_d}^{\vee}=\delta$. This concludes the proof.
\end{proof}\hfill\\

\begin{proof}[\textbf{Proof of Proposition \ref{prop:reg}}]
 Denote the order of an ODO $\ls$ by $O(\ls).$ From Proposition \ref{prop:Greens}, we find that 
    \begin{equation}
        [\bb{G}_{\Ls}]_{r,c}=[\mathbf{com}(\Ls)^{\intercal}]_{r,c}\ast g_{\mathrm{det}(\Ls)}=(-1)^{r+c}\mathrm{det}(\Ls_{c,r})\ast g_{\mathrm{det}(\Ls)}.
    \end{equation}
 Furthermore, we find by summing over all permutations $\tilde{\sigma}:[1\cdots D]\setminus\{c\}\to[1\cdots D]\setminus\{r\}$ that,
    \begin{align}
    \mathrm{det}(\Ls_{c,r})=\sum_{\tilde{\sigma}}\mathrm{sgn}(\tilde{\sigma})\circledast_{d=1,d\neq c}^D\ls_{d,\tilde{\sigma}(d)}
    \Rightarrow O(\mathrm{det}(\Ls_{c,r}))\leq\underset{\tilde{\sigma}}{\mathrm{max}}\sum_{d=1,d\neq c}^DO(\ls_{d,\tilde{\sigma}(d)}).\label{eq:A.3.97}
    \end{align}
     Consider a specific permutation $\hat{\sigma}$ for which the maximum in \eqref{eq:A.3.97} is achieved. Then, we calculate that

     \begin{align}
    N_{c,r}+ O(\mathrm{det}(\Ls_{c,r}))&\leq N_{c,r}+\underset{\tilde{\sigma}}{\mathrm{max}}\sum_{d=1,d\neq c}^DO(\ls_{d,\tilde{\sigma}(d)})=O(\ls_{c,r})+\sum_{d=1, d\neq c}^DO(\ls_{d,\hat{\sigma}(d)})\\
    &\leq\underset{\sigma}{\mathrm{max }}\left(O(\ls_{c,\sigma(c)})+\sum_{d=1, d\neq c}^DO(\ls_{d,\sigma(d)})\right)\\
    &\leq\underset{\sigma}{\mathrm{max}}\sum_{d=1}^DO(\ls_{d,\sigma(d)})=O(\mathrm{det}(\Ls)).\label{eq:A.3.99}
     \end{align}
The last equality in \eqref{eq:A.3.99} holds because
\begin{equation}
    O(\mathrm{det}(\Ls))=O\left(\sum_{\sigma}\mathrm{sgn}(\tilde{\sigma})\circledast_{d=1}^D\ls_{d,\sigma(d)}\right)=\sum_{d=1}^DO(\ls_{d,d})=\underset{\sigma}{\mathrm{max}}\sum_{d=1}^DO(\ls_{d,\sigma(d)}).\label{eq:A.3.100}
\end{equation}
In turn, the last two equalities in \eqref{eq:A.3.100} hold because of the assumption of the proposition. Hence, from \eqref{eq:A.3.99}, if $N_{c,r}\geq2$, then $g_{\mathrm{det}(\Ls)}\in\mathcal{C}^{N_{c,r}+O(\mathrm{det}(\Ls_{c,r}))-2}(\R)$ and $[\bb{G}_{\Ls}]_{r,c}\in\mathcal{C}^{N_{c,r}-2}(\R).$ We used the fact that the Green's function of an ordinary differential operator of order $N\geq2$ always belong in $\mathcal{C}^{N-2}(\R).$ If $N_{c,r}=1,$ the same analysis is made but $[\bb{G}_{\Ls}]_{r,c}$ is reduced (in term of regularity) to the derivative of a locally Lipschitz-continuous function, which is known to be in the function space $\mathcal{L}^{\infty}_{loc}(\R).$
\end{proof}

\subsection{Proofs of Section 3.3}
\label{app:3.3}
This section is richer than the concatenation of the proofs of Section 3.3. We start by providing the proof of Proposition \ref{prop:testfunctions}. Then, we introduce in Proposition \ref{prop:compact} a compactly supported right inverse $\bb{G}_{\bm{\phi}}$ of $\Ls$, which allows us to move from $\M(\R,\R^D)$ to $\M_{\Ls}(\R,\R^D)$. In Proposition \ref{prop:defmap} and Theorem \ref{th:M} we further study the action of the kernel $\bb{G}_{\bm{\phi}}$ on the measure space $\M(\R,\R^D)$. Finally, leveraging the previous results, we provide the proof of Theorem \ref{prop:MLprop}.

Observe that even-though Theorem \ref{th:M} is stated before Theorem \ref{prop:MLprop}, the latter does not use Theorem \ref{th:M} in its proof. Consequently, we are not making a "circular" argumentation there. 

\begin{proof}[\textbf{Proof of Proposition \ref{prop:testfunctions}}]
Consider the function $\zeta:\R\to\R$ expressed as
\begin{equation}
    \zeta(x)=\begin{cases}
        c\mathrm{e}^{\frac{1}{x^2-1}},&x\in]-1,1[\\
        0, & \text{else},
    \end{cases}
\end{equation}
 where $c^{-1}=\int_{-1}^{1}\mathrm{e}^{\frac{-1}{1-x^2}}\mathrm{d}x$. The function $\zeta$ is compactly supported in $[-1,1]$, belongs to $\mathcal{L}_1(\R)$ with $\norm{\zeta}_{\mathcal{L}_1}=1$, and is also in $\mathcal{C}^{\infty}(\R).$ In particular, $\zeta\in\mathcal{D}(\R)$. Then, $\forall f\in\C(\R)$, using the continuity of $f$ we get \eqref{eq:A.4.99}, and $\forall t\in\R $ we get \eqref{eq:A.4.100}.
\begin{align}
   &\underset{R\to0}{\text{lim }}\left\langle\frac{1}{R}\zeta\left(\frac{\cdot}{R}\right),f\right\rangle-f(0)=\underset{R\to0}{\text{lim }}\int_{-R}^R\frac{c}{R}\mathrm{e}^{\frac{-1}{1- \frac{x^2}{R^2}}}(f(x)-f(0))\mathrm{d}x\\
    &\hspace{+4.1cm}=\underset{R\to0}{\text{lim }}\int_{-1}^1c\mathrm{e}^{\frac{-1}{1- x^2}}(f(Rx)-f(0))\mathrm{d}x=0\label{eq:A.4.99}.\\
   & \underset{R\to0}{\text{lim }}\left\langle\frac{1}{R}\zeta\left(\frac{\cdot}{R}-t\right),f\right\rangle=f(t)\label{eq:A.4.100}.
\end{align}
We take advantage of Lemma \ref{Lemma:mat2} to build the desired family $(\bm{\phi}_n)_{n=1}^N$.

\begin{Lemma} 
\label{Lemma:mat2}
Let $\{\bb{p}_n\}_{n=1}^{N}:[1,\infty[\to\R^D$ be a free (independent) family of analytic functions. Then, there exist two lists $\{\bb{t}_{n'}\}_{n'=1}^N$ and $\{\bb{c}_{n'}\}_{n'=1}^N$ of coordinates $\bb{t}_{n'}\in[1,\infty[^D$ and coefficients $\bb{c}_{n'}\in\R^D$ such that the matrix $\mathbf{A}\in\R^{N\times N}$ defined as $\mathbf{A}_{n',n}=\langle\bb{p}_{n},\bb{c}_{n'}\odot\bm{\delta}_{\bb{t}_{n'}} \rangle=\sum_{d=1}^D[\bb{c}_{n'}]_d[\bb{p}_{n}]_{d}([\bb{t}_{n'}]_d)$ is invertible.
\end{Lemma}

\begin{proof}[\textbf{Proof of Lemma \ref{Lemma:mat2}}]
We only provide a brief proof by induction on $N>1$. If $N=2,$ take $\bb{c}_1\in\R^D$ and $\bb{t}_1\in[1,\infty[^D$ such that $\langle\bb{c}_{1}\odot\bm{\delta}_{\bb{t}_{1}} ,\bb{p}_{1}\rangle\neq0$. They must exist, otherwise $\bb{p}_1=\bm{0}$ on $[1,\infty[$ and, since $\bb{p}_1$ is analytic, $\bb{p}_1=\bm{0}$ on $\R$. For $\bb{c}\in\R^D$ and $\bb{t}\in[1,\infty[^D$, the determinant of the matrix 
\begin{equation}
\begin{bmatrix}
    \langle\bb{p}_{1},\bb{c}_{1}\odot\bm{\delta}_{\bb{t}_{1}} \rangle & \langle\bb{p}_{2},\bb{c}_{1}\odot\bm{\delta}_{\bb{t}_{1}} \rangle\\
    \langle\bb{p}_{1},\bb{c}\odot\bm{\delta}_{\bb{t}}\rangle & \langle\bb{p}_{2},\bb{c}\odot\bm{\delta}_{\bb{t}} \rangle
    \end{bmatrix}
\quad\text{is}\quad  \langle\bb{p}_{1},\bb{c}_{1}\odot\bm{\delta}_{\bb{t}_{1}}\rangle\langle\bb{p}_{2},\bb{c}\odot\bm{\delta}_{\bb{t}}\rangle-\langle\bb{p}_{2},\bb{c}_{1}\odot\bm{\delta}_{\bb{t}_{1}} \rangle \langle\bb{p}_{1},\bb{c}\odot\bm{\delta}_{\bb{t}}\rangle.
\end{equation}
Assume by contradiction that there exists no $\bb{c}\in\R^D$ and $\bb{t}\in[1,\infty[^D$ for which the matrix is invertible. Therefore,  $\forall\bb{c}\in\R^D$ and $\forall\bb{t}\in[1,\infty[^D$, its determinant is 0 and
\begin{equation}
\label{eq:3.3.52}
 \langle\bb{p}_{2},\bb{c}\odot\bm{\delta}_{\bb{t}} \rangle=\frac{\langle\bb{p}_{2},\bb{c}_{1}\odot\bm{\delta}_{\bb{t}_{1}}\rangle }{\langle\bb{p}_{1},\bb{c}_{1}\odot\bm{\delta}_{\bb{t}_{1}}\rangle }\langle\bb{p}_{1},\bb{c}\odot\bm{\delta}_{\bb{t}}\rangle.
\end{equation}
It follows from \eqref{eq:3.3.52} that $\bb{p}_1$ and $\bb{p}_2$ are linearly dependent on $[1,\infty[$. Because $\{\bb{p}_{1},\bb{p}_{2}\}$ are analytic, they must be linearly dependent on $\R$, which is a contradiction with the initial assumption of them being a free family. The induction step $N$ reaches the same contradiction over the general formula of a determinant and the assumption that an invertible matrix exists for the $(N-1)$ case. This concludes the proof of Lemma \eqref{Lemma:mat2}.
\end{proof}

 We now resume the proof of Proposition \eqref{prop:testfunctions}. Lemma \ref{Lemma:mat2} yields an invertible matrix $\bb{A}\in\R^{N\times N}$ such that 
\begin{equation}
      \bb{A}_{n',n}=\left\langle \bb{p}_n,\bb{c}_{n'}\odot\bm{\delta}_{\bb{t}_{n'}}\right\rangle\quad\text{with}\quad\bb{t}_{n'}\in[1,\infty[^D. 
\end{equation}
Define $\forall n'\in[1\cdots N]$, $\bm{\zeta}_{n'}(R)=\bb{c}_{n'}\odot\left(\frac{1}{R}\zeta\left(\frac{\cdot}{R}-[\bb{t}_{n'}]_d\right)\right)_{d=1}^D\in\mathcal{D}(\R,\R^D)$ and $\bm{\zeta}(R)=(\bm{\zeta}_{n'}(R))_{n'=1}^N\in\mathcal{D}(\R,\R^D)^N$. Define the matrix $\bb{A}(R)\in\R^{N\times N}$ such that
\begin{equation}
  \bb{A}_{n',n}(R)=\left\langle \bb{p}_n,\bm{\zeta}_{n'}(R)\right\rangle . 
\end{equation}
Since $\mathrm{det}(\cdot)$ is a continuous function in terms of the entries of the matrix it receives, we find that $\underset{R\to 0}{\text{lim}}\mathrm{det}(\bb{A}(R))=\mathrm{det}(\bb{A})\neq0$, and there exists an $R'$ such that $\mathrm{det}(\bb{A}(R'))\neq0.$  Define $\bm{\phi}=(\bm{\phi}_{n'})_{n'=1}^N=\bb{A}(R')^{-1}\bm{\zeta}(R')\in\mathcal{D}(\R,\R^D)^N$. The family $\bm{\phi}$ has the expected properties. First, it verifies Item 2 as a linear combination of $\mathcal{D}(\R,\R^D)$ functions. Second, it verifies Item 3 as a linear combination of functions $\frac{1}{R'}\zeta\left(\frac{\cdot}{R'}-[\bb{t}_{n'}]_d\right)$ whose support is in $[R'(1-[\bb{t}_{n'}]_d),R'(1+[\bb{t}_{n'}]_d)]\subset[0,\infty[$. Finally, concerning Item 1, we calculate that 
\begin{align}
\langle\bb{p}_{n},\bm{\phi}_{n'}\rangle&=\left\langle\bb{p}_{n},\sum_{m=1}^{N}\bb{A}(R')^{-1}_{n',m}\bm{\zeta}_m(R')\right\rangle=\sum_{m=1}^{N}\bb{A}(R')^{-1}_{n',m}\left\langle\bb{p}_{n},\bm{\zeta}_m(R')\right\rangle\\
&=\sum_{m=1}^{N}\bb{A}(R')^{-1}_{n',m}\bb{A}(R')_{m,n}=\delta_{n,n'}.
\end{align}

\end{proof}
\begin{ Proposition }
\label{prop:compact}
Let $\left((\bb{p}_n)_{n=1}^N,(\bm{\phi}_n)_{n=1}^N\right)$ be an $\Ls$-admissible system. Then, the kernel $\bb{G}_{\bm{\phi}}(\cdot,\cdot)$ defined as 
    \begin{equation}
    \bb{G}_{\bm{\phi}}(x,y)=\bb{G}_{\Ls}(x-y)-\sum_{n=1}^{N}\bb{p}_n(x)\left(\bb{G}_{\Ls}^{\star}\ast\bm{\phi}_n\right)^{\intercal}(y)
    \end{equation}
    is such that 
    \begin{align}
        &\forall x\in\R,\forall y\notin[\mathrm{min}(0,x),\mathrm{max}(x,\phi^+)]:&&\quad\bb{G}_{\bm{\phi}}(x,y)=\bm{0},\label{eq:A.4.107}\\
        &\begin{cases}
            \forall y\in]-\infty,0[,\forall x\in]y,\infty[\\
            \forall y\in]\phi^+,\infty[,\forall x\in]-\infty,y[, 
        \end{cases}:&&\quad\bb{G}_{\bm{\phi}}(x,y)=\bm{0}\label{eq:A.4.108}.
    \end{align}
\end{ Proposition }
\begin{proof}[\textbf{Proof of Proposition \ref{prop:compact}}]
 We recall that
\begin{equation}
    \Ls=\bb{U}\ast\bb{D}\ast\bb{V}\quad\text{with}\quad\bb{D}=\mathbf{diag}(\D_1,...,\D_D),
\end{equation}
where $\bb{U}$ and $\bb{V}$ are unimodular. As usual, we denote the causal Green's function of $\D_d$ by $g_{\D_d}$ and 
\begin{equation}
\bb{G}_{\Ls}=\bb{V}^{-1}\ast\bb{D}^{-1}\ast\bb{U}^{-1},\quad\mathrm{with}\quad\bb{D}^{-1}=\mathbf{diag}(g_{\D_1},...,g_{\D_D}).    
\end{equation}
 Observe that, if $\left((\bb{p}_n)_{n=1}^N,(\bm{\phi}_n)_{n=1}^N\right)$ is an $\Ls$-admissible system, then the system
 
 \begin{equation}
     \left((\tilde{\bb{p}}_n)_{n=1}^N,(\tilde{\bm{\phi}}_n)_{n=1}^N\right)=\left((\bb{V}\ast\bb{p}_n)_{n=1}^N,(\bb{V}^{-1,\star}\ast\bm{\phi}_n)_{n=1}^N\right)
 \end{equation}
is $\bb{D}$-admissible, where $(\tilde{\bb{p}}_n)_{n=1}^{N}$ is a basis of $\mathcal{N}_{\bb{D}}$. Consequently, we calculate that 
\begin{align}
        \bb{G}_{\bm{\phi}}(x,y)&=\left[\bb{V}^{-1}\ast\bb{D}^{-1}\ast\bb{U}^{-1}\ast\bb{I}\delta_{y}\right](x)-\sum_{n=1}^{N}\bb{p}_n(\cdot)
        \left[(\bb{V}^{-1}\ast\bb{D}^{-1}\ast\bb{U}^{-1}\ast\bb{I}\delta)^{\star}\ast\bm{\phi}_n\right]^{\intercal}(y)\label{eq:A.4.115}
        \\
        &=\left[\bb{V}^{-1}\ast\bb{D}^{-1}\ast\bb{U}^{-1}\ast\bb{I}\delta_{y}\right](x)-\sum_{n=1}^{N}\bb{p}_n(x)
        \left\langle\bb{V}^{-1}\ast\bb{D}^{-1}\ast\bb{U}^{-1}\ast\bb{I}\delta_y,\bm{\phi}_n\right\rangle^{\intercal}\label{eq:A.4.116}
        \\
        &=\left[\bb{V}^{-1}\ast\bb{D}^{-1}\ast\bb{I}\delta_{y}\ast\bb{U}^{-1}\right](x)-\sum_{n=1}^{N}[\bb{V}^{-1}\ast\tilde{\bb{p}}_n](x)\left\langle\bb{D}^{-1}\ast\bb{I}\delta_{y}\ast\bb{U}^{-1},\tilde{\bm{\phi}}_n\right\rangle^{\intercal}\label{eq:A.4.112},
\end{align}
where to go from \eqref{eq:A.4.115} to \eqref{eq:A.4.116} we used the Bracket property (see Proposition \ref{prop:10}) of convolution. Recall that $\bb{D}^{-1}$ is the diagonal concatenation of the Green's functions of the operators $\D_d$ forming the diagonal of $\bb{D}$. In particular, $\bb{D}^{-1}$ is supported in $[0,\infty[.$ Consequently, $\text{supp}\left(\bb{D}^{-1}\ast\bb{I}\delta_y\right)=[y,\infty[$ and since $\bb{U}^{-1}$ is an MDO, it preserves the support, and  $\text{supp}\left(\bb{D}^{-1}\ast\bb{I}\delta_y\ast\bb{U}^{-1}\right)=[y,\infty[$. Observe further that $\tilde{\bm{\phi}}_n$ is supported in $[0,\phi^+]$. Consequently, if $y>\mathrm{max}(x,\phi^+)\geq\phi^+$, then 
\begin{equation}
    \text{supp}\left(\bb{D}^{-1}\ast\bb{U}^{-1}\ast\bb{I}\delta_{y}\right)\cap\text{supp}(\tilde{\bm{\phi}}_n)\subset[y,\infty[\cap[0,\phi^+]=\emptyset\Rightarrow\left\langle\bb{D}^{-1}\ast\bb{U}^{-1}\ast\bb{I}\delta_{y},\tilde{\bm{\phi}}_n\right\rangle=\bm{0}.
\end{equation}
Since $\bb{V}^{-1}$ is an MDO, it preserves the support and $\text{supp}\left(\bb{V}^{-1}\ast\bb{D}^{-1}\ast\bb{I}\delta_y\ast\bb{U}^{-1}\right)=[y,\infty[$. Therefore, if $y>\mathrm{max}(x,\phi^+)\geq x$, then $\left[\bb{V}^{-1}\ast\bb{D}^{-1}\ast\bb{I}\delta_{y}\ast\bb{U}^{-1}\right](x)=\bm{0}$. Finally, $\bb{G}_{\bm{\phi}}(x,y)=0$ if $y>\mathrm{max}(x,\phi^+).$ Then, we observe that $\bb{D}^{-1}$ is such that $\bb{D}^{-1}=\bb{P}_+$, where each column of $\bb{P}$ is in $\mathcal{N}_{\bb{D}}$, and $\bb{P}_+$ is the causal part of $\bb{P}.$ If $y<\mathrm{min}(0,x)$, then $y<x$ and $y<0$ imply that 
\begin{align}
   &\left[\bb{V}^{-1}\ast\bb{D}^{-1}\ast\bb{I}\delta_{y}\ast\bb{U}^{-1}\right](x)=\left[\bb{V}^{-1}\ast\bb{P}\ast\bb{I}\delta_{y}\ast\bb{U}^{-1}\right](x),\\
   \text{and},\quad &\left\langle\bb{D}^{-1}\ast\bb{I}\delta_{y}\ast\bb{U}^{-1},\tilde{\bm{\phi}}_n\right\rangle=\left\langle\bb{P}\ast\bb{I}\delta_{y}\ast\bb{U}^{-1},\tilde{\bm{\phi}}_n\right\rangle.\label{eq:A.4.114}
\end{align}
Furthermore, $\mathcal{N}_{\bb{D}}$ is translation-invariant and stable under differentiation. Therefore, each column of $\tilde{\bb{P}}=\bb{P}\ast\bb{I}\delta_y\ast\bb{U}^{-1}$ is in $\mathcal{N}_{\bb{D}}$ and $\tilde{\bb{P}}(\cdot)=\sum_{n=1}^{N}\tilde{\bb{p}}_n(\cdot)\left\langle\tilde{\bb{P}},\tilde{\bm{\phi}}_n\right\rangle^{\intercal}$. Finally, the combination of \eqref{eq:A.4.112}, \eqref{eq:A.4.114} with the previous comment on $\mathcal{N}_{\bb{D}}$, yields that 
\begin{align}
        \bb{G}_{\bm{\phi}}(x,y)&=\left[\bb{V}^{-1}\ast\bb{D}^{-1}\ast\bb{I}\delta_{y}\ast\bb{U}^{-1}\right](x)-\sum_{n=1}^{N}[\bb{V}^{-1}\ast\tilde{\bb{p}}_n](x)\left\langle\bb{D}^{-1}\ast\bb{I}\delta_{y}\ast\bb{U}^{-1},\tilde{\bm{\phi}}_n\right\rangle^{\intercal}\\
        &=\left[\bb{V}^{-1}\ast\tilde{\bb{P}}\right](x)-\sum_{n=1}^{N}[\bb{V}^{-1}\ast\tilde{\bb{p}}_n](x)\left\langle\tilde{\bb{P}},\tilde{\bm{\phi}}_n\right\rangle^{\intercal}\\
        &=\left[\bb{V}^{-1}\ast\tilde{\bb{P}}\right](x)-\left[\bb{V}^{-1}\ast\tilde{\bb{P}}\right](x)=\bm{0}.    
\end{align}
This proves \eqref{eq:A.4.107}. In turn, \eqref{eq:A.4.108} is the converse of \eqref{eq:A.4.107}.
\end{proof}

We define the functional $\bb{G}_{\bm{\phi}}:\Dd(\R,\R^D)\to\C^{\infty}(\R,\R^D)$ as 
\begin{equation}
    \forall\bm{\psi}\in\Dd(\R,\R^D):\quad\bb{G}_{\bm{\phi}}\{\bm{\psi}\}(\cdot)=\int_{\R}\bb{G}_{\bm{\phi}}(\cdot,y)\bm{\psi}(y)\mathrm{d}y=(\bb{G}_{\Ls}\ast\bm{\psi})(\cdot)-\sum_{n=1}^N\bb{p}_n(\cdot)\langle\bb{G}_{\Ls}^{\star}\ast\bm{\phi}_n,\bm{\psi}\rangle.
\end{equation}
Likewise, we define the functional $\bb{G}_{\bm{\phi}}^{\star}:\Dd(\R,\R^D)\to\Dd(\R,\R^D)$ as 
\begin{equation}
    \forall\bm{\psi}\in\Dd(\R,\R^D):\quad\bb{G}_{\bm{\phi}}^{\star}\{\bm{\psi}\}(\cdot)=\int_{\R}\bb{G}_{\bm{\phi}}(x,\cdot)^{\intercal}\bm{\psi}(x)\mathrm{d}x=(\bb{G}_{\Ls}^{\star}\ast\bm{\psi})(\cdot)-\sum_{n=1}^N(\bb{G}_{\Ls}^{\star}\ast\bm{\phi}_n)(\cdot)\langle\bb{p}_n,\bm{\psi}\rangle.
\end{equation}
We observe that $\forall(\bm{\psi},\bm{\rho})\in\Dd(\R,\R^D)^2$, it holds that 
\begin{equation}
    \langle \bb{G}_{\bm{\phi}}^{\star}\{\bm{\psi}\},\bm{\rho}\rangle= \langle \bm{\psi},\bb{G}_{\bm{\phi}}\{\bm{\rho}\}\rangle,\quad\quad\text{and}\quad\langle \bb{G}_{\bm{\phi}}\{\bm{\psi}\},\bm{\rho}\rangle= \langle \bm{\psi},\bb{G}_{\bm{\phi}}^{\star}\{\bm{\rho}\}\rangle.
\end{equation}
In Proposition \ref{prop:defmap} we extend by duality the definition of the mapping $\bb{G}_{\bm{\phi}}$ to the whole vector space $\M(\R,\R^D)$.
\begin{ Proposition }
\label{prop:defmap}
The mapping 
\begin{equation}
        \bb{G}_{\bm{\phi}}:\begin{cases}
        \M(\R,\R^D)\to\bb{G}_{\bm{\phi}}\M(\R,\R^D)\subset\mathcal{D}'(\R,\R^D)\\
        \bb{m}\mapsto\bb{G}_{\bm{\phi}}\{\bb{m}\},
        \end{cases}
    \end{equation}
defined as
\begin{equation}
\label{eq:A.4.123}
    \forall\bm{\psi}\in\mathcal{D}(\R,\R^D):\quad\left\langle\bb{G}_{\bm{\phi}}\{\bb{m}\},\bm{\psi}\right\rangle=\left\langle\bb{m},\bb{G}_{\bm{\phi}}^{\star}\{\bm{\psi}\}\right\rangle,
\end{equation}
is well-defined because $\bb{G}_{\bm{\phi}}^{\star}\{\bm{\psi}\}\in \Dd(\R,\R^D)\subset\C_0(\R,\R^D).$
\end{ Proposition }
\begin{proof}[\textbf{Proof of Proposition \eqref{prop:defmap}}]
    We first calculate that 
    \begin{align}
\int_{\R}\bb{G}_{\bm{\phi}}^{\intercal}(x,\cdot)\bm{\psi}(x)\mathrm{d} x&=\int_{\R}\bb{G}_{\Ls}^{\intercal}(x-\cdot)\bm{\psi}(x)\mathrm{d}x-\sum_{n=1}^{N}(\bb{G}_{\bb{L}}^{\star}\ast\bm{\phi}_n)(\cdot)\int_{\R}\langle\bb{p}_n(x),\bm{\psi}(x)\rangle\mathrm{d}x\\
&=\left[\bb{G}_{\bb{L}}^{\star}\ast\left(\bm{\psi}-\sum_{n=1}^Nc_n\bm{\phi}_n\right)\right](\cdot),\label{eq:3.2.68}
 \end{align}
where $c_n=\int_{\R}\langle\bb{p}_n(x),\bm{\psi}(x)\rangle\mathrm{d}x$. Equation \eqref{eq:3.2.68} is the convolution between a distribution and a test function, which is known to be in $\C^{\infty}(\R,\R^{D})$. Finally, assume that $\bm{\psi}$ is compactly supported in $[\psi^-,\psi^+]$. From \eqref{eq:A.4.108}, we find that, if $y\in]-\infty,\text{min}(0,\psi^{-})[\bigcup]\text{max}(\phi^+,\psi^+),\infty[$, then 
\begin{equation}
    \int_{\R}\bb{G}_{\bm{\phi}}^{\intercal}(x,y)\bm{\psi}(x)\mathrm{d}x=\bm{0}\quad\Rightarrow\quad\int_{\R}\bb{G}_{\bm{\phi}}^{\intercal}(x,\cdot)\bm{\psi}(x)\mathrm{d}x\in\C^{\infty}_c(\R,\R^D)\subset\C_0(\R,\R^D).
\end{equation}
\end{proof}

\begin{theorem}
\label{th:M} 
First, the space $\bb{G}_{\bm{\phi}}\M(\R,\R^D)$ has the following properties.
\begin{itemize}
    \item [1.] For all $\bb{m}\in\M(\R,\R^D):\quad\Ls\{\bb{G}_{\bm{\phi}}\{\bb{m}\}\}=\bb{m}$;
    \item [2.] it is a Banach space for the norm $\norm{\Ls\{\cdot\}}_{\Zz',\M}$ (resp. $\norm{\Ls\{\cdot\}}_{\M,\Zz'}$);
    \item [3.] for all $n$ such that $ 1\leq n\leq N:\quad\left\langle\bb{G}_{\bm{\phi}}\{\bb{m}\},\bm{\phi}_n\right\rangle=0.$
\end{itemize}
Second, the space $\M_{\Ls}(\R,\R^D)$ has the following property:
\begin{itemize}
    \item [4.] direct-sum decomposition $\M_{\Ls}(\R,\R^D)=\bb{G}_{\bm{\phi}}\M(\R,\R^D)\oplus\mathcal{N}_{\Ls}$.
\end{itemize}
\end{theorem}

\begin{proof}[\textbf{Proof of Theorem \ref{th:M}}]
\hfill\\
\textbf{Item 1.}
We first calculate that $\langle\bb{p},\Ls^{\star}\{\bm{\psi}\}\rangle=\langle\Ls\{\bb{p}\},\bm{\psi}\rangle=0,$ $\forall\bm{\psi}\in\Dd(\R,\R^D).$ Observe that $\Ls\{\bb{G}_{\bm{\phi}}\{\bb{m}\}\}\in\Dd'(\R,\R^D).$ Therefore, from \eqref{eq:3.2.68}, we calculate that 
\begin{align}
    \langle\Ls\{\bb{G}_{\bm{\phi}}\{\bb{m}\}\},\bm{\psi}\rangle=\langle\bb{G}_{\bm{\phi}}\{\bb{m}\},\Ls^{\star}\{\bm{\psi}\}\rangle=\langle\bb{m},\bb{G}_{\Ls}^{\star}\ast\Ls^{\star}\{\bm{\psi}\}\rangle=\langle\bb{m},(\Ls\ast\bb{G}_{\Ls})^{\star}\ast\bm{\psi}\rangle=\langle\bb{m},\bm{\psi}\rangle\label{eq:A.4.126}.
\end{align}
Equation \eqref{eq:A.4.126} holds $\forall\bm{\psi}\in\Dd(\R,\R^D),$ which is dense in the predual $\C_0(\R,\R^D)$ of $\M(\R,\R^D)$; thus, it directly implies Item 1.
\hfill\\\hfill\\
\textbf{Item 2.}
It follows directly from Item 1 that $\left(\M(\R,\R^D),\norm{\cdot}_{\Zz',\M}\right)$ and $\left(\bb{G}_{\bm{\phi}}\M(\R,\R^D), \norm{\Ls\{\cdot\}}_{\Zz',\M}\right)$ are isometric.
\hfill\\\hfill\\
\textbf{Item 3.}
For $1\leq n'\leq N$, we get from  \eqref{eq:A.4.123} that
\begin{align}
    \left\langle\bb{G}_{\bm{\phi}}\{\bb{m}\},\bm{\phi}_{n'}\right\rangle&=\left\langle\bb{m},\bb{G}_{\bm{\phi}}^{\star}\{\bm{\phi}_{n'}\}\right\rangle=\left\langle\bb{m},\bb{G}_{\bb{L}}^{\star}\ast\left(\bm{\phi}_{n'}-\sum_{n=1}^N\bm{\phi}_n\langle\bb{p}_n,\bm{\phi}_{n'}\rangle\right)\right\rangle\\
    &=\left\langle\bb{m},\bb{G}_{\bb{L}}^{\star}\ast\left(\bm{\phi}_{n'}-\bm{\phi}_{n'}\right)\right\rangle=0.
\end{align}
\hfill\\\hfill\\
\textbf{Item 4.} We first prove that the sum is direct. Suppose that it is not; then, there exists $\bm{0}\neq\bb{G}_{\bm{\phi}}\{\bb{m}\}\in\bb{G}_{\bm{\phi}}\M(\R,\R^D)$ and $\bm{0}\neq\bb{q}\in\mathcal{N}_{\Ls}$, for which we find the contradiction
\begin{equation}
\bb{G}_{\bm{\phi}}\{\bb{m}\}=\bb{q}\Rightarrow\Ls\{\bb{G}_{\bm{\phi}}\{\bb{m}\}\}=\Ls\{\bb{q}\}\Rightarrow\bb{m}=\bm{0}.
\end{equation}
This shows that $\bb{G}_{\bm{\phi}}\M(\R,\R^D)\oplus\mathcal{N}_{\Ls}\subset\M_{\Ls}(\R,\R^D)$. We conclude by noticing that, $\forall \bb{f}\in\M_{\Ls}(\R,\R^D),$ it holds that
\begin{equation}
    \bb{f}=\bb{G}_{\bm{\phi}}\{\Ls\{\bb{f}\}\}+\left(\bb{f}-\bb{G}_{\bm{\phi}}\{\Ls\{\bb{f}\}\}\right),\quad\text{with}\quad\left(\bb{f}-\bb{G}_{\bm{\phi}}\{\Ls\{\bb{f}\}\}\right)\in\mathcal{N}_{\Ls}.
\end{equation}
\end{proof}

\begin{proof}[\textbf{Proof of Theorem \ref{prop:MLprop}}]
\hfill\\
\textbf{Item 1.}
This follows directly from Theorem \ref{th:M}.
\hfill\\\hfill\\
\textbf{Item 2.}
We know from Theorem \ref{th:M} that $\M_{\Ls}(\R,\R^D)=\bb{G}_{\bm{\phi}}\M(\R,\R^D)\oplus\mathcal{N}_{\Ls}.$ Because $\bb{G}_{\bm{\phi}}\M(\R,\R^D)\oplus\mathcal{N}_{\Ls}$ is isometric to $\M(\R,\R^D)\times\R^N$, the latter space admits a predual, and so does $\M_{\Ls}(\R,\R^D)$. 
\hfill\\\hfill\\
\textbf{Item 3.} This follows directly from Item 2 of Proposition \ref{prop:1} and Items 1 and 2 of Proposition \ref{prop:3}.
\hfill\\\hfill\\
\textbf{Item 4.} We know from Item 4 of Proposition \ref{prop:1} that the extreme points of the unit ball in $(\M(\R,\R^D),\norm{\cdot}_{\M,\Zz'})$ are of the form 
\begin{equation}
\bb{w}\odot\bm{\delta}_{\bb{x}},\quad\text{with}\quad\bb{x}\in\R^D,
\end{equation}
where $\bb{w}$ is an extreme point of the centered unit ball in $\Zz'$. Further, it follows from Theorem \eqref{th:M} that
\begin{equation}
    (\M(\R,\R^D),\norm{\cdot}_{\M,\Zz'}),\quad\text{and}\quad(\bb{G}_{\bm{\phi}}\M(\R,\R^D),\norm{\Ls\{\cdot\}}_{\M,\Zz'})
\end{equation}
are isometric and that, consequently, the extreme points in the unit ball of $(\bb{G}_{\bm{\phi}}\M(\R,\R^D),\norm{\Ls\{\cdot\}}_{\M,\Zz'})$ are exactly of the form
\begin{multline}
\label{eq:148}
    \bb{e}=\bb{G}_{\bm{\phi}}\{\bb{w}\odot\bm{\delta}_{\bb{x}}\}=\bb{G}_{\Ls}\ast(\bb{w}\odot\bm{\delta}_{\bb{x}})-\sum_{n=1}^{N}b_n\bb{p}_n,\quad\text{with}\\
    b_n=\langle\bb{w}\odot\bm{\delta}_{\bb{x}},\bb{G}_{\Ls}^{\star}\ast\bm{\phi}_n\rangle.
\end{multline}
Finally, \eqref{eq:148} and the direct sum decomposition in Theorem \ref{th:M} Item 4 show that the extreme points $[\bb{e}]$ of the centered unit ball in  $(\M_{\Ls}(\R,\R^D)/\mathcal{N}_{\Ls}, \norm{\Ls\{\cdot\}}_{\M,\Zz'})$ are exactly
\begin{equation}
[\bb{e}]=\bb{G}_{\Ls}\ast(\bb{w}\odot\bm{\delta}_{\bb{x}})-\sum_{n=1}^{N}b_n\bb{p}_n+\mathcal{N}_{\Ls}=\bb{G}_{\Ls}\ast(\bb{w}\odot\bm{\delta}_{\bb{x}})+\mathcal{N}_{\Ls}.
\end{equation}
This also shows that all $\Ls$-splines $\bb{f}$ are in $\M_{\Ls}(\R,\R^D)$. Moreover if $\bb{f}=\sum_{k=1}^K\bb{G}_{\Ls}\ast(\bb{a}_k\odot\bm{\delta}_{\bb{x}_k})+\bb{q}$ with $\bb{x}_k\neq\bb{x}_{k'}$ for all rows, then we calculate that
\begin{align}
\norm{\Ls\{\bb{f}\}}_{\M,\Zz'}&=
\norm{\Ls\left\{\bb{G}_{\Ls}\ast\left(\sum_{k=1}^K\bb{a}_k\odot\bm{\delta}_{\bb{x}_k}\right)\right\}+\Ls\{\bb{q}\}}_{\M,\Zz'}=
\norm{\sum_{k=1}^K\bb{a}_k\odot\bb{\delta}_{\bb{x}_k}}_{\M,\Zz'}\\
&=\norm{\left(\sum_{k=1}^K\vert [\bb{a}_k]_d\vert\right)_{d=1}^D}_{\Zz'}.
\end{align}
\hfill\\\hfill\\
\textbf{Item 5.}
The argument is the same as the one in the Item 4, upon substitution of $x$ for $\bb{x}$. The norm calculation follows from Item 3 of Proposition \ref{prop:3}.
\end{proof}

\begin{proof}[\textbf{Proof of Proposition \ref{prop:stabilitymat}}]
\hfill\\
\textbf{Item 1.} The calculation of the Green's matrix is straightforward. The equality of native spaces follows from the observation that
\begin{equation}
\tilde{\Ls}\{\bb{f}\}\in\M(\R,\R^D)\Leftrightarrow\\(\bb{S}\Ls)\{\bb{f}\}\in\M(\R,\R^D)\Leftrightarrow\Ls\{\bb{f}\}\in\M(\R,\R^D).
\end{equation}
The equality of null spaces follows likewise.\hfill\\
\textbf{Item 2.} The calculation of the Green's matrix is straightforward. Next, we fix
\begin{equation}
    \bb{A}=\begin{bmatrix}
        0&1\\
        1&0
    \end{bmatrix},\quad
    \bb{L}=\begin{bmatrix}
        \D&0\\
        0&\D^2
    \end{bmatrix}\Rightarrow\quad\bb{G}_{\tilde{\Ls}}(\cdot)=\begin{bmatrix}
        0&(\cdot)_+\\
        (\cdot)_+^0&0
    \end{bmatrix}.
\end{equation}
If $\M_{\Ls}(\R,\R^D)=\M_{\tilde{\Ls}}(\R,\R^D)$, then 
\begin{equation}
    \bb{G}_{\tilde{\Ls}}\ast\begin{bmatrix}
        \delta_0\\
        0
    \end{bmatrix}\in\M_{\Ls}(\R,\R^D)\Rightarrow
    \begin{bmatrix}
        0\\
        (\cdot)_+^0
    \end{bmatrix}\in\M_{\Ls}(\R,\R^D)\Rightarrow
    \begin{bmatrix}
        0\\
        \delta_0'   
    \end{bmatrix}\in\M(\R,\R^D),
\end{equation}
which is a contradiction. The inequality of null spaces is proved likewise.
\end{proof}

\subsection{Proofs of Section 3.4}
\label{app:3.4}
Using the kernel $\bb{G}_{\bm{\phi}}$, we restate Theorem \ref{th:predual} as the precise Theorem \ref{th:Mpred}. The classic extension Theorem \ref{th:extension} taken from \cite{reed1972methods}
[Theorem I.7, p. 9] is necessary for the proof of Theorem \ref{th:Mpred}.

\begin{theorem}
\label{th:extension}
Let $G$ be a bounded linear transformation from a normed space $(\mathcal{X}, \norm{\cdot}_{\mathcal{X}})$ to a complete normed space $(\mathcal{Y}, \norm{\cdot}_{\mathcal{Y}})$. Then, G
has a unique extension to a bounded linear transformation (with the same norm) from the 
completion of $\mathcal{X}$ to $\mathcal{Y}$.
\end{theorem}

\begin{theorem} 
\label{th:Mpred}
Let $(\bb{p}, \bm{\phi})$ be an admissible system in the sense of Definition \ref{def:admissible}. Then, the function space
\begin{equation}
\label{eq:3.2.40}
    \C_{\Ls}(\R,\R^D)=\left\{\bb{g}=\Ls^{\star}\{\bb{v}\}+\sum_{n=1}^Na_n\bm{\phi}_n:\bb{v}\in\C_{0}(\R,\R^D),\bb{a}=(a_n)_{n=1}^N\in\R^N\right\}
\end{equation}
has the following properties
\begin{itemize}
    \item [1.] every $\bb{g}$ has a unique direct-sum representation as in \eqref{eq:3.2.40}, with $\bb{v}=\bb{G}_{\bm{\phi}}^{\star}\{\bb{g}\}$ and, $\forall n\in[1,N],a_n=\langle\bb{p}_n,\bb{g}\rangle$;
    \item [2.] the space $\C_{\Ls}(\R,\R^D)$ is a Banach space for the norm
    \begin{equation}
        \norm{\bb{g}}_{\C_{\Ls}}=\mathrm{max}(\norm{\bb{G}_{\bm{\phi}}^{\star}\{\bb{g}\}}_{\Zz,\infty},\norm{(\langle \bb{p}_n,\bb{g}\rangle)_{n=1}^N)}_{\Yy})=\mathrm{max}(\norm{\bb{v}}_{\Zz,\infty},\norm{\bb{a}}_{\Yy})
    \end{equation}
    or
        \begin{equation}
        \norm{\bb{g}}_{\C_{\Ls}}=\mathrm{max}(\norm{\bb{G}_{\bm{\phi}}^{\star}\{\bb{g}\}}_{\infty,\Zz},\norm{(\langle \bb{p}_n,\bb{g}\rangle)_{n=1}^N)}_{\Yy})=\mathrm{max}(\norm{\bb{v}}_{\infty,\Zz},\norm{\bb{a}}_{\Yy});
    \end{equation}
    \item [3.] the space $\C_{\Ls}(\R,\R^D)$ is the predual of $\M_{\Ls}(\R,\R^D)$; so that $\C_{\Ls}(\R,\R^D)'=\M_{\Ls}(\R,\R^D)$ with 
    \begin{equation}
        \norm{\Ls\{\cdot\}}_{\Zz',\M}+\norm{\left(\langle \cdot,\bm{\phi}_n\rangle\right)_{n=1}^N}_{\Yy'}\text{ being the dual norm of } \mathrm{max}(\norm{\bb{G}_{\bm{\phi}}^{\star}\{\cdot\}}_{\Zz,\infty},\norm{(\langle \bb{p}_n,\cdot\rangle)_{n=1}^N)}_{\Yy})
    \end{equation}
    and 
    \begin{equation}
        \norm{\Ls\{\cdot\}}_{\M,\Zz'}+\norm{\left(\langle \cdot,\bm{\phi}_n\rangle\right)_{n=1}^N}_{\Yy'}\text{ being the dual norm of } \mathrm{max}(\norm{\bb{G}_{\bm{\phi}}^{\star}\{\cdot\}}_{\infty,\Zz},\norm{(\langle \bb{p}_n,\cdot\rangle)_{n=1}^N)}_{\Yy}.
    \end{equation}
\end{itemize} 
\end{theorem}

\begin{proof}[\textbf{Proof of Theorem \ref{th:Mpred}}]
\hfill\\
\textbf{Item 1.}
We first calculate that $\forall\bm{\psi}\in\Dd(\R,\R^D),$ it holds that
\begin{equation}
    \bb{G}_{\bm{\phi}}\{\bm{\psi}\}(\cdot)=\left(\bb{G}_{\Ls}\ast\bm{\psi}\right)(\cdot)-\sum_{n=1}^N\bb{p}_n(\cdot)\left\langle\bb{G}_{\Ls}^{\star}\ast\bm{\phi}_n,\bm{\psi}\right\rangle=\left[\bb{G}_{\Ls}\ast\bm{\psi}\right](\cdot)-\sum_{n=1}^N\bb{p}_n(\cdot)\left\langle\bm{\phi}_n,\bb{G}_{\Ls}\ast\bm{\psi}\right\rangle.\label{eq:A.4.136}
\end{equation}
It follows from \eqref{eq:A.4.136} that, $\forall\bb{v}\in\C_0(\R,\R^D)$ and $\forall(a_n)_{n=1}^N\in\R^N$,
\begin{align}   &\left\langle\bb{G}_{\bm{\phi}}^{\star}\left\{\Ls^{\star}\{\bb{v}\}\right\},\bm{\psi}\right\rangle=\left\langle\bb{v},\Ls\{\bb{G}_{\bm{\phi}}\{\bm{\psi}\}\}\right\rangle=\left\langle\bb{v},\bm{\psi}\right\rangle\quad&&\Rightarrow\quad\bb{G}_{\bm{\phi}}^{\star}\left\{\Ls^{\star}\{\bb{v}\}\right\}=\bb{v},\\
&\left\langle\bb{G}_{\bm{\phi}}^{\star}\left\{\sum_{n=1}^Na_n\bm{\phi}_n\right\},\bm{\psi}\right\rangle= \left\langle\sum_{n=1}^Na_n\bm{\phi}_n,\bb{G}_{\bm{\phi}}\{\bm{\psi}\}\right\rangle=0\quad&&\Rightarrow\quad\bb{G}_{\bm{\phi}}^{\star}\left\{\sum_{n=1}^Na_n\bm{\phi}_n\right\}=\bm{0}\\
& &&\Rightarrow\bb{G}_{\bm{\phi}}^{\star}\left\{\bb{g}\right\}=\bb{v}.\label{eq:A.4.139}
\end{align}
We also calculate that, $\forall\bb{v}\in\Dd(\R,\R^D)$ and $\forall(a_n)_{n=1}^N\in\R^N$,
\begin{align}
\label{eq:A.4.141}
    \langle \bb{p}_n,\Ls^{\star}\{\bb{v}\}\rangle=\langle \Ls\{\bb{p}_n\},\bb{v}\rangle=0,\quad\text{and}\quad\left\langle\bb{p}_{n},\sum_{n'=1}^Na_{n'}\bm{\phi}_{n'}\right\rangle=a_{n}\quad\Rightarrow\langle\bb{p}_{n},\bb{g}\rangle=a_{n}.
\end{align}
Then, from Theorem \ref{th:extension} and from the density 
 of $\mathcal{D}(\R,\R^D)$ in $\C_0(\R,\R^D)$, we conclude that \eqref{eq:A.4.141} holds $\forall\bb{v}\in\C_0(\R,\R^D)$. The desired direct-sum decomposition follows directly.
\hfill\\\hfill\\
\textbf{Item 2.} Item 2 follows directly from Item 1 and \eqref{eq:A.4.139}, \eqref{eq:A.4.141}.
\hfill\\\hfill\\
\textbf{Item 3.} Item 3 follows directly from the direct-sum decompositions in Item 4 of Theorem \ref{th:M} and Item 1 of Theorem \ref{th:Mpred}, and from the two duality results 
\begin{align}
    &\left(\bb{G}_{\bm{\phi}}\M(\R,\R^D),\quad\norm{\Ls\{\cdot\}}_{\M,\Zz'}\right)=\left(\Ls^{\star}\C_0(\R,\R^D),\quad\norm{\bb{G}_{\bm{\phi}}^{\star}\{\cdot\}}_{\infty,\Zz}\right)'\\
    \text{and}&\left(\mathcal{N}_{\Ls},\norm{(\langle \cdot,\bm{\phi}_n\rangle)_{n=1}^N}_{\Yy'}\right)=\left(\text{span}(\bm{\phi}_n)_{n=1}^N,\quad\norm{(\langle \bb{p}_n,\cdot\rangle)_{n=1}^N}_{\Yy}\right)'.
\end{align}
The argumentation in the case of inner norms is the same. 
\end{proof}

\begin{proof}[\textbf{Proof of Corollary \ref{coro:diracpredual}}]We first calculate that 
\begin{equation}
    \bb{G}_{\bm{\phi}}^{\star}\{\bm{\rho}\}=\bb{G}^{\star}_{\Ls}\ast\bm{\rho}-\sum_{n=1}^N\langle\bm{\rho},\bm{p}_n\rangle\bb{G}_{\Ls}^{\star}\ast\bm{\phi}_n,
\end{equation}
wherefrom it follows that $\bb{G}_{\bm{\phi}}^{\star}\{\bm{\rho}\}$, as the convolution between distributions and test functions, is itself a continuous function. In addition, and following the argument made in the proof of Proposition \ref{prop:defmap}, we find that $\bb{G}_{\bm{\phi}}^{\star}\{\bm{\rho}\}$ is compactly supported. Consequently, $\bb{G}_{\bm{\phi}}^{\star}\{\bm{\rho}\}\in\C_0(\R,\R^D)$. Then, we claim that $\Ls^{\star}\{\bb{G}_{\bm{\phi}}^{\star}\{\bm{\rho}\}\}+\sum_{n=1}^N\langle\bm{\rho},\bb{p}_n\rangle\bm{\phi}_n=\bm{\rho}$, which concludes the proof. We calculate that
\begin{align}
    \Ls^{\star}\{\bb{G}_{\bm{\phi}}^{\star}\{\bm{\rho}\}\}+\sum_{n=1}^N\langle\bm{\rho},\bb{p}_n\rangle\bm{\phi}_n&=\Ls^{\star}\left\{\bb{G}_{\Ls}^{\star}\ast\bm{\rho}-\sum_{n=1}^N\langle\bm{\rho},\bb{p}_n\rangle\bb{G}_{\Ls}^{\star}\ast\bm{\phi}_n\right\}+\sum_{n=1}^N\langle\bm{\rho},\bb{p}_n\rangle\bm{\phi}_n=\bm{\rho}.
\end{align}
\end{proof}

\subsection{Proofs of Section 4.2}
\label{app:4.3}
\begin{proof}[\textbf{Proof of Theorem \ref{th:optweakstar}}]
\hfill\\
    The proofs for inner and outer norms are the same, up to the description of the extreme points whose difference has already been explained in Items 4 and 5 of Theorem \ref{prop:MLprop}. Therefore, we only provide a proof for the case of the outer norm. For all topological matters, we equip the search space $\M_{\Ls}(\R,\R^D)$ with its $\text{weak}^{\star}$-topology, whose existence is guaranteed by Theorem \ref{th:predual}. 

    First observe that $\mathcal{J}$ is convex, $\mathrm{weak}^{\star}$ lower-semicontinuous, proper and coercive. Therefore, in vertue of Proposition 8 in \cite{gupta2018continuous}, $\mathcal{V}$ is non empty, $\text{weak}^{\star}$-compact, and the $\text{weak}^{\star}$-closed convex hull of its extreme points $\bb{e}$. In turn, it follows from Corollary 3.8 in \cite{boyer2019representer} that an extreme point $\bb{e}$ must be of the form
\begin{equation}
    \bb{e}=\bb{q}+\sum_{k=1}^{M-N}\alpha_k\bb{e}_k,\quad\bb{q}\in\mathcal{N}_{\Ls},
\end{equation}
where $\alpha_k\in\R$ and $\bb{e}_k$ are extreme points of the unit ball in $\bb{G}_{\bm{\phi}}\M(\R,\R^D)$. Moreover, we know from Theorem \ref{th:M} that $(\bb{G}_{\bm{\phi}}\M(\R,\R^D),\norm{\Ls\{\cdot\}}_{\M,\Zz'})$ is isometric to $(\M(\R,\R^D),\norm{\cdot}_{\M,\Zz'})$. Hence, the extreme points in $\bb{G}_{\bm{\phi}}\M(\R,\R^D)$ are all of the form 
\begin{equation}
\bb{G}_{\bm{\phi}}\{\bb{m}\},\quad\bb{m}=\bb{a}\odot\bm{\delta}_{\bb{x}}\quad\text{or equivalently,}\quad\bb{G}_{\Ls}\ast\bb{m}-\sum_{n=1}^N\bb{p}_n\langle \bb{m},\bb{G}_{\Ls}^{\star}\ast\bm{\phi}_n\rangle,
\end{equation}
where $\bb{x}\in\R$ and $\bb{a}$ is an extreme point of the unit ball in $\Zz'.$  It follows that 
\begin{equation}
    \bb{e}=\bb{G}_{\Ls}\ast\left(\sum_{k=1}^{M-N}\alpha_k\bb{a}_k\odot\bm{\delta}_{\bb{x}_k}\right)+\bb{q}-\sum_{n=1}^N\bb{p}_n\left\langle\sum_{k=1}^{M-N}\alpha_k\bb{a}_k\odot\bb{x}_k,\bb{G}_{\Ls}^{\star}\ast\bm{\phi}_n\right\rangle.
\end{equation}
\end{proof}

\subsection{Proofs of Section 4.3}
\label{app:4.4}

\begin{proof}[\textbf{Proof of Proposition \ref{prop:banachrestricted}}]We provide the proof only in the case of inner norms. For $\bb{P}=\bb{A}^{\intercal}\bb{A}$ and $\bb{m}\in\M(\R,\R^D)$, we establish the continuity bound
\begin{equation}
    \norm{(\bb{I}-\bb{P})\{\bb{T}\{\bb{m}\}\}}_{\Zz',\M}=\underset{\pi}{\mathrm{sup}}\sum_{A\in\pi}\norm{(\bb{I}-\bb{P})\{\bb{T}\{\bb{m}(A)\}\}}_{\Zz'}\leq\norm{(\bb{I}-\bb{P})\bb{T}}\norm{\bb{m}}_{\Zz',\M},
\end{equation}
where $\norm{(\bb{I}-\bb{P})\bb{T}}=\underset{\bb{z}:\norm{\bb{z}}_{\Zz'}\leq1}{\mathrm{max}}\norm{(\bb{I}-\bb{P})\{\bb{T}\{\bb{z}\}\}}_{\Zz'}$. Let $(\bb{f}_n)_{n=1}^{\infty}$ be a sequence that is convergent in the strong topology to some $\bb{f}\in\M_{\Ls}(\R,\R^D)$, with $\bb{f}_n\in\M_{\tilde{\Ls}\vert\bb{A}}(\R,\R^D)$. It follows that $\bb{f}\in\M_{\tilde{\Ls}\vert\bb{A}}(\R,\R^D)$ because
\begin{align}
    \norm{(\bb{I}-\bb{P})\{\tilde{\Ls}\{\bb{f}\}\}}_{\Zz',\M}&\leq\norm{(\bb{I}-\bb{P})\{\tilde{\Ls}\{\bb{f}-\bb{f}_n\}\}}_{\Zz',\M}+\norm{(\bb{I}-\bb{P})\{\tilde{\Ls}\{\bb{f}_n\}\}}_{\Zz',\M}\\
    &\leq\norm{(\bb{I}-\bb{P})\bb{T}}
    \norm{\Ls\{\bb{f}-\bb{f}_n\}}_{\Zz',\M}\underset{n\to\infty}{\to}0.
\end{align}
\end{proof}

\begin{proof}[\textbf{Proof of Corollary \ref{coro:restrictedRT1}}]
We prove the extreme-point representation in \eqref{eq:4.3.76} for inner norms $\norm{\cdot}_{\Ww',\M}$. First observe that, $\forall\bb{f}\in\M_{\tilde{\Ls}\vert\bb{A}}(\R,\R^D)$,
\begin{equation}
\label{eq:4.3.79}
    \norm{(\bb{Q}\Ls)\{\bb{f}\}}_{\Ww',\M}=\norm{\tilde{\bb{S}}\{\tilde{\Ls}\{\bb{f}\}\}}_{\Zz',\M},\quad\text{with}\quad
    \tilde{\bb{S}}=\begin{bmatrix}
        \bb{S}&\bm{0}\\
        \bm{0}&\hat{\bb{S}}
    \end{bmatrix},
\end{equation}
where $\hat{\bb{S}}\in\R^{(D-D')\times(D-D')}$ is any invertible matrix and where $\norm{\tilde{\bb{S}}\{\cdot\}}_{\Zz',\M}$ is an inner norm over $\M(\R,\R^D)$.
Since $\M_{\tilde{\Ls}\vert\bb{A}}(\R,\R^D)$ is a Banach subspace of $\M_{\Ls}(\R,\R^D)=\M_{\tilde{\Ls}}(\R,\R^D)$ equipped with the inner semi-norm $\norm{(\tilde{\bb{S}}\tilde{\Ls})\{\cdot\}}_{\Zz',\M}$, it follows from the same argument as in the proof of Theorem \ref{th:optweakstar} that the extreme points  of $\mathcal{V}$ are of the form 
\begin{equation}
\label{eq:4.3.81}
\bb{p}+\sum_{k=1}^Ka_k\bb{f}_k,\quad\quad a_k\in\R,\quad K\leq(M-N),\quad\bb{p}\in\mathcal{N}_{\Ls},
\end{equation}
where $\bb{f}_k$ is a representative of an extreme point of the unit ball in $\M_{\tilde{\Ls}\vert\bb{A}}(\R,\R^D)/\mathcal{N}_{\Ls}$. 
The extreme points of the centered unit ball in 
\begin{equation}
    (\M(\R,\R^{D'}), \norm{\bb{S}\{\cdot\}}_{\Ww',\M})\subset(\M(\R,\R^{D}), \norm{\tilde{\bb{S}}\{\cdot\}}_{\Zz',\M})
\end{equation}
are exactly of the form
\begin{equation}
\label{eq:4.3.82}
    \bb{S}^{-1}\{\bb{a}\delta_x\}=\bb{S}^{-1}\{\bb{a}\}\delta_x,\quad\text{with}\quad x\in\R,\quad\bb{a}\in\R^{D'},
\end{equation}
 where $\bb{a}$ is an extreme point of the centered unit ball in $(\R^{D'},\norm{\cdot}_{\Ww'}).$ In addition, 
\begin{equation}
    \bb{A}\tilde{\Ls}:\M_{\tilde{\Ls}\vert\Ww'}(\R,\R^D)/\mathcal{N}_{\Ls}\to\M(\R,\R^{D'})
\end{equation}
is an isometry. The inverse of an element of the form \eqref{eq:4.3.82} is thus the equivalence class represented by 
\begin{equation}
\label{eq:4.3.84}
(\bb{G}_{\tilde{\Ls}}\bb{A}^{\intercal})\{\bb{S}^{-1}\{\bb{a}\delta_x\}\}=(\bb{G}_{\Ls}\bb{T}^{-1}\bb{A}^{\intercal}\bb{S}^{-1})\{\bb{a}\delta_x\}=(\bb{G}_{\Ls}\bb{Q}^{\dagger})\{\bb{a}\delta_x\}
\end{equation}
which, therefore, exactly represents the extreme points of the unit ball in $\M_{\tilde{\Ls}\vert\bb{A}'}(\R,\R^D)/\mathcal{N}_{\Ls}$. The desired representation follows from \eqref{eq:4.3.81} and \eqref{eq:4.3.84}.
Finally, the proof of the extreme-point representation in \eqref{eq:4.3.76} for outer norms $\norm{\cdot}_{\M,\Ww'}$ has a unique  feature. It is that the norm $\norm{\tilde{\bb{S}}\{\cdot\}}_{\M,\Zz'}$ in \eqref{eq:4.3.79}, now defined as 
\begin{equation}
\label{eq:206}
    \norm{\tilde{\bb{S}}\{\bb{m}\}}_{\M,\Zz'}:=\norm{\underset{\pi}{\text{sup}}\sum_{A\in\pi}\vert\tilde{\bb{S}}\{\bb{m}\}(A)\vert}_{\Zz'},
\end{equation}
 is \emph{not} an outer norm anymore. It is not an inner norm either. Nevertheless, a short calculation would show that the norm in \eqref{eq:206} is $\mathrm{weak}^{\star}$-continuous and equivalent to $\norm{\cdot}_{1,\M}=\norm{\cdot}_{\M,1}.$
These properties are required to ensure that $\mathcal{V}$ is non empty. Finally, the description of extreme points is proved as for inner norms, with an isometry-based argumentation.
\end{proof}

\begin{proof}[\textbf{Proof of Corollary \ref{coro:restrictedRT2}}]
Let $\bb{f}\in\M_{\tilde{\Ls}\vert\bb{A}}(\R,\R^D)$. It follows from
\begin{align}
    \Ls\{\bb{f}\}=(\bb{T}^{-1}\bb{A}^{\intercal}\bb{S}^{-1})\{\bb{m}\}\Leftrightarrow\tilde{\Ls}\{\bb{f}\}=(\bb{A}^{\intercal}\bb{S}^{-1})\{\bb{m}\}&\Leftrightarrow(\bb{A}\tilde{\Ls})\{\bb{f}\}=\bb{S}^{-1}\{\bb{m}\}\\
    &\Leftrightarrow(\bb{Q}\Ls)\{\bb{f}\}=\bb{m}
\end{align}
that $(\bb{f},\bb{m})\in\M_{\tilde{\Ls}\vert\bb{A}}(\R,\R^{D})\times\M(\R,\R^{D'})$ verifies \eqref{eq:4.3.69} if and only if $\bb{m}=(\bb{Q}\circ\Ls)\{\bb{f}\}$. The equality of the solution sets is thus verified.
\end{proof}

\section{Convolution of Distributions}
\label{app:B}
\subsection{Convolution of Scalar Distributions}
The convolution $g\ast\psi\in\Dd(\R)$ of a scalar distribution $g\in\Dd'(\R)$ and a test function $\psi\in\Dd(\R)$ is defined as, 
\begin{equation}
\forall x\in\R:\quad\quad (g\ast\psi)(x)=\langle g^{\vee},\psi(\cdot+x)\rangle:=\langle g^{\vee}\ast\delta_{x},\psi\rangle.
\end{equation}
The convolution $g_1\ast g_2$ of two scalar distribution $(g_1,g_2)\in\mathcal{D}'(\R)^2$ is defined as \cite{beffa2024weakly}[Section 3.2]
\begin{equation}
    \forall\psi\in\mathcal{D}(\R):\quad\langle g_1\ast g_2,\psi\rangle=\langle g_1((\cdot)_1)\otimes g_2((\cdot)_2),\psi((\cdot)_1+(\cdot)_2)\rangle:=\int_{\R^2}g_1(x)g_2(y)\psi(x+y)\mathrm{d}x\mathrm{d}y,
\end{equation}
where the second equality is the informal statement of the central one. Alternatively, one can rely on Fubini's theorem for distributions \cite{treves2016topological} to establish the tensor-product-based definition
\begin{equation}
    \langle g_1\ast g_2,\psi\rangle=\langle g_2,g_1^{\vee}\ast\psi\rangle=\langle g_1,g_2^{\vee}\ast\psi\rangle.
\end{equation}
We summarize in Proposition \ref{prop:9} three properties of the convolution, all extracted from \cite{beffa2024weakly}[Section 3.2].
\begin{ Proposition }
\label{prop:9}
    Let $(g_1,g_2,g_3)\in\Dd'(\R)^3$ be three scalar distributions, either all compactly supported or all one-sided, and let $\ls$ be an ODO. Then the following holds:
    \begin{itemize}
        \item [] \emph{Commutativity} $\forall\psi\in\Dd(\R):\quad\langle g_1\ast g_2,\psi\rangle=\langle g_2\ast g_1,\psi\rangle$.
        \item [] \emph{Associativity} $\forall\psi\in\Dd(\R):\quad\langle (g_1\ast g_2)\ast g_3,\psi\rangle=\langle g_1\ast( g_2\ast g_3),\psi\rangle$.
        \item [] \emph{Differentiation}$\forall\psi\in\Dd(\R):\quad\langle \ls\{g_1\ast g_2\},\psi\rangle=\langle \ls\{g_1\}\ast g_2,\psi\rangle=\langle g_1\ast\ls\{g_2\},\psi\rangle$.
    \end{itemize}
\end{ Proposition }

\subsection{Convolution of Matrix-Valued Distributions}
First, the action of a matrix-valued distribution $\bb{G}\in\Dd'(\R,\R^{D_1\times D_2})$ on a matrix-valued test function $\bm{\Psi}\in\Dd(\R,\R^{D_1\times D})$ is such that 
\begin{multline}
    \langle\bb{G},\bm{\Psi}\rangle\in\R^{D_2\times D },\quad\text{and}\quad\forall1\leq d_2\leq D_2,\forall1\leq d\leq D:\\
    \left[\langle\bb{G},\bm{\Psi}\rangle\right]_{d_2,d_1}=\sum_{d_1}^{D_1}\langle\bb{G}^{\intercal}_{d_2,d_1},\Psi_{d_1,d}\rangle.
\end{multline}
The convolution of a matrix-valued distribution $\bb{G}\in\Dd'(\R,\R^{D_2\times D_1})$ with a matrix-valued test function $\bm{\Psi}\in\Dd(\R,\R^{D_1\times D})$ is defined as 
\begin{multline}
\bb{G}\ast\bb{\Psi}\in\Dd(\R,\R^{D_2\times D})\quad\text{with}\quad
    \forall1\leq d_2\leq D_2,1\leq d\leq D:\\
    [\bb{G}\ast\bm{\Psi}]_{d_2,d}=\sum_{d_1=1}^{D_1}\mathrm{G}_{d_2,d_1}\ast\Psi_{d_1,d}.
\end{multline}
In turn, the convolution of two matrix-valued distributions  $\bb{G}\in\Dd'(\R,\R^{D_1\times D_2})$ and $\bb{H}\in\Dd'(\R,\R^{D_2\times D_3})$ is defined as  
\begin{multline}
\bb{G}\ast\bb{H}\in\Dd'(\R,\R^{D_1\times D_3})\quad\text{with}\quad\forall1\leq d_1\leq D_1,1\leq d_3\leq D_3:\\
[\bb{G}\ast\bb{H}]_{d_1,d_3}=\sum_{d_2=1}^{D_2}\mathrm{G}_{d_1,d_2}\ast\mathrm{H}_{d_2,d_3}.
\end{multline}
We summarize in Proposition \ref{prop:10} four important properties of the convolution of matrix-valued distributions. 
\begin{ Proposition }
\label{prop:10}
Let $\bb{G}\in\Dd'(\R,\R^{D_1\times D_2})$, and $\bb{H}\in\Dd'(\R,\R^{D_2\times D_3})$, and $\bb{J}\in\Dd'(\R,\R^{D_3\times D_4})$ be three matrix-valued distributions either all compactly supported or all one-sided, and let $\Ls$ be an MDO of size $(D_1\times D_1)$. Then, the following holds $\forall\bb{\Psi}\in\Dd(\R,\R^{D_1\times D})$.
\begin{itemize}
    \item[] \emph{Adjoint}$\quad\left\langle\bb{G}\ast\bb{H},\bb{\Psi}\right\rangle=\left\langle\bb{H},\bb{G}^{\intercal,\vee}\ast\bb{\Psi}\right\rangle\quad\Rightarrow\quad\bb{G}^{\star}=\bb{G}^{\intercal,\vee}.$
    \item [] \emph{Associativity} $\quad\left\langle \bb{G}\ast\left(\bb{H}\ast\bb{J}\right),\bb{\Psi}\right\rangle=\left\langle \left(\bb{G}\ast\bb{H}\right)\ast\bb{J},\bb{\Psi}\right\rangle.$
    \item [] \emph{Differentiation} $\quad\left\langle \Ls\{\bb{G}\ast\bb{H}\},\bb{\Psi}\right\rangle=\left\langle \Ls\{\bb{G}\}\ast\bb{H},\bb{\Psi}\right\rangle=\left\langle \bb{G}\ast\bb{H},\Ls^{\star}\{\bb{\Psi}\}\right\rangle.$
    \item [] \emph{Bracket} $\left(\bb{G}^{\intercal}\ast\bm{\Psi}\right)(x)=\left\langle\bb{G}^{\vee}\ast\bb{I}\delta_{x},\bm{\Psi}\right\rangle.$
    
\end{itemize}
\end{ Proposition }
\begin{proof}[\textbf{Proof of Proposition \ref{prop:10}}]
\hfill\\
    \textbf{Adjoint.}
    We calculate that, $\forall1\leq d_3\leq D_3,\forall1\leq d\leq D,$
    \begin{align}
        \left[\left\langle\bb{G}\ast\bb{H},\bm{\Psi}\right\rangle \right]_{d_3,d}&=\sum_{d_1=1}^D\left\langle\left(\bb{G}\ast\bb{H}\right)^{\intercal}_{d_3,d_1},\Psi_{d_1,d}\right\rangle=\sum_{d_1,d_2=1,1}^{D_1,D_2}\left\langle \mathrm{H}_{d_2,d_3}\ast\mathrm{G}_{d_1,d_2},\Psi_{d_1,d}\right\rangle\\
        &=\sum_{d_1,d_2=1,1}^{D_1,D_2}\left\langle \mathrm{H}_{d_2,d_3},\mathrm{G}_{d_1,d_2}^{\vee}\ast\Psi_{d_1,d}\right\rangle=\sum_{d_2=1}^{D_2}\left\langle \mathrm{H}_{d_2,d_3},[\bb{G}^{\intercal,\vee}\ast\bm{\Psi}]_{d_2,d}\right\rangle\\
        &=\left\langle\bb{H},\bb{G}^{\intercal,\vee}\ast\bb{\Psi}\right\rangle.
    \end{align}
    \hfill\\
    \textbf{Associativity.} We calculate that, $\forall1\leq d_4\leq D_4$ and $\forall1\leq d\leq D,$
    \begin{align}
        \left[\langle \bb{G}\ast\left(\bb{H}\ast\bb{J}\right),\bb{\Psi}\rangle\right]_{d_4,d}&=\sum_{d_1=1}^{D_1}\langle[(\bb{H}\ast\bb{J})^{\intercal}\ast\bb{G}^{\intercal}]_{d_3,d_1},\Psi_{d_1,d}\rangle\\
        &=\sum_{d_1=1}^{D_1}\langle[(\bb{J}^{\intercal}\ast\bb{H}^{\intercal})\ast\bb{G}^{\intercal}]_{d_3,d_1},\Psi_{d_1,d}\rangle\label{eq:B.2.162}\\
        &=\sum_{d_1=1}^{D_1}\langle[\bb{J}^{\intercal}\ast(\bb{H}^{\intercal}\ast\bb{G}^{\intercal})]_{d_3,d_1},\Psi_{d_1,d}\rangle\\
        &=\sum_{d_1=1}^{D_1}\langle[\bb{J}^{\intercal}\ast(\bb{G}\ast\bb{H})^{\intercal}]_{d_3,d_1},\Psi_{d_1,d}\rangle\label{eq:B.2.163}\\
        &=\left[\langle (\bb{G}\ast\bb{H})\ast\bb{J},\bb{\Psi}\rangle\right]_{d_4,d},
    \end{align}
where the associativity of the convolution for scalar-valued distribution allowed us to deduce \eqref{eq:B.2.163} from \eqref{eq:B.2.162}.
\hfill\\\hfill\\
\textbf{Differentiation.} Observe that the action of the MDO $\Ls$ is exactly the same as the convolution with the distribution $\Ls\{\bb{I}\delta\}$. It follows that 
\begin{align}
    \left\langle \Ls\{\bb{G}\ast\bb{H}\},\bb{\Psi}\right\rangle&=\left\langle \Ls\{\bb{I}\delta\}\ast(\bb{G}\ast\bb{H}),\bb{\Psi}\right\rangle=\left\langle (\Ls\{\bb{I}\delta\}\ast\bb{G})\ast\bb{H},\bb{\Psi}\right\rangle=\left\langle \Ls\{\bb{G}\}\ast\bb{H},\bb{\Psi}\right\rangle.
\end{align}
We use the adjoint property to prove the second and third equalities of the Differentiation statement.
\hfill\\ \hfill\\
\textbf{Bracket.}
We  calculate that, $\forall1\leq d_2\leq D_2,\forall1\leq d,\leq D$,
   \begin{align}
       [\left(\bb{G}^{\intercal}\ast\bm{\Psi}\right)]_{d_2,d}(x)&=\sum_{d_1=1}^{D_1}(\mathrm{G}_{d_1,d_2}\ast\Psi_{d_1,d})(x)=\sum_{d_1=1}^{D_1}\left\langle\mathrm{G}_{d_ 1,d_2}^{\vee}\ast\delta_{x},\Psi_{d_1,d}\right\rangle\\
       &=\sum_{d_1=1}^{D_1}\left\langle\left[\bb{G}^{\vee}\ast\bb{I}\delta_{x}\right]_{d_1,d_2},\Psi_{d_1,d}\right\rangle=\left[\left\langle\left(\bb{G}^{\vee}\ast\bb{I}\delta_{x}\right),\bb{\Psi}\right\rangle\right]_{d_2,d}\\
       &=\left[\left\langle\bb{G}^{\vee}\ast\bb{I}\delta_{x},\bb{\Psi}\right\rangle\right]_{d_2,d}.
   \end{align}
\end{proof}
\end{document}